\documentstyle[12pt]{article}
\newtheorem{thm}{Theorem}[section]
\newtheorem{prop}[thm]{Proposition}
\newtheorem{claim}[thm]{Claim}

\newtheorem{thm-defi}[thm]{Theorem/Definition}
\newtheorem{example}[thm]{Example}
\newtheorem{cor}[thm]{Corollary}

\newtheorem{new-lemma}[thm]{Lemma}
\newtheorem{defi}[thm]{Definition}
\newtheorem{rem}[thm]{Remark}

\newtheorem{condition}[thm]{Condition}

\newcommand{\B}{{\cal B}}

\newcommand{\E}{{\cal E}}
\newcommand{\F}{{\cal F}}

\renewcommand{\H}{{\cal H}}

\newcommand{\LB}{{\cal L}}

\newcommand{\M}{{\cal M}}
\newcommand{\N}{{\cal N}}

\newcommand{\W}{{\cal W}}
\newcommand{\X}{{\cal X}}
\newcommand{\Y}{{\cal Y}}
\newcommand{\Z}{{\cal Z}}

\newcommand{\one}{{1}}
\newcommand{\PP}{{\Bbb P}}

\newcommand{\Integers}{{\Bbb Z}}
\newcommand{\ComplexNumbers}{{\Bbb C}}
\newcommand{\RationalNumbers}{{\Bbb Q}}
\newcommand{\LieAlg}[1]{{\frak #1}}
\newcommand{\linsys}[1]{{\mid}#1{\mid}}

\newcommand{\IsomRightArrow}{\stackrel{\cong}{\rightarrow}}
\newcommand{\LongIsomRightArrow}{\stackrel{\cong}{\longrightarrow}}
\newcommand{\RightArrowOf}[1]{\stackrel{#1}{\rightarrow}}

\newcommand{\LongLeftArrowOf}[1]{\stackrel{#1}{\longleftarrow}}

\newcommand{\LongRightArrowOf}[1]{\stackrel{#1}{\longrightarrow}}
\newcommand{\LongIsomRightArrowOf}[1]{
\stackrel
{\stackrel{#1}{\cong}}
{\longrightarrow}
}

\newcommand{\StructureSheaf}[1]{{\cal O}_{#1}}
\newcommand{\EndProof}{\hfill  $\Box$}
\newcommand{\restricted}[2]{#1_{\mid_{#2}}}

\newcommand{\rank}{{\rm rank}}
\newcommand{\coker}{{\rm coker}}
\newcommand{\Pic}{{\rm Pic}}

\newcommand{\Hom}{{\rm Hom}}
\newcommand{\Aut}{{\rm Aut}}
\newcommand{\End}{{\rm End}}

\newcommand{\Abs}[1]{\mid\!#1\!\mid}

\newcommand{\SheafHom}{{\cal H}om}

\newcommand{\SheafExt}{{\cal E}xt}
\newcommand{\SheafTor}{{\cal T}or}
\newcommand{\RelExt}{{\cal E}xt}

\newcommand{\Ideal}[1]{{\cal I}_{#1}}

\newcommand{\Wedge}[1]{\stackrel{#1}{\wedge}}

\input amssymb.sty

\begin{document}
\begin{center}
\begin{Large}
{\bf 
\noindent
On the monodromy of moduli spaces of sheaves on K3 surfaces II
}
\end{Large}
\\
Eyal Markman
\footnote{Partially supported by NSF grant number DMS-9802532}
\end{center}

{\scriptsize 
\tableofcontents 
} 

\section{Introduction}
\label{sec-introduction}

Let $S$ be a projective K3 surface. 
Moduli spaces of stable coherent sheaves on
$S$ are parametrized by their Mukai vectors in the Mukai lattice. 
This lattice is the cohomology $H^*(S,\Integers)$, endowed with 
a modification of the Poincare pairing (see equation 
(\ref{eq-mukai-pairing}) below). Denote by $S^{[n]}$ the Hilbert scheme 
of length $n$ subschemes of $S$. 
Set $v$ to be the Mukai vector of 
the ideal sheaf of $n$ points on $S$.
Then the Hilbert scheme $S^{[n]}$ is isomorphic to the moduli space $\M(v)$, 
of stable sheaves with Mukai vector $v$. 

We denote the isometry group of the Mukai lattice by $\Gamma$ and
let $\Gamma_v$ be the stabilizing subgroup of the Mukai vector $v$.
We began, in part I of the paper, the construction
of a representation

\begin{equation}
\label{eq-action-of-stabilizer}
\gamma \ : \ \Gamma_v \ \ \longrightarrow \ \ 
\Aut(H^*(S^{[n]},\Integers)_{\rm free}),
\end{equation}
acting via ring automorphisms. Above,
$H^*(S^{[n]},\Integers)_{\rm free}$ denotes the torsion free summand of the 
integral cohomology.
We reduced the construction to two classes of isometries
$\tau_{v_0}$ and $\sigma_{u_0}$ in $\Gamma_v$, reviewed below. 

The main result of the paper is Theorem \ref{thm-Gamma-v-acts-motivicly}.
It relates the representation (\ref{eq-action-of-stabilizer})
to the monodromy representation of the Hilbert scheme. 
The Mukai lattice is endowed with a weight $2$ Hodge structure. 
Isometries preserving this structure are called {\em Hodge isometries}. 
Theorem \ref{thm-Gamma-v-acts-motivicly} relates also the 
automorphism $\gamma(g)$, for any Hodge isometry  $g\in \Gamma_v$,  
to an autoequivalence of the bounded derived category $D(S)$ of
coherent sheaves on $S$. 

In part one of the paper, we reduced the proof of Theorem 
\ref{thm-Gamma-v-acts-motivicly} to 
the case of two sequences of isometries $\tau_n$, $\sigma_n$, 
$n\geq 2$. 
Recall, that every simple and rigid sheaf $L$ on $S$, is 
spherical, in the sense of Seidel-Thomas, 
and gives rise to an auto-equivalence of $D(S)$
(see \cite{seidel-thomas} or 
part one of the paper \cite{markman-part-one}). 
The corresponding isometry of the Mukai lattice, is the reflection
with respect to the $-2$ Mukai vector $v(L)$. 
The isometry $\tau_n$, in the first sequence, is the reflection
with respect to $v(L)$, where $L$ is a line bundle with $c_1(L)^2=2n-4$. 
When $L$ is such a line bundle,  then
$v(L)$ is perpendicular to the Mukai vector $v$ of $S^{[n]}$
and $\tau_n$ stabilizes $v$. 

The verification of the results of Theorem
\ref{thm-Gamma-v-acts-motivicly} for $\tau_n$, is carried out in section 
\ref{sec-self-dual-stratified-elemetary-transformations} (Theorem
\ref{thm-class-of-correspondence-in-stratified-elementary-trans}). 
Let us consider the case $n=1$, where the line bundle $L$ is
$\StructureSheaf{S}(\Sigma)$ for a smooth rational curve
$\Sigma\subset S$. The corresponding monodromy automorphism $\gamma(\tau_1)$ 
of $H^*(S^{[1]})$ is induced by the correspondence
\[
\Z_0 \ + \ [\Sigma \times \Sigma] \ \ \ \subset \ \ \ S^{[1]}\times S^{[1]},
\]
where $\Z_0$ is the diagonal. This is a classical reflection of a K3,
with respect to a $-2$ curve. 
Observe, that the stratification 
\[
S^{[1]} \ \supset \ \Sigma
\]
is a Brill-Noether stratification, with respect to the dimension of
$H^1(I_p\otimes L)$, where $I_p$ is the ideal sheaf of a point $p\in S$. 
This observation is key to the proof 
of the general case $n\geq 2$. 
Let us consider the case of a very ample line-bundle $L$,
giving rise to an embedding $S\hookrightarrow \PP^{n-1}$. 
The analogous Brill-Noether stratification is
\[
S^{[n]} = (S^{[n]})^0 \ \supset \ (S^{[n]})^1 \ \supset \ \cdots \ \supset \
(S^{[n]})^\mu.
\]
$(S^{[n]})^t$ is the locus, in which the $n$-points span a
projective linear subspace of dimension $\leq n-1-t$. 
This stratification is a special case of those considered in
\cite{markman-reflections}. The stratum $(S^{[n]})^t\setminus (S^{[n]})^{t+1}$
is a $G(t,2t)$-bundle over the open Brill-Noether stratum 
of a moduli space $M_t$ of sheaves of rank $t+1$. We set
\[
\Z_t \ \subset \ S^{[n]}\times S^{[n]}
\]
to be the closure of the fiber product 
of $(S^{[n]})^t\setminus (S^{[n]})^{t+1}$ with itself over $M_t$. 
$\Z_0$ is the diagonal.
Theorem \ref{thm-class-of-correspondence-in-stratified-elementary-trans} 
proves, in particular, that the cohomology class 
of $\sum_{t=0}^\mu \Z_t$ is equal to the class
$\gamma(\tau_n)$, with $\gamma$ given in
(\ref{eq-action-of-stabilizer}). 
The formula (\ref{eq-Chow-theoretic-formula-for-gamma-tau}) 
for $\gamma(\tau_n)$ expresses the class of 
$\sum_{t=0}^\mu \Z_t$ in terms of the Chern classes of the ideal sheaf 
of the universal subscheme in $S\times S^{[n]}$. 
We prove in section
\ref{sec-monodromy-of-A1-singularities}, that the endomorphism
of $H^*(S^{[n]},\Integers)$, induced by the correspondence 
$\sum_{t=0}^\mu \Z_t$, 
is a monodromy operator. 

The verification of the results of Theorem
\ref{thm-Gamma-v-acts-motivicly}, for 
the second sequence of isometries $\sigma_n$, 
is carried out in section \ref{sec-stratified-elementary-trans-plus-2-vectors}
(Theorem \ref{thm-reflection-sigma-satisfies-main-conj}).
The isometry $\sigma_n$ is a reflection, with respect to a $+2$ vector
$u_n$ orthogonal to $v$. We set $u_n$ to be the Mukai vector of
$I_{p+q}\otimes L^{-1}$, where 
$I_{p+q}$ is the ideal sheaf of a length $2$ subscheme, and 
$L$ is a line bundle with $c_1(L)^2=2n$. The 
isometry $\sigma_n$ corresponds to a {\em contravariant} 
equivalence $\Phi \ : \ D(S) \rightarrow D(S)^{op}$. 
It is the composition, of a covariant auto-equivalence, with the duality 
functor. $\Phi$ induces a birational involution of $S^{[n]}$,
which we now describe. 
When $n=1$, this is the regular Galois involution of a double
cover $S\rightarrow \PP^2$ branched along a sextic. When $n\geq 2$,
assume, for simplicity, that 
$L$ is very ample, giving rise to an embedding
$S\hookrightarrow \PP^{n+1}$. 
A generic length $n$ subscheme $D$ spans a $\PP^{n-1}$,
intersecting $S$ along a length $2n$ subscheme $Z$ containing $D$. 
The birational involution sends $D$ to the complementary 
subscheme $Z\setminus D$.

The line bundle $L$ induces a Brill-Noether stratification of
$S^{[n]}$, via the dimension of $H^1(I_D\otimes L)$. 
Again, we get a
reducible correspondence $\sum_{t=0}^\mu\Z_t$ in
$S^{[n]}\times S^{[n]}$. Here 
$\Z_0$ is the closure of the graph of the birational involution. 
Theorem \ref{thm-reflection-sigma-satisfies-main-conj} proves, in particular, 
that the cohomology class 
$\sum_{t=0}^\mu \Z_t$ is equal to the class of the composition 
$D_{S^{[n]}}\circ \gamma(\sigma_n)$, with $\gamma$ given in
(\ref{eq-action-of-stabilizer}) and $D_{S^{[n]}}$ is the duality operator,
multiplying $H^{2i}(S^{[n]})$ by $(-1)^i$.

\medskip
The paper is organized as follows. 
We state our main Theorem \ref{thm-Gamma-v-acts-motivicly}
in section \ref{sec-summary-of-part-one}. 
In section \ref{sec-stratified-elemetary-transformations}
we review the definition of stratified Mukai elementary transformations
of holomorphic symplectic varieties. These are birational transformations,
corresponding to stratification, 
such as the Brill-Noether stratifications of $S^{[n]}$ described above.
Associated to such a transformation, is a lagrangian correspondence
$\Z\subset S^{[n]}\times S^{[n]}$,
analogous to the Steinberg correspondence in the cartesian square 
of the cotangent bundle of a flag variety. 
The correspondence $\Z$ is reducible, with singular irreducible components.
In section \ref{sec-stratified-elemetary-transformations} we
relate the $K$-group of $\Z$ to that of a birational model
with smooth and transversal irreducible components.
In section \ref{sec-self-dual-stratified-elemetary-transformations} 
we treat the sequence $\{\tau_n\}$ of reflections 
with respect to line bundles
(Theorem \ref{thm-class-of-correspondence-in-stratified-elementary-trans}).
In section \ref{sec-stratified-elementary-trans-plus-2-vectors} we treat 
the sequence $\{\sigma_n\}$ of isometries
(Theorem \ref{thm-reflection-sigma-satisfies-main-conj}). 
The proof of each of these theorems is a rather involved K-theoretic 
calculation. The proofs rely heavily on our results in
\cite{markman-reflections,markman-diagonal}. 
In section \ref{sec-monodromy-of-A1-singularities} 
we relate the above mentioned correspondences to 
monodromy operators.

\smallskip
{\em Acknowledgments:} It is a pleasure to acknowledge fruitful conversations 
with Andrei C\u{a}ld\u{a}raru, Igor Dolgachev, Jim Humphreys, 
Daniel Huybrechts, Eduard Looijenga, Vikram Mehta, Ivan Mirkovic, 
Shigeru Mukai, 
Hiraku Nakajima, and Kota Yoshioka. I thank Yoshinori Namikawa
for pointing out a mistake, in an earlier version of the proof of Proposition
\ref{prop-rational-singularities-of-correspondence} in section 
\ref{sec-rational-singularities}.

\section{Statement of the Main Theorem}
\label{sec-summary-of-part-one}

Let $S$ be a K3 surface and $\LB$ a primitive ample line bundle on $S$. 
The Todd class of $S$ is $1+2\omega$, where $\omega$ is the fundamental class
in $H^4(S,\Integers)$. Its square root is $1+\omega$. 
Given a coherent sheaf $F$ on $S$ of rank $r$, we denote by 
\[
v(F)\ := \ ch(F)\sqrt{td_S} = (r,c_1(F),\chi(F)-r) 
\] 
its {\em Mukai vector} in  
\[
H^*(S,\Integers) \ = \ H^0(S,\Integers) \oplus H^2(S,\Integers) 
\oplus H^4(S,\Integers). 
\] 
Mukai endowed the cohomology group $H^*(S,\Integers)$ 
with a weight 2 polarized Hodge structure. 
The bilinear form is
\begin{eqnarray}
\label{eq-mukai-pairing}
\langle (r',c',s'), (r'',c'',s'') \rangle & = & c'c''-r's''-r''s',
\ \ \mbox{or equivalently}, \\
\nonumber
\langle \alpha, \beta\rangle & = & -\int_S \alpha^\vee\wedge\beta, \ \ 
\mbox{where},
\end{eqnarray}
\[
(\bullet)^\vee \ : \ H^*(S,\Integers) \ \rightarrow \ H^*(S,\Integers)
\]
is the duality acting by $-1$ on the second cohomology (sending the Mukai 
vector $(r,c_1,s)$ to $(r,-c_1,s)$). 
The Hodge filtration is induced by that of $H^2(S,\Integers)$. 
In other words, $H^{2,0}(S)$ is defined to be also the $(2,0)$-subspace of 
the complexified Mukai lattice. 

Denote by $\M(v):=\M_\LB(v)$ the moduli space of 
Gieseker-Simpson $\LB$-stable 
sheaves with Mukai vector $v$ of non-negative rank $r(v)\geq 0$. 
Assume, that the vector $v$ is primitive (indivisible). 
Mukai constructed a natural homomorphism
\[
\theta_v \ : \ v^\perp \ \ \rightarrow \ \ H^2(\M_\LB(v),\Integers)
\]
given by
\begin{equation}
\label{eq-mukai-homomorphism}
\theta_v(x) \ := \ \frac{1}{\rho}\left[p_{\M_*}((ch\E)\cdot\sqrt{td_S}\cdot 
\pi_S^*(x^\vee))
\right]_1,
\end{equation}
where $\E$ is a quasi-universal family of similitude $\rho$. 
In this paper, $\E$ will always be a universal family and $\rho=1$. 
Note, that the homomorphism $\theta_v$ extends to the whole Mukai lattice, but
the extension depends on the choice of $\E$. 
The following theorem is due to Mukai, O'Grady and Yoshioka:

\begin{thm}
\cite{yoshioka-irreducibility}
\label{thm-irreducibility}
Let $v$ be a Mukai vector of positive rank with $\langle v,v\rangle \geq 0$. 
For a generic $\LB$, 
$\M_\LB(v)$ is a smooth, non-empty, irreducible symplectic, 
projective variety of dimension $\dim(v)=\langle v,v\rangle+2$. 
The homomorphism (\ref{eq-mukai-homomorphism}) is an isomorphism of weight 2 
Hodge structures with respect to Beauville's bilinear form on 
$H^2(\M_\LB(v),\Integers)$ when $\dim(v)\geq 4$.  
\end{thm}

Let us review the construction of the class $\gamma(g)$
in (\ref{eq-action-of-stabilizer}). 
We work with more general moduli spaces $\M(v)$ 
as in Theorem \ref{thm-irreducibility}, admitting a universal sheaf
$\E_v$. However, we will review the definition of $\gamma(g)$
only for two kinds of isometries $g$. 
%
The first kind $\tau$ is the reflection with respect to
the $-2$ Mukai vector $v_0:=(1,0,1)$ of the trivial line bundle
\begin{equation}
\label{eq-tau}
\tau(v) \ \ := \ \ v + (v_0,v)\cdot v_0. 
\end{equation}
The second one $\sigma$ is the reflection with respect to the $+2$ 
Mukai vector $u_0=(1,0,-1)$ of the ideal sheaf of two points. 
\begin{equation}
\label{eq-sigma}
\sigma(w) \ = \ w - (w,u_0)u_0.
\end{equation}
The corresponding reflections $\sigma$ and $\tau$ 
commute and satisfy the relation 
\begin{equation}
\label{eq-relation-between-sigma-and-tau}
\sigma \ \ = \ \  (- \tau)^\vee.
\end{equation}

When $v$ is the Mukai vector of the Hilbert scheme $S^{[n]}$, 
the homomorphism $\gamma$ sends an isometry $g$ to a class
$\gamma_g(\E_v,\E_v)$ in the cohomology 
$H^{4n}(S^{[n]}\times S^{[n]},\Integers)$. 
Let $m$ be the complex dimension of a moduli space $\M(v)$. 
Dropping the assumption that
$\tau$ stabilize $v$, we define the class
$\gamma_{\tau}(\E_v,\E_{\tau(v)})$ in 
$H^{2m}(\M(v)\times \M(\tau(v)),\Integers)$ via
\begin{equation}
\label{eq-Chow-theoretic-formula-for-gamma-tau}
\gamma_{\tau}(\E_v,\E_{\tau(v)}) \ \ := \ \ 
c_m\left[- \ 
\pi_{13_!}
\left(
\pi_{12}^*(\tau[\E_v])^\vee\otimes
\pi_{23}^*(\E_{\tau(v)})
\right)
\right],
\end{equation}
where 
\begin{equation}
\label{eq-tau-of-universal-sheaf}
\tau[\E_v] \ \ := \ \ 
\E_v \ - \ 
p^!p_!\E_v,
\end{equation}
$\pi_{ij}$ is the projection from $\M(v)\times S \times \M(\tau(v))$
to the product of the $i$-th and $j$-th factors, $\pi_{ij_!}$ is the 
K-theoretic pushforward, $\otimes$ is the K-theoretic product 
and $p:\M(v)\times S\rightarrow \M(v)$ is the projection. 
Note, that $\tau[\E_v]$ is a K-theoretic representative, of a 
relative Fourier-Mukai transform, via a natural lift of
$\tau$ to an autoequivalence of the derived category of $S$
(see part I \cite{markman-part-one}). Lemma 
4.4 
in \cite{markman-part-one} implies the following equality.
\begin{equation}
\label{eq-chern-character-of-tau-E-v}
ch(\tau[\E_v])\cdot \sqrt{td_S} \ \ = \ \
(\tau\otimes id_{\M(v)})\left(ch(\E_v)\cdot \sqrt{td_S}\right).
\end{equation}

Similarly, we define the class
$\gamma_{\sigma}(\E_v,\E_{\tau(v)})$ in 
$H^{2m}(\M(v)\times \M(\sigma(v)),\Integers)$ via
\begin{equation}
\label{eq-Chow-theoretic-formula-for-gamma-sigma}
\gamma_{\sigma}(\E_v,\E_{\sigma(v)}) \ \ := \ \ 
c_m\left[- \ 
\pi_{13_!}
\left(
\pi_{12}^*(\sigma[\E_v])^\vee\otimes
\pi_{23}^*(\E_{\sigma(v)})
\right)
\right],
\end{equation}
where 
\begin{equation}
\label{eq-sigma-of-universal-sheaf}
\sigma[\E_v] \ \ := \ \ -\tau[\E_v]^\vee \ \ = \ \ 
\left[p^!p_!\E_v \ - \ \E_v 
\right]^\vee.
\end{equation}

We recall, next, the normalization of the chern character of the 
universal sheaf, introduced in \cite{markman-part-one}. 
Assume, that $\dim\M(v)>2$. 
Let $\eta_v$ be a $\RationalNumbers$-Cartier divisor on $\M(v)$, 
such that $(v,v)\cdot c_1(\eta_v)=\theta_v(v)$, where $\theta(v)$
is given in (\ref{eq-mukai-homomorphism}). The class
\begin{equation}
\label{eq-invariant-normalized-class-of-chern-character-of-universal-sheaf}
ch(p^*\eta_v)\cdot ch(\E_v) \cdot\pi_S^*\sqrt{td_S}
\end{equation}
in $H^*(\M(v)\times S,\RationalNumbers)$ 
is independent of the choice of a universal sheaf $\E_v$.

Theorem \ref{thm-Gamma-v-acts-motivicly} is the main result
of this paper. Acting by $(-1)^i$ on $H^{2i}(S^{[n]},\Integers)$, 
we get a graded-involution of the cohomology ring of 
$S^{[n]}$. We call it the Duality involution and denote it by 
$D_{S^{[n]}}$. The group $\Gamma$ has a natural orientation character 
$cov:\Gamma\rightarrow \{\pm 1\} \subset \ComplexNumbers^\times$
(see part I \cite{markman-part-one}).
This character sends $\tau$ to $1$ and $\sigma$ to $-1$. 
We regard the values of $cov$ in $\Integers/2\Integers$, 
when they appear in an exponent.
Consequently, we get another 
representation
\begin{eqnarray}
\label{eq-monodromy-is-gamma-times-cov}
\gamma\cdot cov  \ : \  \Gamma_v & \rightarrow & 
\Aut(H^*(S^{[n]},\Integers)_{\rm free}).
\\
\nonumber
g & \mapsto & \gamma(g)\cdot (D_{S^{[n]}})^{cov(g)}.
\end{eqnarray}

\begin{thm}
\label{thm-Gamma-v-acts-motivicly}
Let $v:=(1,0,1-n)$ be the Mukai vector of the Hilbert scheme $S^{[n]}$,
and $g$ an isometry in $\Gamma_v$.
\begin{enumerate}
\item
\label{thm-item-gamma-g-is-a-ring-isomorphism}
The class $\gamma_g(\E_v,\E_v)$ induces an automorphism $\gamma_g$ of the 
cohomology ring $H^*(\M(v),\Integers)_{\rm free}$, 
independent of the choice of $\E_v$.
If $g$ is a Hodge isometry, then 
$\gamma_g$ is an isomorphism of Hodge structures. 
\item
\label{thm-item-normalized-universal-chern-character-is-invariant}
The class 
(\ref{eq-invariant-normalized-class-of-chern-character-of-universal-sheaf}) 
is invariant with respect to the automorphism
\[
(\gamma_g\otimes g) \ : \ H^*(\M(v)\times S,\RationalNumbers) 
\ \LongIsomRightArrow \ 
H^*(\M(v)\times S,\RationalNumbers).
\]
\item
\label{thm-image-of-gamma-is-a-subgroup-of-monodromy-group}
The image of $\gamma\cdot cov$ 
is a subgroup of the monodromy group.
\end{enumerate}
\end{thm}

Theorem \ref{thm-Gamma-v-acts-motivicly} follows from
Theroems \ref{thm-class-of-correspondence-in-stratified-elementary-trans} 
and \ref{thm-reflection-sigma-satisfies-main-conj} 
and Lemma \ref{lemma-correspondence-is-a-monodromy-operator} 
via the reduction in part one of the paper \cite{markman-part-one}.

Finaly, we will need to use the following characterization of the class
$\gamma_g(\E_{1},\E_{2})$ 
(lemma 
5.2
in part one \cite{markman-part-one}). 
Let $(S_1,\LB_1)$ and $(S_2,\LB_2)$ be polarized K3 surfaces, 
$\M_{\LB_1}(v_1)$ and $\M_{\LB_2}(v_2)$ compact moduli spaces of 
stable sheaves, and $\E_i$ a universal sheaf over 
$\M_{\LB_i}(v_i)\times S_i$. 

\begin{new-lemma}
\label{lemma-recovering-f}
Suppose that $f: H^*(\M_{\LB_1}(v_1),\RationalNumbers) \rightarrow 
H^*(\M_{\LB_2}(v_2),\RationalNumbers)$ 
is a {\em ring} isomorphism, 
$g : H^*(S_1,\RationalNumbers) \rightarrow H^*(S_2,\RationalNumbers)$ 
a linear homomorphism, 
and  $f\otimes g$ 
maps the class 
(\ref{eq-invariant-normalized-class-of-chern-character-of-universal-sheaf})
in $H^*(\M_{\LB_1}(v_1)\times S_1,\RationalNumbers)$ to the 
class 
(\ref{eq-invariant-normalized-class-of-chern-character-of-universal-sheaf})
in $H^*(\M_{\LB_2}(v_2)\times S_2,\RationalNumbers)$.
Then
$
[f] \ = \ \gamma_g(\E_1,\E_2).
$
In particular, given $g$, a ring isomorphism $f$, satisfying
the condition above, is {\em unique} (if it exists).
\end{new-lemma}

\section{Stratified elementary transformations}
\label{sec-stratified-elemetary-transformations}
Let $X$ be a smooth variety with a holomorphic symplectic structure. 
A {\em stratified Mukai elementary transformation} is a certain 
birational transformation of $X$, resulting in a smooth variety $X'$ 
with a holomorphic symplectic structure \cite{markman-reflections}. 
We review the definition 
in section \ref{sec-review-of-stratified-elementary-transformations}.
A lagrangian correspondence $\Z\subset [X\times X']$ is introduced in
Definition \ref{def-Z-for-cotangent-bundles}. $\Z$ is 
analogous to the Steinberg correspondence, in the cartesian square 
of the cotangent bundle of a flag variety \cite{chriss-ginzburg}. 

The homomorphism 
\[
\Z_* \ : \ H^*(X,\Integers) \ \rightarrow \ H^*(X',\Integers),
\]
induced by the correspondence, will be a monodromy operator, 
in self-dual cases, where $X$ is naturally isomorphic to $X'$.
We will study the homomorphism $\Z_*$ via its K-theoretic analogue. 
This will be facilitated by resolutions $\beta_i:\Z^{[1]}_i\rightarrow \Z_i$ 
of the irreducible components of $\Z$. 
We will see in section \ref{sec-resolution-of-correspondences}, 
that the correspondence $\Z$ is the image of a
reducible variety $\Z^{[1]}$, which is a divisor,
with normal crossings, in a smooth variety. 
The $i$-th component of $\Z^{[1]}$ is the resolution $\Z^{[1]}_i$ of $\Z_i$. 
The main result of this section is proposition
\ref{prop-rational-singularities-of-correspondence}. 
It implies, that the composition 
\[
K^0(\Z) \ \LongRightArrowOf{\beta^!} \ K^0(\Z^{[1]}) \ 
\LongRightArrowOf{\beta_!} K_0(\Z)
\]
is equal to the natural homomorphism, 
from the Grothendiek K-group of vector bundles 
on $\Z$, to the Grothendiek K-group of coherent sheaves on $\Z$. 
This will enable us to carry-out calculations with $\Z^{[1]}$, 
instead of the more singular correspondence $\Z$. 
Proposition \ref{prop-rational-singularities-of-correspondence}
is related to the statement, that the components $\Z_i$, 
as well as their intersections, all have {\em rational singularities}
(see proposition \ref{prop-rational-singularities-in-the-Grassmannian-case}). 

In section \ref{sec-circular-complexes} we study the
prototypical case, where $X=T^*G(r,{\Bbb C}^h)$ 
is the cotangent bundle of a Grassmannian. 
In that case, we relate the
correspondence $\Z$ to the variety of circular complexes. 
We conclude, as a consequence, that the irreducible components of $\Z$
have rational singularities.
In section \ref{sec-resolution-of-correspondences} we construct the 
normal crossing model $\Z^{[1]}$ of $\Z$. 
In section \ref{sec-deformation-to-the-cone} we relate the general 
stratified elementary transformation, to the prototypical one, 
via  a deformation. This deformation enables us to reduce the proof of
Proposition \ref{prop-rational-singularities-of-correspondence} to
the case, where $X=T^*G(r,{\Bbb C}^h)$.
The latter case of Proposition 
\ref{prop-rational-singularities-of-correspondence} 
is proven in section \ref{sec-rational-singularities}. 
The action of the correspondence $\Z$, on the cohomology of $X$, is
briefly discussed in section 
\ref{sec-induced-action-on-cohomology}.

\subsection{Review of stratified elementary transformations}
\label{sec-review-of-stratified-elementary-transformations}

The prototypical example, of a stratified Mukai elementary transformation,
consists of the pair $X=T^*G(r,H)$ 
and $X'=T^*G(r,H^*)$, where $H$ is a vector space. 
Assume, that the dimension $h$ of $H$ is $\geq 2r$. 
Both $X$ and $X'$ are resolutions of the closure $\overline{\N}_r$ of 
a nilpotent orbit $\N_r$ in $\End(H)$. 
The orbit $\N_r$ consists of elements $\eta$ satisfying
\[
\eta^2=0 \ \ \mbox{and} \ \ \rank(\eta) = r. 
\]
Equivalently, the Jordan normal form of $\eta$ consists of 
$r$ two-by-two nilpotent blocks (and $h-2r$ zeroes along the diagonal). 
The closure $\overline{\N}_r$ is the union of $\N_t$, $0\leq t \leq r$.
The resolution morphism $T^*G(r,H) \rightarrow \overline{\N}_r$
sends a pair $(W,\psi: H/W\rightarrow W)$ to the composition
$H\rightarrow H/W \RightArrowOf{\psi} W \hookrightarrow H$ of $\psi$
with the natural projection and inclusion. 

The cotangent bundle $T^*G(r,H)$ is stratified by the rank of the
homomorphism $\psi$. Let $T^*G(r,H)^t$ be the locus, where $\psi$ has 
$t$-dimensional cokernel. This stratification has a recursive nature. 
Accordingly, the  birational transformation 
\begin{equation}
\label{eq-dual-springer-resolutions}
T^*G(r,H) \ \rightarrow \ \overline{\N}_r \ \leftarrow T^*G(r,H^*)
\end{equation}
is best described, in terms of a recursive structure, which we call
a {\em stratified dualizable collection} 
(Definition 2.3 in \cite{markman-reflections}, which is review below). 
Set $X(k):= T^*G(r-k,H)$, $0\leq k \leq r$, where $T^*G(0,H)$ is a point.
Set $\mu(k):=r-k$. 
Denote by $B^{[i]}X(k)$, $1\leq i \leq \mu(k)$ the following 
iterated blow-up of $X(k)$. 
$B^{[\mu(k)]}X(k)$ is the blow-up along the smallest stratum $X(k)^{\mu(k)}$
(the zero section). 
$B^{[i]}X(k)$ is the blow-up of $B^{[i+1]}X(k)$ along the
strict transform $B^{[i+1]}X(k)^i$ of $X(k)^i$. 
$B^{[i+1]}X(k)^i$ is smooth. In fact, it is a $G(i,h-2r+2i)$-bundle 
\begin{equation}
\label{eq-prototypical-f-k-i}
f_{k,i} \ : \ B^{[i+1]}X(k)^i \ \ \longrightarrow \ \ B^{[1]}X(k+i).
\end{equation}
The morphism $f_{k,i}$ sends a point $(W,\psi:H/W\rightarrow W)$ in 
$X(k)^i\setminus X(k)^{i+1}$, 
with $W\in G(r-k,H)$ and $\psi$ of rank $r-k-i$, 
to the pair 
$({\rm Im}(\psi),\bar{\psi}:H/{\rm Im}(\psi)\rightarrow {\rm Im}(\psi))$. 

The varieties $B^{[1]}T^*G(r,H)$ and $B^{[1]}T^*G(r,H^*)$ are isomorphic. 
In other words, the birational isomorphism 
(\ref{eq-dual-springer-resolutions}) is resolved as a
sequence of blow-ups, followed by a dual sequence of  blow-down
operations. Note, that the grassmannian bundle $X^{i}\setminus X^{i+1}$
(with bundle map $f_{0,i}$) is replaced by the dual grassmannian 
bundle. 

A {\em stratified dualizable collection} is a
collection $\{X(k)^t,f_{k,t}\}$, of connected symplectic varieties $X(k)$, 
each with a stratification
\begin{equation}
\label{eq-stratification-of-X-k}
X(k) = X(k)^0 \ \supset \ X(k)^1 \  \supset \ \cdots \  \supset \ 
X(k)^{\mu(k)}.
\end{equation}
Set $\mu:=\mu(0)$. Then $\mu(k)=\mu-k$.
It is convenient to arrange the data in an upper triangular 
$(\mu\!+\!1)\times(\mu\!+\!1)$-matrix with symplectic diagonal entries:

\begin{equation}
\label{eq-diagram-of-X-v}
\begin{array}{cccccc}
X(0)\supset      & X(0)^1 \supset       & \cdots             & \supset X(0)^t  
&   \cdots       & \supset X(0)^{\mu}
\\
                 & \downarrow           & 
\\
                 & X({1}) \supset & X({1})^1 \supset & 
&   \cdots       & \supset X({1})^{\mu-1}
\\
                 &                      & \downarrow
\\
                 &                      & X({2}) \supset   & 
&   \cdots       & \supset X({2})^{\mu-2}
\\
                 &                      &                        &
&                & \vdots
\\
                 &                      &                        &
&                & \downarrow
\\
                 &                      &                        &
&                & X(\mu)
\end{array}
\end{equation}

\noindent
Let $n$ be the codimension of $X^1$ in $X$. 
Every entry $X(k)^t$ admits a rational morphism to the symplectic diagonal 
entry $X(k+t)$. Each of these morphisms is resolved as a Grassmannian 
foliation 
\begin{equation}
\label{eq-f-k-t}
f_{k,t}:B^{[t+1]}X(k)^t\rightarrow B^{[1]}X(k+t)
\end{equation}
of iterated blow-ups of the strata, with $G(t,n+2k+2t-1)$-fibers. 
In particular, 
\[
f_{k,1}\ : \ B^{[2]}X(k)^1 \ \ \longrightarrow \ \ B^{[1]}X(k+1)
\]
is a $\PP^{n+2k}$-bundle, which we denote by
\begin{equation}
\label{eq-P-W-k+1}
\PP{W}_{k+1} \ \ := \ \ B^{[2]}X(k)^1.
\end{equation}
We do {\em not} assume, that the projective bundle $\PP{W}_{k+1}$
is the projectivization of a vector bundle. 
We continue this abuse of notation, and denote by
$G(t,W_{k+1})$ the relative grassmannian bundle of 
$\PP^{t-1}$ subspaces in the fibers of $\PP{W}_{k+1}$.
All the grassmannian bundles $f_{k,t}$, with 
the same value of $k+t$, are related by identifications
(which are part of the data of the collection)
\begin{equation}
\label{eq-relation-between-various-grassmannian-bundles}
B^{[t+1]}X(k)^t \ \ = \ \ G(t,W_{k+t}). 
\end{equation}

The codimension of $X(k)^t$ is equal to the dimension of the
grassmannian fibers of $f_{k,t}$. This implies the equality
\[
\dim X(k) \ \ = \ \ \dim X(0) \ - \ 2k(n+k-1).
\]
It follows, that the length $\mu$, of the stratification in
(\ref{eq-diagram-of-X-v}), is bounded by the inequality
\[
\dim X(0) \ \ \geq \ \ 2\mu(n+\mu-1),
\]
which is the condition, that $\dim X(\mu)\geq 0.$ 

The morphisms $f_{k,t}$ are compatible, with respect to the stratifications, 
in the following sense. 
We have the following equality of Cartier divisors on $B^{[t+1]}X(k)^t$.
\begin{equation}
\label{eq-conditioned-equality-of-catier-divisor}
B^{[t+1]}X(k)^t \cap B^{[t+1]}X(k)^r \ \ \ = \ \ \
f_{k,t}^{-1}\left(
B^{[1]}X(k+t)^{r-t}
\right),
\end{equation}
for $t<r\leq \mu(k)$.

The collection satisfies one additional condition, which we now recall
(Condition 2.5 in \cite{markman-reflections}). The exceptional divisor
$B^{[t]}X(k)^t$ is isomorphic to the 
projectivized normal bundle of $B^{[t+1]}X(k)^t$ in 
$B^{[t+1]}X(k)$. The latter is isomorphic to the 
projectivization $\PP{T}^*_{f_{k,t}}$ 
of the relative cotangent bundle of the Grassmannian fibration
(\ref{eq-f-k-t}) (see Lemma 2.2 in \cite{markman-reflections}).
The relative cotangent bundle is a homomorphism bundle from the
tautological quotient bundle to the tautological sub-bundle.
The condition asserts, that the stratification  
(\ref{eq-stratification-of-X-k}) induces on the exceptional divisor
$B^{[t]}X(k)^t$ the same 
determinantal stratification of the homomorphism bundle. 
Denote by $(\PP{T}^*f_{k,t})^i$ the locus, where the homomorphism has corank
$\geq i$. Then the condition asserts  the equality
\begin{equation}
\label{cond-induced-stratification-on-exceptional-divisor-is-determinantal}
B^{[t]}X(k)^t \ \cap \ B^{[t]}X(k)^i \ \ = \ \ 
(\PP{T}^*_{f_{k,t}})^i, \ \ \ 0\leq i \leq t-1.
\end{equation}
The last condition relates the 
bundles $\PP{W}_k$, for different values of $k$, along the intersection
of the exceptional divisors (Lemma 
\ref{lemma-grassmannian-fibers-are-homologous-to-sub-grassmannians} below). 

Let $V$ be a vector space of dimension $N\geq 2t$.
A flag $U\subset W \subset V$, with $\dim(U)=t-k$ and
$\dim(W)=N-(t-k)$, determines the subvariety 
$G(k,W/U)\subset G(t,V)$. 
The flag variety $Flag(t-k,N-(t-k),V)$ 
is a component of the Hilbert scheme
of $G(t,V)$, parametrizing subvarieties isomorphic to $G(k,N-2(t-k))$.

Consider the following commutative diagram, for $1\leq k \leq t-1$
and $t\leq \mu$. 
\[
\begin{array}{ccccccc}
G(k,W_k) & = & 
B^{[k+1]}X^k & \LongRightArrowOf{f_{0,k}} & B^{[1]}X(k)
\\
& & \cup & & \cup
\\
& & B^{[k+1]}X^{k\cap t} &
\LongRightArrowOf{g} &
B^{[1]}X(k)^{t-k} &
\RightArrowOf{\beta} &
B^{[t-k+1]}X(k)^{t-k}
\\
& & \hspace{1ex} \ \downarrow \ \phi & & 
f_{k,t-k}\circ\beta \ \downarrow \ \hspace{7ex}
&   \stackrel{f_{k,t-k}}{\swarrow} 
\\
G(t,W_t) & = & 
B^{[t+1]}X^t & \LongRightArrowOf{f_{0,t}} & B^{[1]}X(t)
\end{array}
\]

\noindent
Above, $g$ is the restriction of the Grassmannian bundle
$f_{0,k}$. The morphism $\phi$ is the composition of the inclusion
$B^{[k+1]}X^{k\cap t}\subset B^{[k+1]}X^{t}$, followed by the blow-down
morphism $B^{[k+1]}X^{t}\rightarrow B^{[t+1]}X^{t}$.

\begin{new-lemma}
\label{lemma-grassmannian-fibers-are-homologous-to-sub-grassmannians}
There exists a natural morphism
\[
\eta \ : \ B^{[1]}X(k)^{t-k} \ \ \longrightarrow \ \ 
Flag(t-k,\rank(W_t)-(t-k),W_t).
\]
Moreover, $\eta$ determines a $G(k,n+2k-1)$-subbundle
of $\beta^*f_{k,t-k}^*G(t,W_t)$, which is naturally isomorphic to 
the $G(k,n+2k-1)$-bundle $g$. Consequently, $\phi$ embeds 
each fiber of $g$ in a fiber of $f_{0,t}$. 
In particular, the generic fiber of $f_{0,k}$ is homologous 
in $X(0)$ to a Grassmannian subvariety in a generic fiber of $f_{0,t}$.
\end{new-lemma}

\noindent
{\bf Proof:}
The condition 
(\ref{cond-induced-stratification-on-exceptional-divisor-is-determinantal})
identifies $B^{[1]}X(k)^{t-k}$ with 
the top iterated blow-up $B^{[1]}(\PP{T}^*_{f_{k,t-k}})$, of the
relative projectivized cotangent bundle.
The morphism $\eta$ is the natural extension, of the morphism,
which sends a point $(U,\psi:V/U\rightarrow U)$
in the open stratum
$\PP{T}^*_{f_{k,t-k}}\setminus (\PP{T}^*_{f_{k,t-k}})^1$
where $\psi$ is surjective, to the flag 
$U\subset \ker[V\rightarrow V/U\RightArrowOf{\psi} U]\subset V$
in the fiber $V$ of $W_t$.

The morphism $\phi$ factors through 
$B^{[t]}X^t\cap B^{[t]}X^k$. The latter is identified with 
the stratum 
$(\PP{T}^*_{f_{0,t}})^k$, by equality 
(\ref{cond-induced-stratification-on-exceptional-divisor-is-determinantal}).
The identification, of the two $G(k,n+2k-1)$-bundles, now follows from
Lemma 2.7 in \cite{markman-reflections}.
\EndProof

\medskip
We study in this section a correspondence $Z \subset X\times X'$. 
$Z$ has pure dimension equal to that of $X$. 
In the prototypical example above,
the correspondence $Z$ admits three equivalent descriptions. 

\begin{defi}
\label{def-Z-for-cotangent-bundles}
{\rm 
\begin{enumerate}
\item
\label{correspondence-is-a-fiber-product}
$Z$ is the fiber product of $T^*G(r,H)$ and $T^*G(r,H^*)$ over 
$\overline{\N}_r$. 
\item
\label{correspondence-is-union-of-conormal-varieties}
$Z$ is the union of the irreducible components $Z_i$, $0\leq i \leq r$.
The component $Z_i$ is the conormal variety in
$T^*[G(r,H)\times G(r,H^*)]$, to the incidence subvariety 
\begin{equation}
\label{eq-incidence-variety-in-product-of-grassmannians}
I_i \ := \ \{(U,W) \ : \  \dim(U\cap W^\perp)\geq r-i\}
\end{equation}
of $G(r,H)\times G(r,H^*)$.
\item
\label{correspondence-fiber-products-of-dual-grassmanian-bundles}
$Z$ is the union of the irreducible components $Z_i$, $0\leq i \leq r$.
The component $Z_i$ is the closure, in $X \times X'$, 
of the fiber product of dual grassmannian fibrations
\[
[X^i\setminus X^{i+1}] \times_{X(i)\setminus X(i)^1} 
[(X')^i\setminus (X')^{i+1}].
\]
Note, that the base $X(i)\setminus X(i)^1$ is the open stratum 
in both $X(i)$ and $X'(i)$.
\end{enumerate}
}
\end{defi}

For a general stratified dualizable collection $\{X(k)^t,f_{k,t}\}$,
only description 
\ref{correspondence-fiber-products-of-dual-grassmanian-bundles}
is available (with $r=\mu$). We use description 
\ref{correspondence-fiber-products-of-dual-grassmanian-bundles} 
as the initial definition in the prototypical example.
The equivalence, of descriptions
\ref{correspondence-is-union-of-conormal-varieties} and 
\ref{correspondence-fiber-products-of-dual-grassmanian-bundles},
follows from Claim \ref{claim-transversality}, 
and the analogous equivalence for the variety of circular complexes, 
proven in \cite{strickland}. 
Note, that the {\em conormal variety} of $I_i$ is defined to be the closure, 
in the cotangent bundle of $G(r,H)\times G(r,H^*)$, of 
the conormal bundle of the smooth locus of $I_i$. 
The equivalence of descriptions
\ref{correspondence-is-a-fiber-product} and 
\ref{correspondence-fiber-products-of-dual-grassmanian-bundles}
is proven in Corollary \ref{cor-Z-is-a-fiber-product}. 
Description \ref{correspondence-is-union-of-conormal-varieties}
will be used in section \ref{sec-circular-complexes},
in order to prove that the varieties $Z_i$
are not too singular. We will not use description 
\ref{correspondence-is-a-fiber-product}.

\subsection{Circular complexes}
\label{sec-circular-complexes}

We consider in this section the prototypical case $X:=T^*G(r,H)$. 
Let $Z_{i,j}$ be the intersection $Z_i\cap Z_j$. 
We prove in this section the following proposition.

\begin{prop}
\label{prop-rational-singularities-in-the-Grassmannian-case}
The varieites $Z_i$, as well as $Z_{i,j}$, are all normal, 
Cohen-Macaulay, and have rational singularities. 
\end{prop}

We first exhibit the incidence variety $I_i$, given in 
(\ref{eq-incidence-variety-in-product-of-grassmannians}), 
as a determinantal variety. 
Let $\tau_H$ and $q_H$ be the tautological sub and quotient bundles over 
$G(r,H)$. Denote by $\pi_1: G(r,H)\times G(r,H^*)\rightarrow G(r,H)$
the projection and by $\pi_2$ the projection onto $G(r,H^*)$.
Let $g:\pi_1^*\tau_H \rightarrow \pi_2^*(\tau_{H^*})^*$ be the composition 
\begin{equation}
\label{eq-homomorphism-between-tautological-subundle-of-dual-grassmannians}
\pi_1^*\tau_H \hookrightarrow (H)_{G(r,H)\times G(r,H^*)} \rightarrow 
\pi_2^*(\tau_{H^*})^*
\end{equation}
of the natural inclusion and quotient homomorphisms. 
The incidence variety $I_i$ is the $i$-th 
degeneracy locus of the section $g$. 

\begin{claim}
\label{claim-transversality} 
The section $g$ is transversal 
to the determinantal stratification of the total space of the bundle 
$\Hom(\pi_1^*\tau_H,\pi_2^*(\tau_{H^*})^*)$. 
\end{claim}

\noindent
{\bf Proof:}
We need to show, that $g$ induces a submersive morphism to $\Hom(U,Q)$,
upon a choice of a local trivialization of 
$\pi_1^*\tau_H$ and $\pi_2^*(\tau_{H^*})^*$, 
modeled after the vector spaces $U$ and $Q$.
We begin with a construction, providing a more canonical formulation of the 
transversality. 

Let $U$ and $Q$ be $r$-dimensional vector spaces,
$Inj(U,H)\subset \Hom(U,H)$ the open subset of injective homomorphisms, 
and $Sur(H,Q)\subset \Hom(H,Q)$ the open subset of surjective homomorphisms.
Note, that $Inj(U,H)$ and $Sur(H,Q)$ are the frame bundles of
$\tau_{G(r,H)}$ and $q_{G(h-r,H)}$ respectively. 
We identify $G(r,H^*)$ with $G(h-r,H)$, replacing $\tau^*_{G(r,H^*)}$ with
$q_{G(h-r,H)}$. 
Consider the bundle map
\begin{eqnarray}
\label{eq-product-of-frame-bundles}
Inj(U,H) \times Sur(H,Q) & 
\LongRightArrowOf{p}
& G(r,H)\times G(h-r,H)
\\
\nonumber
(\alpha,\beta) & \mapsto & ({\rm Im}(\alpha),\ker(\beta))
\end{eqnarray}
The pullback of $\Hom(\pi_1^*\tau,\pi_2^*q)$ to 
$Inj(U,H) \times Sur(H,Q)$
trivializes
\begin{eqnarray*}
p^*\Hom(\pi_1^*\tau,\pi_2^*q) & \LongIsomRightArrow & 
Inj(U,H) \times Sur(H,Q)\times\Hom(U,Q)
\\
(\alpha,\beta,\eta:{\rm Im}(\alpha)\rightarrow H/\ker(\beta))
& \mapsto & (\alpha,\beta,\bar{\beta}\circ\eta\circ\alpha),
\end{eqnarray*}
where $\bar{\beta}:H/\ker(\beta)\rightarrow Q$ is induced by $\beta$. 
Choose a local section $s$ of the bundle map $p$, given in
(\ref{eq-product-of-frame-bundles}). Let
$p^*(g)$ 
be the pullback of $g$ via the bundle map $p$.
We summarize the construction in the diagram
\[
\begin{array}{ccccc}
\Hom(U,Q) & \LongLeftArrowOf{f} & 
p^*\Hom(\pi_1^*\tau,\pi_2^*q) & 
\stackrel{\LongLeftArrowOf{p^*g}}{\longrightarrow}& 
Inj(U,H) \times Sur(H,Q)
\\
& & 
\downarrow & 
& 
p \ \downarrow \ \uparrow \ s
\\
& & \Hom(\pi_1^*\tau,\pi_2^*q) & 
\stackrel{\LongLeftArrowOf{g}}{\longrightarrow}
& G(r,H)\times G(h-r,H).
\end{array}
\]
Above, $f$ is the fibration induced by the trivialization. 
The composition $f\circ p^*g\circ s$
pulls back the determinantal stratification of $\Hom(U,Q)$
to the incidence stratification of $G(r,H)\times G(h-r,H)$. 
It remains to prove, that the 
composition $f\circ p^*g\circ s$ is a smooth morphism.

The morphism $p$ is equivariant, with respect to the diagonal action of 
the group $GL(H)$ on $Inj(U,H)\times Sur(H,Q)$
and $G(r,H)\times G(h-r,H)$. The incidence varieties $I_i$ in 
$G(r,H)\times G(h-r,H)$ are closures of $GL(H)$-orbits and 
the complement of $I_{r-1}$ is the open dense orbit. 
Let $Y_s$ be the open subset, over which the composition 
$f\circ p^*g\circ s$ is submersive. Then $Y_s=Y_{s'}$
if the sections $s$ and $s'$ are conjugate under the $GL(H)$-action. 
It follows, that $Y_s$ is $GL(H)$-invariant. The equality, 
$Y_s=G(r,H)\times G(h-r,H)$ would follow, once we prove that 
$Y_s$ contains the smallest stratum $I_0$ (the flag variety). 

The composition $f\circ p^*g$ sends $(\alpha,\beta)$ to $\beta\circ\alpha$.
An elementary calculation checks, that $f\circ p^*g$ is submersive.
The composition  $f\circ p^*g\circ s$ is submersive at $x$, if and only if
the following equality holds
\[
{\rm Im}(d_xs) \ + \ \ker d_{s(x)}(f\circ p^*g) \ \ = \ \ T_{s(x)}
[Inj(U,H)\times Sur(H,Q)]. 
\]
We will prove this equality, 
for $x$ in $I_0$, by proving the containment
\[
\ker d_{s(x)}(f\circ p^*g) \ \ \supset \ \ \ker(d_{s(x)}p).
\] 
Observe, that the image $s(x)$, of $x\in I_0$, is a pair $(\alpha,\beta)$,
such that $\beta\circ\alpha=0$. 
The above containment follows immediately from the identification 
of each of the two subspaces below. 
\begin{eqnarray*}
\ker(d_{(\alpha,\beta)}p) & = & 
\left\{(\dot{\alpha},\dot{\beta}) \ : \ {\rm Im}(\dot{\alpha}) \subset  
{\rm Im}(\alpha) \ \ 
\mbox{and} \ \ \ker(\dot{\beta}) \supset \ker(\beta)\right\}
\\
\ker d_{(\alpha,\beta)}(f\circ p^*g) & = & 
\left\{
(\dot{\alpha},\dot{\beta}) \ : \ 
\beta\circ \dot{\alpha} + \dot{\beta}\circ \alpha \ = \ 0
\right\}.
\end{eqnarray*}
\EndProof

\smallskip
Proposition 
\ref{prop-rational-singularities-in-the-Grassmannian-case}
follows from the analogous result for determinantal varieties,
stated in Theorem \ref{thm-mehta-trivedi}.
Let $V$ and $U$ be finite dimensional vector spaces of ranks $n$ and $m$ 
respectively. Denote by $W$ the subvariety of $\Hom(U,V)\times \Hom(V,U)$ 
of pairs $(\varphi_1,\varphi_2)$ satisfying 
\[
\varphi_1\circ\varphi_2=0 \ \ \mbox{and} \ \ \varphi_2\circ\varphi_1=0.
\]
The group $GL(U)\times GL(V)$ acts on $W$ by 
\[
(a,b)\cdot (\varphi_1,\varphi_2)
\ = \ 
(b\circ\varphi_1\circ a^{-1},a\circ\varphi_2\circ b^{-1}).
\]
Let $W_0(k_1,k_2)$ be the $GL(U)\times GL(V)$ orbit in $W$, 
consisting of pairs $\varphi_1,\varphi_2$ of ranks $(k_1,k_2)$. The closure
$W(k_1,k_2)$ consists of pairs, satisfying the rank inequalities
$\rank(\varphi_i)\leq k_i$. The ideal $I(k_1,k_2)$ of $W(k_1,k_2)$,
in the affine coordinate ring of $\Hom(U,V)\times \Hom(V,U)$, is
the sum of 
1) the ideal of $W$, 
2) the pull back of the ideal of the determinantal locus $\Hom(U,V)^{n-k_1}$, 
and 
3) the pull back of the ideal of the determinantal locus $\Hom(V,U)^{m-k_2}$ 
(see \cite{strickland}).

Assume now, that $n\leq m$. Note, that the inequality $k_1+k_2\leq n$ holds.
When $k_1+k_2= n$, then $W(k_1,k_2)$ admits another description.
$W(t,n-t)$ is the conormal variety to the determinantal variety 
$\Hom(U,V)^{n-t}$, where $\varphi_1$ has a cokernel of 
dimension $\geq n-t$ (see \cite{strickland}). 
Furthermore, for $0<t_1<t_2<n$, we have the scheme theoretic  
equality 
\[
W(t_1,n-t_1) \ \cap \ W(t_2,n-t_2) \ \ = \ \ W(t_1,n-t_2).
\]
Consequently, given an ascending sequence 
$t_1<t_2< \cdots < t_\ell$, we get the equality
\begin{equation}
\label{eq-intersections-of-components-of-the-variety-W}
W(t_1,n-t_1) \ \cap \ W(t_2,n-t_2) \ \cap \ \cdots \ \cap \ W(t_\ell,n-t_\ell)
\ \ = \ \ W(t_1,n-t_\ell).
\end{equation}

Claim \ref{claim-transversality} and 
description \ref{correspondence-is-union-of-conormal-varieties},
of the correspondences $Z_t$ in definition \ref{def-Z-for-cotangent-bundles}, 
imply that the variety $Z_{t_1,t_2}$ has, locally, the structure 
of a smooth fibration over $W(t_1,n-t_2)$ (in the special case $n=m=r$). 
Proposition 
\ref{prop-rational-singularities-in-the-Grassmannian-case}
follows from the following theorem of Mehta and Trivedi.

\begin{thm}
\label{thm-mehta-trivedi}
\cite{mehta-trivedi,magyar}
The varieties $W(k_1,k_2)$ are all normal, Cohen-Macauley, and have
rational singularities. 
\end{thm}

\subsection{Resolutions of the Steinberg correspondences}
\label{sec-resolution-of-correspondences}

Let $\{X(r)^t,f_{r,t}\}$, $0\leq r \leq \mu$, $0\leq t\leq \mu-r$, be
a dualizable collection, with dual collection
$\{X'(r)^t,f'_{r,t}\}$. 
Set $X:=X(0)^0$ and $X':=X'(0)^0$. 
Let $n$ be the codimension of $X^1$ in $X$. 
We recall, in this section, the construction of the resolutions 
$Z^{[1]}_t\rightarrow Z_t$ of the correspondences in definition
\ref{def-Z-for-cotangent-bundles}.

The fiber product 
\begin{equation}
\label{eq-Delta-t-is-a-fiber-product}
\Delta_t^{[t+1]} \ := \ 
B^{[t+1]}X^t\times_{B^{[1]}X(t)}B^{[t+1]}(X')^t
\end{equation}
is a bundle over $B^{[1]}X(t)$ with fibers of the form 
$G(t,H)\times G(t,H^*)$, 
where $H$ is a $(2t+n-1)$-dimensional vector space. 
Note, that when $n=1$, the $G(t,2t)$-bundle
$B^{[t+1]}X^t\rightarrow B^{[1]}X(t)$ is 
self-dual and $X'(0)^t = X(0)^t$, for all $t$.

The product $G(t,H)\times G(t,H^*)$ is stratified 
by the dimension of the intersection $W'\cap (W'')^\perp$,
$W'\in G(t,H)$, $W''\in G(t,H^*)$.
This is the determinantal stratification of the global section 
of $\Hom(\pi_1^*\tau_{G(t,H)},\pi_2^*\tau_{G(t,H^*)}^*)$, given in 
(\ref{eq-homomorphism-between-tautological-subundle-of-dual-grassmannians}).
The incidence variety 
\begin{equation}
\label{eq-incidence-variety-I-t-i}
(I_t)^i \ \ \subset \ \ \Delta_t^{[t+1]}, \ \ \ \ 0\leq i < t,
\end{equation}
is defined to be the stratum, 
where the dimension of the intersection $W'\cap (W'')^\perp$
is $\geq (\!t-i\!)$.
It is the relative version of the one given in 
(\ref{eq-incidence-variety-in-product-of-grassmannians}). 
Denote by 
\begin{equation}
\label{eq-top-iterate-of-component-of-correspondence}
\Delta_t^{[i]}, \ \ 1\leq i \leq t-1,
\end{equation}
the iterated blow-up of $\Delta_t^{[t+1]}$, starting from the smallest 
stratum $(I_t)^0$ 
(the bundle of flag varieties $Flag(t,\chi(v)+t,H)$ over
$B^{[1]}X(t)$) and proceeding in
decreasing dimension of the intersection $W'\cap (W'')^\perp$. 
The fiber of $\Delta_t^{[1]}$ over $\Delta_t^{[t+1]}$
is the space of complete colineations 
$(\varphi_1,\varphi_2, \dots, \varphi_k)$ from
$W'$ to $(W'')^*$. The first homomorphism 
$\varphi_1:W'\rightarrow (W'')^*$ is the composition
(\ref{eq-homomorphism-between-tautological-subundle-of-dual-grassmannians}). 
Each successive colineation $\varphi_{i+1}$ is in 
$\PP\Hom(ker(\varphi_i),\coker(\varphi_i))$ and the last one
$\varphi_k$ is an isomorphism. Truncating the predetermined $\varphi_1$,
we conclude that the fiber is the space of complete
colineations from $W'\cap (W'')^\perp$ to $H/[W'+(W'')^\perp]$.

{\em Note: $\Delta_t^{[i]}$ is not contained in 
$B^{[i]}X\times B^{[i]}X'$, for $1\leq i\leq t-1$. 
The proper transform of $\Delta_t^{[t+1]}$ in
$B^{[t]}X\times B^{[t]}X'$ is the fiber product 
$B^{[t]}X^t\times_{B^{[1]}X(t)} B^{[t]}(X')^t$.
The latter is a fiber product of bundles of projectivized 
cotangent bundles of grassmannians. Its dimension is larger than 
the dimension of $\Delta_t^{[t+1]}$. 
Consequently, the proper transform of $\Delta_t^{[t+1]}$ in
$B^{[i]}X\times B^{[i]}X'$ has dimension larger than that of 
$\Delta_t^{[t+1]}$.
}

Denote by $\Z_t$ the set theoretic image of $\Delta_t^{[t+1]}$ in
$X\times X'$,  
endowed with the {\em reduced induced subscheme structure}. Let
\[
\Z \ := \ \cup_{t=1}^\mu\Z_t   \ \subset \ X\times X'
\]
be the set theoretic union, 
endowed with the {\em reduced induced subscheme structure}. 
It is precisely the correspondence,
in definition \ref{def-Z-for-cotangent-bundles} part 
\ref{correspondence-fiber-products-of-dual-grassmanian-bundles}.


There is a natural way to glue the top iterates $\Delta_t^{[1]}$
as a reduced and reducible connected variety 
\[
\Z^{[1]} \ := \ \cup_{i=1}^\mu  \Delta_t^{[1]} 
\]
with normal crossings. 
$\Z^{[1]}$ can be defined abstractly, by constructing 
the natural isomorphism, between the variety
defined in part 
\ref{lemma-item-pull-back-of-BN-divisors-to-correspondence} of Lemma 
\ref{lemma-structure-of-correspondence}
and the variety
defined in part 
\ref{lemma-item-pullback-of-incidence-variety-is-intersection-of-components} 
of that Lemma. 
We will, however, use an easier definition; we construct $\Z^{[1]}$ 
as a divisor with normal crossings in an auxilary smooth variety
(Lemma \ref{lemma-two-descriptions-of-the-components-of-the-correspondence}
and Lemma \ref{lemma-structure-of-correspondence}
part \ref{lemma-item-top-correspondence-has-normal-crossing}). 

Let
$\X\rightarrow B$ be a family of smooth symplectic
varieties over a smooth one-dimensional base $B$.
Assume the following condition. 

\begin{condition}
\label{cond-extension-class-is-non-trivial}
\begin{enumerate}
\item
The fiber over $0\in B$ is isomorphic to $X$.
\item
The extension class in $H^1(X,TX)$ of
\[
0 \rightarrow TX \rightarrow \restricted{T\X}{X}\rightarrow 
\StructureSheaf{X}
\rightarrow 0
\]
pairs with the symplectic structure to give a class 
$\alpha\in H^1(X,T^*X)$. 
Assume, that this class restricts to a 
non-vanishing class in $H^1(\PP^{n},T^*\PP^{n})$ on every $\PP^{n}$ fiber of 
$[X^1\setminus X^2]\rightarrow 
[X(1)\setminus X(1)^1]$. 
\end{enumerate}
\end{condition}

The class $\alpha$, in the above condition, restricts to 
a non-vanishing class in $H^{1,1}$ of every grassmannian fiber of 
$[X^t\setminus X^{t+1}]\rightarrow 
[X(t)\setminus X(t)^1]$, for $1<t\leq \mu$ as well, by Lemma
\ref{lemma-grassmannian-fibers-are-homologous-to-sub-grassmannians}. 

A deformation, satisfying Condition 
\ref{cond-extension-class-is-non-trivial}, exists whenever $X$ 
is an irreducible symplectic 
projective variety (see section 2.5 in \cite{markman-reflections}). 
In the prototypical example $X=T^*G(r,H)$, such a deformation is given by
the total space $\X$ of the vector bundle $E$ in 
the non-trivial extension 
\begin{equation}
\label{eq-non-trivial-extension-is-deformation-of-cotangent-bundle}
0 \rightarrow T^*G(r,H) \rightarrow E  \rightarrow \StructureSheaf{G(r,H)}
 \rightarrow 0.
\end{equation}
The base $B$ is $H^0(\StructureSheaf{G(r,H)})$. The generic fiber 
is the twisted cotangent bundle of $G(r,H)$. 

We will need an explicit embedding of the vector bundle $E$, as a 
subbundle of the trivial $\End_0(H)$-bundle over $G(r,H)$.
$\End_0(H)$ is the algebra of traceless endomorphisms. The tautological flag
$\tau_{G(r,H)}\subset H$ determines a bundle of parabolic subalgebras 
in $\End_0(H)$. $T^*G(r,H)$ is the nilpotent radical. Conceptually, 
$E$ is the extension of the bundle of nilpotent radicals by the 
bundle of centers of the Levi subalgebras. 
Given an endomorphism $\varphi$, leaving 
a subspace $W$ invariant, let 
$\restricted{\varphi}{W}$ be its restriction to $W$ and 
$\bar{\varphi}$ the induced endomorphism of $H/W$. 
The fiber of the 
bundle $E$ is given explicitly, at a point $W$ of $G(r,H)$, by 
\begin{eqnarray}
\label{eq-E-is-a-subbundle-of-End-H}
E & \subset & \End_0(H)_{G(r,H)}
\\
\nonumber
E_W & := &
\{
\varphi \in \End_0(H) \ : \ \varphi(W) \subset W, \ 
\restricted{\varphi}{W}\in {\rm span}(I_W), \ 
\bar{\varphi} \in {\rm span}(I_{H/W})
\}. 
\end{eqnarray}
The homomorphism $E\rightarrow \StructureSheaf{G(r,H)}$ is the
natural homomorphism $\varphi\mapsto \restricted{\varphi}{W}$
from $E$ to the center of $\End(\tau_{G(r,H)})$.

Set $\X^t:= X^t$ for $0<t\leq \mu$.  
Denote by $B^{[1]}\X$ the top iterated blow-up of $\X$, 
with respect to the stratification of the special fiber $X$. 
Then we have a dual iterated  blow-down
$B^{[1]}\X\rightarrow \X'$, where $\X'$ is a deformation of
$X'$ satisfying condition
\ref{cond-extension-class-is-non-trivial} 
(\cite{markman-reflections} Section 2.5). 
In particular, $B^{[1]}\X$ and $B^{[1]}\X'$
are isomorphic. Let 
\begin{equation}
\label{eq-top-iterated-blow-up-of-Y}
B^{[1]}\Y \ \subset \ B^{[1]}\X\times_B B^{[1]}\X'
\end{equation}
be the graph of the isomorphism. Denote by
\begin{equation}
\label{eq-Y}
\Y
\end{equation} 
the closure in $\X\times_B\X'$, of the graph
of the isomorphism between $\X\setminus X$ and $\X'\setminus X'$. 
We define $\Z^{[1]}$ as the fiber of $B^{[1]}\Y$ over $0\in B$. 
\begin{equation}
\label{eq-definition-of-top-iterated-blow-up-of-Z}
\Z^{[1]} \ \subset \ B^{[1]}\Y. 
\end{equation}
Set theoretically, $\Z$ is the fiber of $\Y$ over $0\in B$. 
We will later show, that 
the scheme structure of the fiber of $\Y$ over $0\in B$ agrees with
that of $\Z$ (Proposition
\ref{prop-rational-singularities-of-correspondence} 
part \ref{prop-item-Z-is-a-fiber}). 
Let $B^{[1]}\Y^t$ be the image of the divisor 
$B^{[1]}\X^t$, under the isomorphism
$B^{[1]}\Y\cong B^{[1]}\X$. Clearly, $B^{[1]}\Y^t$
is an irreducible component of $\Z^{[1]}$, 
which will be denoted by $\Z_t^{[1]}$ as well.

\begin{new-lemma}
\label{lemma-two-descriptions-of-the-components-of-the-correspondence}
The irreducible component $\Z_t^{[1]}$, 
of the fiber $\Z^{[1]}$ of $B^{[1]}\Y$ over $0\in B$, 
is isomorphic to $\Delta_t^{[1]}$. 
\end{new-lemma}

\noindent
{\bf Proof:}
See the proof of Theorem 1.2 in \cite{markman-reflections}. 
The key ingredient is Proposition 2.9 in \cite{markman-reflections}. 
\EndProof

\begin{rem}
\label{rem-notation-Z-vs-Delta}
{\rm
(Summary of notation)
We will use both $\Delta_t^{[1]}$ and
$\Z_t^{[1]}$ to denote the irreducible component
of the normal crossing model $\Z^{[1]}$ of the
correspondence. This is justified by Lemma 
\ref{lemma-two-descriptions-of-the-components-of-the-correspondence}.
The varieties $\Delta_t^{[1]}$, $\Z_t^{[1]}$, $B^{[1]}\Y^t$, 
$B^{[1]}\X^t$, and $B^{[1]}(\X')^t$ are all isomorphic.
On the other hand, the intermediate blow-up varieties
$\Delta_t^{[i]}$, $2\leq i \leq t-1$, and the exceptional divisors 
$B^{[j]}\X^t$, $2\leq j \leq t$, are {\em not} isomorphic 
for any pair of intermediate indices $(i,j)\neq (1,1)$
(see section \ref{sec-incidence-divisor}).
This difference is the main reason, that we chose the letter
$\Delta$, rather than $\Z$, in the notation of the varieties
$\Delta_t^{[i]}$. The letter $\Z$ will {\em never} be used for an 
intermediate blow-up. The letter $\Z$ is 
reserved for the correspondence $\Z$ or its components $\Z_t$, 
as well as for the  
normal crossing model $\Z^{[1]}$ or its components $\Z^{[1]}_t$.
Similarly, we will not consider the intermediate blow-ups of $\Y$. 
}
\end{rem}

\begin{new-lemma}
\label{lemma-structure-of-correspondence}
\begin{enumerate}
\item
\label{lemma-item-components-of-correspondence-are-smooth}
The varieties $\Delta_t^{[i]}$ are smooth, for $1\leq i \leq t-1$.
\item
\label{lemma-item-top-correspondence-has-normal-crossing}
$\Z^{[1]}$ has pure dimension $\dim(X)$. It can be
embedded as a divisor with normal crossing
in a smooth variety of dimension $\dim(X)+1$. In particular,
if $i\neq j$, then the pairwise scheme-theoretic intersections 
\[
\Delta_{i\cap j}^{[1]} \ := \ \Delta_i^{[1]}\cap\Delta_j^{[1]}
\] 
are smooth divisors in each. 
\item 
\label{lemma-item-fibers-of-beta-are-embedded-in-a-single-component}
Let $z$ be a point of $\Z$ and $t$ the maximal index, among the indices of
components $\Z_i$ containing $z$. 
The fiber of $\beta:\Z^{[1]}\rightarrow \Z$ over $z$ 
is scheme-theoretically equal to the fiber of $ \Delta^{[1]}_t$ over $z$.
\item
\label{lemma-item-pull-back-of-BN-divisors-to-correspondence}
The composition
$\Delta_t^{[1]} \ \rightarrow \ \Delta_t^{[t+1]} \ \rightarrow \ 
B^{[t+1]}X^t \ \hookrightarrow \ B^{[t+1]}X$
pulls back (scheme theoretically) 
the divisor $B^{[t+1]}X^i$, \ $i\geq t+1$, 
to the divisor $\Delta_t^{[1]}\cap\Delta_i^{[1]}$.
\item
\label{lemma-item-pullback-of-incidence-variety-is-intersection-of-components}
The intersection $\Delta_t^{[1]}\cap\Delta_i^{[1]}$, $0\leq i < t$, 
is the exceptional divisor of $\Delta_t^{[1]}$ corresponding,
under the iterated blow-up morphism 
$\Delta_t^{[1]}  \rightarrow  \Delta_t^{[t+1]}$, 
to the incidence variety $(I_t)^i\subset \Delta_t^{[t+1]}$,
given in (\ref{eq-incidence-variety-I-t-i}). 
\item
\label{lemma-item-pushforward-of-ideal-sheaf-is-ideal-sheaf-of-incidence-div}
The sheaf theoretic pushforward of the line-bundle
$\StructureSheaf{\Delta^{[1]}_t}(-\sum_{i=0}^{t-1}
\Delta_i^{[1]}\cap\Delta^{[1]}_t)$, from $\Delta^{[1]}_t$ to 
$\Delta^{[t+1]}_t$, 
is isomorphic to the ideal sheaf $\StructureSheaf{}(-I)$ 
of the incidence divisor $I$ in the
fiber product of the two dual grassmannian bundles. The higher direct images
vanish. 
\end{enumerate}
\end{new-lemma}

Note: We postpone the statement of additional properties of the 
correspondences to section 
\ref{sec-more-on-the-structure-of-the-correspondences} (lemmas
\ref{lemma-relationship-between-two-pullcacks-of-U_t-to-Delta-t} and
\ref{lemma-two-degeneracy-loci-on-Delta-t-are-identical}). 
We will prove these lemmas only in the special case of moduli spaces, 
but lemma \ref{lemma-two-degeneracy-loci-on-Delta-t-are-identical}
actually holds in general, and
lemma \ref{lemma-relationship-between-two-pullcacks-of-U_t-to-Delta-t}
holds, whenever the projective bundle (\ref{eq-P-W-k+1}) 
comes from a vector bundle. 

\noindent
{\bf Proof of lemma \ref{lemma-structure-of-correspondence}:}
Part \ref{lemma-item-components-of-correspondence-are-smooth}) of the lemma 
follows from the definition 
(\ref{eq-top-iterate-of-component-of-correspondence}) of  
$\Delta_t^{[1]}$ as the iterared blow-up 
of the fiber product $\Delta_t^{[t+1]}$ of dual Grassmannian bundles.

Part \ref{lemma-item-top-correspondence-has-normal-crossing}) 
follows from the definition of the correspondence $\Z^{[1]}$,
as the fiber of $B^{[1]}\Y$ over $0\in B$. $\Z^{[1]}$ is thus isomorphic 
to the fiber of $B^{[1]}\X$ over $0\in B$. 
The normal crossing property, of the special fibers 
of $B^{[t]}\X$, is proven by descending induction on $t$. 
It is clear for $t=\mu+1$. 
The blow-up $B^{[t]}\X\rightarrow B^{[t+1]}\X$
is centered along the smooth subvariety 
$B^{[t+1]}X^t$ of the proper transform of 
the fiber of $B^{[t+1]}\X$ over $0\in B$. 
It suffices to prove, that the intersection of $B^{[t+1]}X^t$,
with the union of the exceptional divisors 
$B^{[t+1]}\X^i$, $i>t$, is a divisor with normal crossing in 
$B^{[t+1]}X^t$ (Lemma \ref{lemma-criterion-for-normal-crossing}). 
The intersection of 
$B^{[t+1]}X^t$ with each of the exceptional divisors 
$B^{[t+1]}\X^i$, $i>t$, is equal to the intersection of 
$B^{[t+1]}X^t$ with $B^{[t+1]}X^i$. Thus, it suffices to prove, that 
the union of all exceptional divisors in 
$B^{[i+1]}X(k)^{i}$ 
is a divisor with normal crossing
in $B^{[i+1]}X(k)^{i}$.
The intersection of 
$B^{[i+1]}X(k)^i$ with $B^{[i+1]}X(k)^j$, $j>i$, 
is the pull back of the exceptional divisor
$B^{[1]}X(k+i)^{j-i}$, by equality
(\ref{eq-conditioned-equality-of-catier-divisor}) in the definition of 
a stratified dualizable collection. 
The pull-back to $B^{[i+1]}X(k)^i$, of a divisor with normal crossing in
$B^{[1]}X(k+i)$, remains a divisor with normal crossing, since the
morphism $B^{[i+1]}X(k)^i\rightarrow B^{[1]}X(k+i)$ is smooth. 
It suffices to prove, that the union of all 
exceptional divisors in $B^{[1]}X(k+i)$ is a divisor with normal crossing.

We prove, that the union of all 
exceptional divisors in $B^{[t]}X(k)$, is a divisor with normal crossing.
The proof is by descending induction on $t+k$. The statement is trivial for 
$t+k=\mu+1$. $B^{[t]}X(k)$ is the blow-up of
$B^{[t+1]}X(k)$ along $B^{[t+1]}X(k)^{t}$. 
The union of exceptional divisors in $B^{[t+1]}X(k)$ is a divisor with
normal crossing, by the induction hypothesis.
The union of exceptional divisors in $B^{[t+1]}X(k)^{t}$ 
is the pullback of those in $B^{[1]}X(k+t)$, which is a divisor with
normal crossing, by the induction hypothesis. 
Lemma \ref{lemma-criterion-for-normal-crossing} implies 
the statement for $B^{[t]}X(k)$. This completes the proof of part
\ref{lemma-item-top-correspondence-has-normal-crossing}.

\medskip
Proof of 
Part \ref{lemma-item-fibers-of-beta-are-embedded-in-a-single-component}) 
It suffices to prove, that the fiber of $\beta$ over $z$ is 
a subscheme of $\Z^{[1]}_t$. 
Consider the following commutative diagram of morphisms
\[
\begin{array}{ccccccc}
\Z^{[1]} & \hookrightarrow & 
B^{[1]}\Y & \hookrightarrow &(B^{[1]}\X)\times \X' & \rightarrow &
B^{[1]}\X 
\\
\downarrow & & \downarrow & & \downarrow & & \downarrow
\\
\Z & \hookrightarrow &
\Y & \hookrightarrow & \X\times \X' & \rightarrow & \X.
\end{array}
\]
The top middle morphism is the embedding of
$B^{[1]}\Y$ as the graph of the composition
\[
B^{[1]}\X \ \LongIsomRightArrow \ B^{[1]}\X' \ \rightarrow \X'.
\]
The right hand square is cartesian. Hence, each fiber of
$\Z^{[1]}\rightarrow \Z$ embedds in the corresponding
fiber of $B^{[1]}\X\rightarrow \X$. 
Now, $\Z^{[1]}_t\setminus\cup_{i=t+1}^\mu\Z^{[1]}_i$ is isomorphic
to the inverse image of $\X^t\setminus\cup_{i=t+1}^\mu\X^i$ in $B^{[1]}\X$.

\medskip
Proof of part \ref{lemma-item-pull-back-of-BN-divisors-to-correspondence})
The intersection $\Delta_t^{[1]}\cap \Delta_i^{[1]}$ is
identified with $B^{[1]}\Y^t\cap B^{[1]}\Y^i$, by Lemma 
\ref{lemma-two-descriptions-of-the-components-of-the-correspondence}. 
The isomorphism $B^{[1]}\Y^t \ \cong B^{[1]}\X^t$ translates it to the 
intersection $B^{[1]}\X^t\cap B^{[1]}\X^i$. 
The stratification of $\X$ is defined in terms of 
the stratification of its special fiber $X$. 
Thus, $B^{[t+1]}\X^t$ is isomorphic to $B^{[t+1]}X^t$. 
The exceptional divisor 
$B^{[t+1]}\X^t\cap B^{[t+1]}\X^i$, in $B^{[t+1]}\X^t$,
is isomorphic to the divisor 
$B^{[t+1]}X^t\cap B^{[t+1]}X^i$. The pull-back (total-transform) of 
$B^{[t+1]}\X^t\cap B^{[t+1]}\X^i$ 
in $B^{[1]}\X^t$ is $B^{[1]}\X^t\cap B^{[1]}\X^i$. 

The proof of part
\ref{lemma-item-pullback-of-incidence-variety-is-intersection-of-components}
is identical to that of lemma
\ref{lemma-two-degeneracy-loci-on-Delta-t-are-identical} 
part \ref{lemma-item-generic-rank-of-alpha-t}.
The proof of part
\ref{lemma-item-pushforward-of-ideal-sheaf-is-ideal-sheaf-of-incidence-div} 
is identical to that of lemma
\ref{lemma-two-degeneracy-loci-on-Delta-t-are-identical} 
part 
\ref{lemma-item-pushforward-of-ideal-sheaf-is-ideal-sheaf-of-incidence-div+2}.
\EndProof

\begin{new-lemma}
\label{lemma-criterion-for-normal-crossing}
Let $D$ be a divisor with normal crossing in a smooth variety $X$,
$D_0$ one of the irreducible components of $D$, and 
$D'$ the union of the other components. 
Assume, that $Y$ is a smooth subvariety of $D_0$, such that the intersection 
$Y\cap D'$ is a divisor with normal crossing in $Y$. Let $\tilde{X}$
be the blow-up of $X$ along $Y$ and $E$ the exceptional divisor. 
Then the union of $E$, with the proper transform 
of $D$, is a divisor with normal crossing in $\tilde{X}$. 
In particular, the proper transform of $D'$, as well as its
union with $E$, are divisors with normal crossings.
\end{new-lemma}
\subsection{Deformation to the normal bundle of the smallest stratum}
\label{sec-deformation-to-the-cone}

The main result of this section is Proposition
\ref{prop-rational-singularities-of-correspondence}. 
Let $\beta : \Z^{[1]}\rightarrow \Z$ be the natural morphism. 

\begin{prop}
\label{prop-rational-singularities-of-correspondence}
\begin{enumerate}
\item
\label{prop-item-K-theoretic-push-forward-is-O-Z}
The direct image  $\beta_*\StructureSheaf{\Z^{[1]}}$ is 
$\StructureSheaf{\Z}$ and 
the higher direct image sheaves all vanish.
\[
R^i_\beta\StructureSheaf{\Z^{[1]}} \ = \ 
\left\{
\begin{array}{cl}
\StructureSheaf{\Z} & 
\mbox{if} \ i=0,
\\
0 & \mbox{if} \ i>0.
\end{array}
\right.
\]
\item
\label{prop-item-Z-is-a-fiber}
$\Z$ has the scheme structure of the fiber of $\Y$ over $0$.
\end{enumerate}
\end{prop}

In this section, we reduce the proof of the proposition to the 
prototypical case $X=T^*G(r,H)$. The latter case is proven 
in section \ref{sec-rational-singularities}. 
The reduction to the special case is obtained using the 
construction of the deformation of a variety to the normal cone of a 
subvariety \cite{fulton}. This construction 
deforms the correspondence $\Z$ to (a relative version of) 
the correspondence $Z$ in the prototypical example. 

\begin{rem}
\label{rem-suffices-to-prove-at-Z-mu}
{\rm
The statements of the proposition are local in $\Y$ (or in $\Z$). 
The set theoretic version of part \ref{prop-item-Z-is-a-fiber} is clear. 
Given a point $z$ of $\Z$, let $\Z_t$ be the component, of maximal index,
containing $z$. Then $z$ is a smooth point of $\Z_t$. 
We prove the proposition in a neighborhood of an arbitrary point in the locus 
$\Z_\mu\subset \Y$, i.e., when $t=\mu$. 
The proof for other values of $t$ is identical and is omitted.
}
\end{rem}

Let $\X(t)\rightarrow B(t)$ be the deformation of $X(t)$, 
of the type considered in
section \ref{sec-resolution-of-correspondences}.
We do not assume any relation among the bases $B(t)$, 
for different values of $t$. 
We construct first the deformation of $\X(t)$ to the normal bundle 
of $X(t)^{\mu-t}$. 
Let $\overline{\M}(t)$ be the blow-up of $\X(t)\times \PP^1$ along 
the subscheme $X(t)^{\mu-t}\times \{\infty\}$.
Denote the proper transform of 
$\X(t)\times \{\infty\}$ by $\widetilde{\X}(t)$.
Note, that the proper transform of $X(t)^{\mu-t}\times \PP^1$
in $\overline{\M}(t)$ is isomorphic to $X(t)^{\mu-t}\times \PP^1$
and is disjoint from $\widetilde{\X}(t)$.
Set 
\[
\M(t) \ \ := \ \ \overline{\M}(t) \ \setminus \ \widetilde{\X}(t).
\]
We have the diagram
\begin{equation}
\label{eq-embedding-of-X-mu-times-P-1}
\begin{array}{ccccc}
X(t)^{\mu-t}\times \PP^1 & & \LongRightArrowOf{\epsilon} & & \M(t)
\\
& pr \ \searrow \ \hspace{2ex}
& & \hspace{1ex} \ \swarrow \ \rho
\\
& & \PP^1
\end{array}
\end{equation}
The fiber $\rho^{-1}(\infty)$ is the normal bundle $N_{X(t)^{\mu-t}/\X(t)}$,
which is isomorphic to an extension
\[
0 \rightarrow T^*_{f_{t,\mu-t}} \rightarrow E(t) \rightarrow 
(N_{X(t)/\X(t)}\restricted{)}{X(t)^{\mu-t}} \rightarrow 0.
\]
The bundle $N_{X(t)/\X(t)}$ is trivial and the extension restricts to each 
grassmannian fiber
of $f_{t,\mu-t}$, as the non-trivial extension 
(\ref{eq-non-trivial-extension-is-deformation-of-cotangent-bundle}),
by Condition \ref{cond-extension-class-is-non-trivial}.
Moreover, the morphism $\epsilon$ embedds $X(t)^{\mu-t}\times \{\infty\}$
 as the zero section of $N_{X(t)^{\mu-t}/\X(t)}$
(see \cite{fulton} section 5.1).

We have a similar deformation of $X(t)$ to the normal bundle 
of $X(t)^{\mu-t}$ in $X(t)$. Set $\overline{M}(t)$ to be the blow-up
of $X(t)\times \PP^1$ along $X(t)^{\mu-t}\times \{\infty\}$ and
let $M(t)$ be the complement of the proper transform of 
$X(t)\times \{\infty\}$
\[
\begin{array}{ccccc}
X(t)^{\mu-t}\times \PP^1 & & \LongRightArrowOf{e} & & M(t)
\\
& pr \ \searrow \ \hspace{2ex}
& & \hspace{1ex} \ \swarrow 
\\
& & \PP^1
\end{array}
\]
Observe, that $M(t)$ is the fiber of $\M(t)$ over $0\in B(t)$. 

Let $\M(k)^t$, $0\leq k\leq \mu$ and $0<t\leq \mu-k$,
be the strict transform of $\X(k)^t\times \PP^1$ in $\M(k)$.
Define $M(k)^t$ similarly. Note, in fact, that $\M(k)^t=M(k)^t$.
Hence, the equality $B^{[t+1]}\M(k)^t=B^{[t+1]}M(k)^t$. 
Recall, that the relative cotangent bundle
$(T^*f_{k,\mu-k})$ restricts to each grassmannian fiber 
of $f_{k,\mu-k}$, as a homomorphism bundle. 
Let $(T^*f_{k,\mu-k})^t$ be the determinantal subvariety of 
homomorphisms with cokernel of dimension $\geq t$. 
One of the conditions, in the definition of a dualizable collection,
implies the scheme-theoretic equality
\begin{equation}
\label{eq-stratification-deforms-to-the-normal-cone}
\M(k)^t \ \cap \ \rho^{-1}(\infty) \ \ = \ \ (T^*f_{k,\mu-k})^t
\end{equation}
(Condition 2.5 in Definition 2.3 in \cite{markman-reflections},
recalled in equation 
(\ref{cond-induced-stratification-on-exceptional-divisor-is-determinantal}) 
above). 
Note, that $\M(k)^{\mu-k}$ is the image of $X(k)^{\mu-k}\times \PP^1$ via
(\ref{eq-embedding-of-X-mu-times-P-1}). 
Moreover, the varieties $M(k)^t$, for different values of $k$, 
are related as follows. 
The grassmannian bundle $f_{k,t}$, given in 
(\ref{eq-f-k-t}), extends to a grassmannian bundle 
\[
f_{k,t} \ : \ B^{[t+1]}M(k)^t \ \ \longrightarrow \ \ 
B^{[1]}M(k+t). 
\]
Over $\infty$, $f_{k,t}$ restricts to 
$B^{[t+1]}(T^*f_{k,\mu-k})^t$ as the relative version of the
Grassmannian bundle (\ref{eq-prototypical-f-k-i}).
Both $B^{[1]}M(k+t)$ and $B^{[t+1]}M(k)^t$ are smooth varieties.

We construct next a triple $\Z\Z\subset \Gamma \subset \M\times \M'$, 
which is a deformation, over $\PP^1$, 
of $\Z\subset \Y \subset \X\times \X'$ to the prototypical case. 
Similarly, we construct a deformation 
$\Z\Z^{[1]}\subset \Gamma^{[1]}\subset B^{[1]}\M\times B^{[1]}\M'$
of $\Z^{[1]}\subset B^{[1]}\Y\subset B^{[1]}\X$.

Denote by $\M'(k)$ 
the analogue of $\M(k)$, 
for the varieties $\X'(k)$ and $X'(k)$ dual to $\X(k)$ and $X(k)$. 
Let $B^{[1]}\X(k) \rightarrow B^{[1]}\X'(k)$ be the isomorphism
(\ref{eq-top-iterated-blow-up-of-Y}).
It is easy to check, that the isomorphism extends to an isomorphism
\[
B^{[1]}\M(k) \ \ \rightarrow \ \ B^{[1]}\M'(k).
\]
Let $\Gamma^{[1]}(k)$ be its graph in $B^{[1]}\M(k)\times B^{[1]}\M'(k)$. 
Set $\Z\Z^{[1]}(k)$ to be the fiber of 
$\Gamma^{[1]}(k)$ over $0\in B$.
Define $\Gamma(k) \subset \M(k)\times \M'(k)$
to be the closure of the graph of the birational isomorphism. 
Given $x\in \PP^1$, 
denote by $\Gamma_x(k)$ the fiber of $\Gamma(k)$ over $x$. 
If $x\neq \infty$, then $\Gamma_x(k)$ is precisely
$\Y(k)$, given in (\ref{eq-Y}). 
The fiber $\Gamma_\infty(k)$ is a locally trivial 
family, over $X(\mu)$, whose fiber is the variety $\Y$, 
for the case where $X=T^*G(\mu,2\mu+n-1)$ is the
cotangent bundle of a fiber of $f_{0,\mu}:X^\mu\rightarrow X(\mu)$. 
This description of $\Gamma_\infty(k)$ follows from the equality 
(\ref{eq-stratification-deforms-to-the-normal-cone}).

Set $\Gamma:=\Gamma(0)$. 
Let $\Z\Z$ be the closure of $\Z\times {\Bbb A}^1$ in $\Gamma$. 
Denote by $\zeta:\Z\Z\rightarrow \PP^1$ the natural morphism. 
Note, that set theoretically, $\Z\Z$ is equal to the fiber of $\Gamma$ over 
$0$ in $B$. 
Denote its components by $\Z\Z_t$, $0\leq t \leq \mu$. 
Note, that $\Z\Z_\mu$ is the image of the embedding
\begin{equation}
\label{eq-embedding-of-Z-mu-times-P-1}
\Z_\mu\times \PP^1 \ \ \hookrightarrow \ \ \M\times \M',
\end{equation}
corresponding to the embeddings $\epsilon$ in 
(\ref{eq-embedding-of-X-mu-times-P-1}) and its
analogue $\epsilon'$ for $\M'$. 
Given $x\in \PP^1$, denote by $\zeta^{-1}(x)$ the fiber of $\Z\Z$ over $x$. 
If $x\neq \infty$, then $\zeta^{-1}(x)$ is precisely
$\Z$. The fiber $\zeta^{-1}(\infty)$ is a locally trivial 
family, over $X(\mu)$, whose fiber is the variety $\Z$, 
for the case where $X=T^*G(\mu,2\mu+n-1)$ is 
the cotangent bundle of a fiber of $f_{0,\mu}:X^\mu\rightarrow X(\mu)$.

\medskip
Let us prove Proposition \ref{prop-rational-singularities-of-correspondence}, 
under the assumption, that it holds 
when $X$ is the cotangent bundle of the Grassmannian. 
We prove the proposition in a neighborhood of a point in
$\Z_\mu$ (see Remark \ref{rem-suffices-to-prove-at-Z-mu}). 

{\bf Proof of part \ref{prop-item-K-theoretic-push-forward-is-O-Z}
of proposition 
\ref{prop-rational-singularities-of-correspondence}:}
The morphisms $\zeta:\Z\Z\rightarrow \PP^1$ and
$\beta\circ\zeta : \Z\Z^{[1]}\rightarrow \PP^1$
are flat, since $\Z\Z$ and $\Z\Z^{[1]}$ are reduced and 
each of their irreducible components surjects onto $\PP^1$. 
We may apply Lemma \ref{lemma-on-support-of-pushforward}, with $X=\Z\Z^{[1]}$,
$Y=\Z\Z$, $f=\beta$, and $C=\PP^1$.
The prototypical case of Part 
\ref{prop-item-K-theoretic-push-forward-is-O-Z} of proposition 
\ref{prop-rational-singularities-of-correspondence}
and Lemma \ref{lemma-on-support-of-pushforward} imply, that the support
$\Sigma\Sigma$ of $R^i\beta_*\StructureSheaf{\Z\Z^{[1]}}$ 
does not intersect the fiber $\zeta^{-1}(\infty)$. 

Let $\Sigma$ be the support of $R^i\beta_*\StructureSheaf{\Z^{[1]}}$. 
The open subset $\Sigma\Sigma^0:=\Sigma\Sigma\cap [\Z\times {\Bbb A}^1]$ 
of $\Sigma\Sigma$ is equal to $\Sigma\times {\Bbb A}^1$. 
Assume, that the intersection $\Sigma\cap \Z_\mu$ is not empty, 
and let $z$ be a point in it. Then $\{z\}\times {\Bbb A}^1$ 
is contained in $\Sigma\Sigma$. Since the latter is closed, then 
it contains $\{z\}\times \PP^1$. This contradicts the emptyness of
$\Sigma\Sigma\cap \zeta^{-1}(\infty)$. 
Hence, $\Sigma\cap \Z_\mu$ is empty.

{\bf Proof of part
\ref{prop-item-Z-is-a-fiber} of proposition 
\ref{prop-rational-singularities-of-correspondence}:}
Let $\Z\Z'$ be the fiber of $\Gamma$ over $0\in B$ and
$\zeta':\Z\Z'\rightarrow \PP^1$ the natural morphism. 
Then $\Z\Z$ is the reduced induced subscheme of $\Z\Z'$.
We need to prove, that $\Z\Z'$ is reduced along $\Z\Z_\mu$.

The fiber of $\Z\Z'$ over $\infty\in \PP^1$ is the fiber 
of $\Gamma$ over 
$(0,\infty)\in B\times \PP^1$, which is also the fiber over $0$ of 
$\Gamma_\infty$. The fiber over $0$ of 
$\Gamma_\infty$ is reduced, by the prototypical case. 
Hence, the fibers of $\Z\Z$ and $\Z\Z'$ over $\infty\in \PP^1$
are the same
\begin{equation}
\label{eq-equality-of-fibers-over-infinity}
\zeta^{-1}(\infty) \ \ = \ \ (\zeta')^{-1}(\infty).
\end{equation}

Let $N$ be the sheaf of nilpotent ideals, in the natural short exact sequence 
\[
0\rightarrow N \rightarrow \StructureSheaf{\Z\Z'} 
\RightArrowOf{\rho} 
\StructureSheaf{\Z\Z}\rightarrow 0.
\]
Since $\Z\Z$ is flat over $\PP^1$, then the sheaf 
$\SheafTor^1(\StructureSheaf{\Z\Z},\StructureSheaf{\zeta^{-1}(\infty)})$ 
vanishes and
the restriction of the 
above sequence to the fiber $\zeta^{-1}(\infty)$ remains
exact
\[
0\rightarrow \restricted{N}{\zeta^{-1}\infty} \rightarrow 
(\StructureSheaf{\Z\Z'}\restricted{)}{\zeta^{-1}\infty} 
\RightArrowOf{\rho_\infty} 
(\StructureSheaf{\Z\Z}\restricted{)}{\zeta^{-1}\infty}\rightarrow 0.
\]
Now $\rho_\infty$ is an isomorphism, by the equality
(\ref{eq-equality-of-fibers-over-infinity}).
Hence, the support $\Sigma\Sigma$ of $N$ is disjoint from
the the fiber $\zeta^{-1}(\infty)$.
The rest of the proof is identical to that of part 
\ref{prop-item-K-theoretic-push-forward-is-O-Z}.
\EndProof

\begin{new-lemma}
\label{lemma-on-support-of-pushforward}
Let  $C$ be a smooth curve, $X$ and $Y$ reduced schemes,
$\pi_1:X\rightarrow C$ and $\pi_2:Y\rightarrow C$ flat morphisms, 
and $f:X\rightarrow Y$ a surjective projective $C$-morphism, 
i.e., $\pi_1=\pi_2\circ f$.
Let $\infty\in C$ be a closed point 
and $f_\infty:X_\infty\rightarrow Y_\infty$
the restriction of $f$ to the fibers of $\pi_1$ and $\pi_2$ over $\infty$. 
Assume, that $f_{\infty,*}\StructureSheaf{X_\infty}=\StructureSheaf{Y_\infty}$
and $R^i f_{\infty,*}\StructureSheaf{X_\infty}$ vanishes, for $i>0$.
Then the higher direct image sheaves
$R^i f_*\StructureSheaf{X}$, $i>0$, all have supports,
which are disjoint from the fiber $Y_\infty$. The same is true for
the support of the quotient $(f_*\StructureSheaf{X})/\StructureSheaf{Y}$.
\end{new-lemma}

\noindent
{\bf Proof:}
We need to compare the sheaves
$(R^if_*\StructureSheaf{X})\otimes \pi_2^*\StructureSheaf{\infty}$ and 
$R^if_*(\pi_1^*\StructureSheaf{\infty})$.
Consider the short exact sequence of sheaves over $C$
\begin{equation}
\label{eq-exact-seq-sky-scraper-sheaf-at-infinity}
0\rightarrow \StructureSheaf{C}(-\infty)\RightArrowOf{\iota} 
\StructureSheaf{C}
\rightarrow \StructureSheaf{\infty}\rightarrow 0.
\end{equation}
It yields the commutative diagram of sheaves over $Y$, where the top row
is the right exact sequence, obtained by tesoring 
(\ref{eq-exact-seq-sky-scraper-sheaf-at-infinity}) with 
$R^if_*\StructureSheaf{X}$.
The  pullback of the sequence 
(\ref{eq-exact-seq-sky-scraper-sheaf-at-infinity}) to $X$ is exact,
since $\pi_1$ is flat.
The bottom row is part of the long exact sequence
of higher direct images, obtained from the 
pullback of (\ref{eq-exact-seq-sky-scraper-sheaf-at-infinity}) to $X$.
\[
\begin{array}{ccccccc}
(R^if_*\StructureSheaf{X})\otimes \pi_2^*\StructureSheaf{C}(-\infty)
& \rightarrow & 
R^if_*\StructureSheaf{X} 
& \rightarrow & 
(R^if_*\StructureSheaf{X})\otimes \pi_2^*\StructureSheaf{\infty}
& \rightarrow & 0
\\
\cong \ \downarrow \ \hspace{1ex} & & 
= \ \downarrow \ \hspace{1ex} & & 
\hspace{1ex} \ \downarrow \ \eta
\\
R^if_*(\pi_1^*\StructureSheaf{C}(-\infty))
& \rightarrow & 
R^if_*\StructureSheaf{X} 
& \RightArrowOf{j} & 
R^if_*(\pi_1^*\StructureSheaf{\infty}) 
& \rightarrow & 
R^{i+1}f_*(\pi_1^*\StructureSheaf{C}(-\infty))
\end{array}
\]
The left vertical homomorphism is an isomorphism, by the projection
formula. It follows, that the right vertical homomorphism $\eta$
is injective onto the image of $j$.

Assume, $i>0$. Then $R^if_*(\pi_1^*\StructureSheaf{\infty})$ vanishes, by
assumption. Since $\eta$ is injective, then 
$(R^if_*\StructureSheaf{X})\otimes \pi_2^*\StructureSheaf{\infty}$
vanishes as well. 

The $i=1$ case implies, that 
$R^1f_*(\pi_1^*\StructureSheaf{C}(-\infty))$
vanishes in a neighborhood of the fiber $Y_\infty$. Hence, the
homomorphism
$j : f_*\StructureSheaf{X}\rightarrow f_{*,\infty}\StructureSheaf{X_\infty}$
is surjective. Using the assumed equality
$f_{*,\infty}\StructureSheaf{X_\infty}=\StructureSheaf{Y_\infty}$, 
we get the bottom short exact sequence in the following
commutative diagram.
\[
\begin{array}{ccccccc}
0 \rightarrow & 
\pi_2^*\StructureSheaf{C}(-\infty)  & \rightarrow &
\StructureSheaf{Y}   & \rightarrow &
\StructureSheaf{Y_\infty}
& \rightarrow 0
\\
& \downarrow & & \downarrow & & \hspace{1ex} \ \downarrow \ = 
\\
0 \rightarrow & 
f_*\pi_1^*\StructureSheaf{C}(-\infty) & \rightarrow &
f_*\StructureSheaf{X} & \RightArrowOf{j} &
\StructureSheaf{Y_\infty}
& \rightarrow 0
\end{array}
\]
The surjectivity of $f$ implies, that the left and middle vertical 
homomorphisms are injective. The snake lemma imlies, that the homomorphism
$\iota$ in (\ref{eq-exact-seq-sky-scraper-sheaf-at-infinity}) induces
the isomorphism
\[
\left[
(f_*\StructureSheaf{X})/\StructureSheaf{Y}
\right]
\otimes \pi_2^*\StructureSheaf{C}(-\infty)
\ \ \ \cong \ \ \ 
(f_*\StructureSheaf{X})/\StructureSheaf{Y}.
\]
Thus, the support of $(f_*\StructureSheaf{X})/\StructureSheaf{Y}$
does not intersect $Y_\infty$.
\EndProof

\subsection{Proof of the prototypical case of Proposition
\ref{prop-rational-singularities-of-correspondence}}
\label{sec-rational-singularities}

We prove in this section the prototypical case of 
proposition \ref{prop-rational-singularities-of-correspondence}, 
when $X=T^*G(r,H)$. This completes the proof of  proposition
\ref{prop-rational-singularities-of-correspondence}. 
The following Lemma introduces a complex, which is quasi-isomorphic
to the structure sheaf $\StructureSheaf{\Z}$. 

\begin{new-lemma}
\label{lemma-complex-qi-to-structure-sheaf}
Let $Z_0$, \dots, $Z_r$ be closed subschemes of a scheme $X$. Set 
$Z_{t_0,t_1,\dots,t_\ell}:=\cap_{i=0}^\ell Z_{t_i}$. 
Let $(E_\bullet,d_\bullet)$ be the complex
\begin{equation}
\label{eq-complex-qi-to-structuresheaf}
\oplus_{t_0=0}^r\StructureSheaf{Z_{t_0}} \rightarrow \cdots \rightarrow
\oplus_{t_0<t_1<\cdots<t_\ell}\StructureSheaf{Z_{t_0,t_1,\dots,t_\ell}}
\RightArrowOf{d_\ell} \cdots \rightarrow 
\StructureSheaf{Z_{0,1,\dots,r}},
\end{equation}
where the differential $d_\ell$ is given by
\[
d_\ell(x)_{t_0,t_1,\dots,t_{\ell+1}} \ \ \ = \ \ \ 
\sum_{j=0}^{\ell+1}(-1)^j
(x_{t_0,\dots,\hat{t}_j,\dots,t_{\ell+1}}
\restricted{)}{Z_{t_0,\dots,t_{\ell+1}}}.
\]
Then the sheaf cohomologies $\H^i(E_\bullet)$ vanish,
for $i>0$, and 
$\H^0(E_\bullet)=\StructureSheaf{\cup_{i=0}^r Z_i}$.
\end{new-lemma}

\noindent
{\bf Proof:}
The proof is by induction on $r$. The case $r=0$ is clear.
Let $(Q_\bullet,\bar{d}_\bullet)$ be the analogous complex corresponding
to the subvarieties $Z_0$, \dots, $Z_{r-1}$. 
There is a surjective complex homomorphism
$\pi:E_\bullet\rightarrow Q_\bullet$, where 
$\pi_\ell:E_\ell\rightarrow Q_\ell$ 
maps the summand $\StructureSheaf{Z_{t_0,\dots,t_\ell}}$
to $0$, if $t_\ell=r$, and restricts as the identity
of the summands common to $E_\ell$ and $Q_\ell$. 
The kernel of $\pi$ is the complex 
$\StructureSheaf{Z_r}\rightarrow F_\bullet$, where
$F_\bullet$ is the complex
\[
\oplus_{t_0=0}^{r-1}\StructureSheaf{Z_{t_0,r}}\rightarrow \cdots
\rightarrow \oplus_{t_0<t_1<\cdots<t_\ell<r}
\StructureSheaf{Z_{t_0,\dots,t_\ell,r}}\rightarrow \cdots
\rightarrow \StructureSheaf{Z_{0,1,\dots,r}}
\]
corresponding to the subvarieties
$Z_0\cap Z_r$, \dots, $Z_{r-1}\cap Z_r$. 

The exact sequence
\[
0 \rightarrow F_\bullet[-1] \rightarrow \ker(\pi)  \rightarrow
\StructureSheaf{Z_r}\rightarrow 0
\]
yields the long exact sequence
\[
0 \rightarrow \H^0(\ker(\pi)) \rightarrow \StructureSheaf{Z_r}
\RightArrowOf{\delta_1} \H^0(F_\bullet) 
\rightarrow \H^1(\ker(\pi)) \rightarrow 0.
\]
The cohomology sheaf $\H^0(F_\bullet)$ is 
$\StructureSheaf{\cup_{i=0}^{r-1}Z_{i,r}}$ and 
$\H^i(F_\bullet)$ vanishes, for $i>0$, by the
induction hypothesis. Furthermore, the connecting
homomorphism $\delta_1$ is the surjective restriction homomorphism.
Thus, $\H^0(\ker(\pi))$ is the ideal sheaf
$\Ideal{\cup_{i=0}^{r-1}Z_{i,r}}$ in $\StructureSheaf{Z_r}$
and $\H^i(\ker(\pi))$ vanishes, for $i>0$. 

The induction hypothesis implies also that $\H^i(Q_\bullet)$ vanishes, 
for $i>0$, and $\H^0(Q_\bullet)=\StructureSheaf{\cup_{i=0}^{r-1}Z_i}$.
The Lemma now follows from the long exact sequence, associated to
the short exact sequence
\[
0\rightarrow \ker(\pi) \rightarrow E_\bullet 
\RightArrowOf{\pi}
Q_\bullet \rightarrow 0.
\]
\EndProof

\medskip
Given an ascending sequence $t_1< t_2 < \cdots < t_\ell$, 
let $\Z^{[1]}_{t_1\cap \cdots \cap t_\ell}$ be the intersection 
$\cap_{i=1}^\ell \Z^{[1]}_{t_i}$ and denote the morphism
$\Z^{[1]}_{t_1\cap \cdots \cap t_\ell}
\rightarrow \Z_{t_1\cap \cdots \cap t_\ell}$ by $\beta$. 
We set $\Z_{t_1,t_\ell}:=\Z_{t_1}\cap\Z_{t_\ell}$.

\begin{new-lemma}
\label{lemma-K-theoretic-push-forward-from-Z-i-j}
The subschemes $\Z_{t_1,t_\ell}$ and $\Z_{t_1\cap \cdots \cap t_\ell}$ are 
equal.
The direct image of $\StructureSheaf{\Z^{[1]}_{t_1\cap \cdots \cap t_\ell}}$,
via $\beta:\Z^{[1]}_{t_1\cap \cdots \cap t_\ell}
\rightarrow \Z_{t_1,t_\ell}$, 
is $\StructureSheaf{\Z_{t_1,t_\ell}}$ and all the higher direct images vanish.
\[
R^i_\beta\StructureSheaf{\Z^{[1]}_{t_1\cap \cdots \cap t_\ell}} \ = \ 
\left\{
\begin{array}{cl}
\StructureSheaf{\Z_{t_1,t_\ell}} & 
\mbox{if} \ i=0,
\\
0 & \mbox{if} \ i>0.
\end{array}
\right.
\]
\end{new-lemma}

{\bf Proof: 
} 
When $\ell=1$, the statement follows from the fact, 
that the morphism $\beta$ is a resolution of rational singularities
(Proposition 
\ref{prop-rational-singularities-in-the-Grassmannian-case}).

Assume $\ell\geq 2$. 
The scheme-theoretic equality
$
\Z_{t_1}\cap \Z_{t_2} \cap \cdots \cap \Z_{t_\ell}  =  
\Z_{t_1,t_\ell}
$
is a consequence of claim \ref{claim-transversality} 
and the corresponding equality 
(\ref{eq-intersections-of-components-of-the-variety-W}) 
for the components of the variety of circular complexes.
Hence, the morphism $\beta$
from $\Z^{[1]}_{t_1\cap \cdots \cap t_\ell}$ surjects onto $\Z_{t_1,t_\ell}$.

Recall, that the morphism $\beta_t:\Z^{[1]}_t\rightarrow \Z_t$ 
factors through the resolution
$\Delta^{[t+1]}_t\rightarrow \Z_t$
given in (\ref{eq-Delta-t-is-a-fiber-product}). 
$\Delta^{[t+1]}_t$ is a fiber product, of 
dual grassmannian bundles over $B^{[1]}X(t)$. 
A point $x$ in $B^{[1]}X(t)$ determines a 
$(h-2r+2t)$-dimensional sub-quotient $F_x$ of $H$. 
For $t_1<t_2$, we have 
the incidence variety $(I_{t_2})^{t_1}$, 
given in (\ref{eq-incidence-variety-I-t-i}),
in the fiber product $\Delta^{[t_2+1]}_{t_2}$. This incidence variety 
consists of pairs $(W',W'')$
in $G(t_2,F_x)\times G(t_2,F_x^*)$ satisfying the following condition:
\[
\dim(W'\cap (W'')^\perp) \ \ \geq \ \ t_2-t_1. 
\]
The variety $\Z^{[1]}_{t_2}$ is isomorphic to the top iterated blow-up 
$\Delta^{[1]}_{t_2}$ of $\Delta^{[t_2+1]}_{t_2}$ along the proper 
transforms of 
$(I_{t_2})^{t_1}$, starting with $t_1=0$ and ending with $t_1=t_2-1$
(Lemma \ref{lemma-two-descriptions-of-the-components-of-the-correspondence}).
Denote by $B^{[i]}\Delta^{[t_2+1]}_{t_2}$, $0\leq i\leq t_2-1$, 
the intermediate iterated blow-up of $\Delta^{[t_2+1]}_{t_2}$ 
along the proper transforms of $(I_{t_2})^{t_1}$, 
starting with $t_1=0$ and ending with $t_1=i$. 
Let $B^{[i]}(I_{t_2})^{t_1}$ be the proper transform of $(I_{t_2})^{t_1}$.
Then $B^{[t_1-1]}(I_{t_2})^{t_1}$ is smooth. 
The exceptional
divisor $B^{[t_1]}(I_{t_2})^{t_1}$ is a projective bundle over 
$B^{[t_1-1]}(I_{t_2})^{t_1}$. 
The divisor $B^{[t_2-1]}(I_{t_2})^{t_1}$ is $\Z^{[1]}_{t_1\cap t_2}$,
by Lemma \ref{lemma-structure-of-correspondence} part
\ref{lemma-item-pullback-of-incidence-variety-is-intersection-of-components}.
It is an iterated blow-up of $B^{[t_1]}(I_{t_2})^{t_1}$ along
smooth subvarieties. 
We claim, that the morphism 
$\gamma:B^{[t_1-1]}(I_{t_2})^{t_1} \rightarrow \Z_{t_1,t_2}$
is a resolution of singularities. The morphism $\gamma$ is surjective, 
because the morphism 
$B^{[t_2-1]}(I_{t_2})^{t_1}=\Z^{[1]}_{t_1\cap t_2}\rightarrow \Z_{t_1,t_2}$
factors through $B^{[t_1-1]}(I_{t_2})^{t_1}$. 
We proceed to identify the dense open subset, along which $\gamma$ is
an isomorphism. 
Let $A$ be the open subset of $\Delta_{t_2}^{[t_2-1]}$, over 
the complement $B^{[1]}X(t_2)\setminus \cup_{k=1}^{\mu-t_2}B^{[1]}X(t_2)^k$,
of the exceptional divisors in $B^{[1]}X(t_2)$. 
The morphism $\Delta_{t_2}^{[t_2-1]}\rightarrow \Z_{t_2}$ embeds $A$ 
as an open subset of $\Z_{t_2}$
(Lemma \ref{lemma-structure-of-correspondence} part 
\ref{lemma-item-pull-back-of-BN-divisors-to-correspondence}). 
Denote by $\widetilde{A}$ the inverse image of 
$A$ in $B^{[t_1-1]}(I_{t_2})^{t_1}$.
Let $B$ be the complement 
$B^{[t_1-1]}(I_{t_2})^{t_1}\setminus B^{[t_1-1]}(I_{t_2})^{t_1-1}$. 
Then $B$ maps isomorphically onto 
$(I_{t_2})^{t_1}\setminus (I_{t_2})^{t_1-1}$.
Hence, $\widetilde{A}\cap B \cap B^{[t_1-1]}(I_{t_2})^{t_1}$ 
is a Zariski dense open subset of 
$B^{[t_1-1]}(I_{t_2})^{t_1}$, along which the morphism $\gamma$ 
restricts to an isomorphism.

When $\ell=2$, 
we obtained a factorization of the morphism $\beta$ 
as the composition of 

1) the resolution 
\begin{equation}
\label{eq-resolution-of-Z-t-1-t-ell}
B^{[t_1-1]}(I_{t_\ell})^{t_1} \ \longrightarrow \ \Z_{t_1,t_\ell},
\end{equation}

2) the projective bundle 
$
B^{[t_1]}(I_{t_\ell})^{t_1} \ \longrightarrow \ 
B^{[t_1-1]}(I_{t_\ell})^{t_1}, \ \mbox{and}
$

3) a sequence of blow-ups along smooth subvarieties.

\noindent
Lemma \ref{lemma-K-theoretic-push-forward-from-Z-i-j} 
follows from the fact, that 
the variety $\Z_{t_1,t_\ell}$ has rational singularities, by Proposition
\ref{prop-rational-singularities-in-the-Grassmannian-case}.

When $\ell>2$, it suffices to prove, that 
$\beta:\Z^{[1]}_{t_1}\cap \Z^{[1]}_{t_2} \cap \cdots \cap \Z^{[1]}_{t_\ell}
\rightarrow \Z_{t_1,t_\ell}$ factors through the resolution
(\ref{eq-resolution-of-Z-t-1-t-ell}), via a smooth fibration
with rational fibers. 
Lemma \ref{lemma-structure-of-correspondence} part 
\ref{lemma-item-pullback-of-incidence-variety-is-intersection-of-components}
shows, that the intersection 
$\Z^{[1]}_{t_1}\cap \Z^{[1]}_{t_2} \cap \cdots \cap \Z^{[1]}_{t_\ell}$
is the intersection of $\ell-1$ exceptional
divisors, in the top iterated blow-up
$B^{[t_\ell-1]}(\Delta_{t_\ell}^{[t_\ell-1]})$ of the 
incidence stratification of $\Delta_{t_\ell}^{[t_\ell-1]}$. 
Set $s_i:=t_\ell-t_i$. 
Claim \ref{claim-transversality} relates 
the incidence stratum $(I_{t_\ell})^{t_i}$,
in the fiber product of dual
grassmannian bundles $\Delta_{t_\ell}^{[t_\ell-1]}$,
to the determinantal locus $\Hom(U,Q)^{s_i}$ 
(of corank $\geq s_i$), in the homomorphism bundle 
$\Hom(U,Q)$ of two rank $t_\ell$ vector bundles. 
Moreover, the iterated blow-up $B^{[t_1-1]}(I_{t_\ell})^{t_1}$
corresponds, under Claim \ref{claim-transversality},
to the iterated blow-up $B^{[s_1+1]}\Hom(U,Q)^{s_1}$.

Let $\phi$ be the morphism, from the intersection
$\cap_{i=1}^{\ell-1}B^{[1]}\Hom(U,Q)^{s_i}$, onto $\Hom(U,Q)^{s_1}$.
It remains to prove, that  $\phi$ 
factors through $B^{[s_1+1]}\Hom(U,Q)^{s_1}$
via a smooth fibration with rational fibers. 
The latter fact follows from the results of Vainsencher
\cite{vainsencher}, about blowing-up determinantal ideals. 
Let us describe the fiber of $\phi$ explicitly. 
Assume, for simplicity, that $U$ and $Q$ are vector spaces. 
Then $B^{[s_1+1]}\Hom(U,Q)^{s_1}$ is a bundle
over $G(s_1,U)\times G(t_1,Q)$, whose fiber over the pair of subspaces
$(U_1,Q_1)$, is the space of completed linear maps from 
$U/U_1$ to $Q_1$. 
The intersection $\cap_{i=1}^{\ell-1}B^{[1]}\Hom(U,Q)^{s_i}$
is a bundle over
$Flag(s_1, s_2, \dots, s_{\ell-1}, U)\times 
Flag(t_1, t_2, \dots, t_{\ell-1}, Q)$, whose fiber over the pair of flags
$U=U_0\supset U_1 \supset \cdots \supset U_{\ell-1}$ and
$0=Q_0\subset Q_1\subset \cdots \subset Q_{\ell-1}$,
is the product of the space of completed linear maps from 
$U/U_1$ to $Q_1$, with the product of the spaces
of complete collineations from
$U_{i-1}/U_i$ to $Q_i/Q_{i-1}$, for $i\geq 2$. 
The spaces of complete collineations are smooth compactifications, of 
the dense open 
$GL(U_{i-1}/U_i)\times GL(Q_i/Q_{i-1})$-orbit in 
$\PP\Hom(U_{i-1}/U_i,Q_i/Q_{i-1})$, obtained by an iterated blow-up
of the determinantal stratification.
Clearly, the forgetful morphism $\phi$, 
from $\cap_{i=1}^{\ell-1}B^{[1]}\Hom(U,Q)^{s_i}$ to
$B^{[s_1+1]}\Hom(U,Q)^{s_1}$, is a smooth fibration with rational fibers. 
This completes the proof of 
Lemma \ref{lemma-K-theoretic-push-forward-from-Z-i-j}. 
%
\EndProof

%
\medskip
{\bf Proof of part \ref{prop-item-K-theoretic-push-forward-is-O-Z} of 
Proposition \ref{prop-rational-singularities-of-correspondence}.  \/} 
Let $E_\bullet$ be the complex
\[
\oplus_{t=0}^r\StructureSheaf{\Delta_t^{[1]}} 
\rightarrow \cdots
\rightarrow
\oplus_{t_0< t_1 < \cdots < t_\ell}
\StructureSheaf{\Delta^{[1]}_{t_0\cap \cdots \cap t_\ell}} \rightarrow \cdots 
\rightarrow
\StructureSheaf{\Delta^{[1]}_{0\cap 1 \cap \cdots \cap r}},
\]
of coherent sheaves on $\Z^{[1]}$ analogous to the complex 
(\ref{eq-complex-qi-to-structuresheaf}).
The structure sheaf $\StructureSheaf{\Z^{[1]}}$
is the $0$-th sheaf cohomology $\H^0(E_\bullet)$ and 
$\H^i(E_\bullet)$ vanishes, for $i>0$, by
Lemma \ref{lemma-complex-qi-to-structure-sheaf}. 
Similarly, $\StructureSheaf{\Z}$ is $\H^0(F_\bullet)$, where
the complex $F_\bullet$ is the analogue of $E_\bullet$,
with $\Z^{[1]}$ replaced by $\Z$ and
$\Delta^{[1]}_t$ by $\Z_t$. 
The pushforward of 
$\StructureSheaf{\Delta^{[1]}_{t_0\cap \cdots \cap t_\ell}}$
via the morphism
\[
\beta \ : \ \Delta^{[1]}_{t_0\cap \cdots \cap t_\ell} \ \longrightarrow \ 
\Z_{t_0\cap \cdots \cap t_\ell} 
\]
is $\StructureSheaf{\Z_{t_0\cap \cdots \cap t_\ell}}$ 
and the higher direct images vanish, by 
Lemma \ref{lemma-K-theoretic-push-forward-from-Z-i-j}. 
The vanishing of the higher direct images in the Lemma implies, that 
$R^i_{\beta_*}\StructureSheaf{\Z^{[1]}}$ is isomorphic to 
the sheaf cohomology $\H^i(\beta_*E_\bullet)$, while 
the equality of direct images implies the equality 
$\beta_*E_\bullet=F_\bullet$. 
Part \ref{prop-item-K-theoretic-push-forward-is-O-Z} of Proposition 
\ref{prop-rational-singularities-of-correspondence} follows 
from the vanishing of the sheaf cohomologies $\H^i(F_\bullet)$,
for $i>0$, which is proven in Lemma
\ref{lemma-complex-qi-to-structure-sheaf}.

{\bf Proof of part \ref{prop-item-Z-is-a-fiber} of 
Proposition \ref{prop-rational-singularities-of-correspondence}.  \/} 
In Step 1 we find equations for $\Z$, 
using the equations for the variety of circular complexes
\cite{strickland}. 
In Step 2 we exhibit a set of equations 
(\ref{eq-tau-H-is-phi-invariant})-(\ref
{eq-tau-perp-H-is-phi-dual-invariant}), 
defining a scheme $\Y'$ over the affine line, 
whose reduced induced subscheme is $\Y$. 
In Step 3 we show, that the fiber of $\Y'$ over $0$ is
equal to $\Z$. It follows,
that the fiber of $\Y$ is a closed subscheme of $\Z$. 
We already know, that $\Z$ is the reduced induced subscheme of the
fiber of $\Y$. Hence, the two are equal.

{\bf Step 1:} (Equations for $\Z$) 
Let $p:V\rightarrow G(r,H)\times G(r,H^*)$ be the vector bundle
$\Hom(\pi_2^*\tau^*_{G(r,H^*)},\pi_1^*\tau_{G(r,H)})$. 
Denote by $\Psi$ the tautological section of the pullback $p^*V$ 
of $V$ to $V$. Let $\Z'$ be the subscheme of $V$ defined by the two equations
\begin{eqnarray}
\label{eq-g-Psi-is-zero}
p^*g\circ \Psi & = & 0, 
\\
\label{eq-Psi-g-is-zero}
\Psi\circ p^*g & = & 0,
\end{eqnarray}
where the homomorphism 
$g:\pi_1^*\tau_{G(r,H)}\rightarrow \pi_2^*\tau^*_{G(r,H^*)}$ is given in 
(\ref{eq-homomorphism-between-tautological-subundle-of-dual-grassmannians}). 
Claim \ref{claim-transversality} and the results of Strickland imply,
that $\Z'$ is {\em reduced} \cite{strickland}. 
Moreover, the subvariety of $\Z'$, defined by the equations
$\Wedge{k+1}{p^*g}=0$ and 
$\Wedge{r-k+1}{\Psi}=0$, for $0\leq k \leq r$, 
is an irreducible component $\Z'_{k}$ of $\Z'$, which is 
isomorphic to the conormal variety $\Z_k$ of the incidence variety 
$I_{k}\subset G(r,H)\times G(r,H^*)$ given in equation
(\ref{eq-incidence-variety-in-product-of-grassmannians})
(see \cite{strickland} or the description of the ideal of the variety
$W(k,r-k)$ in section \ref{sec-circular-complexes} above). 

The conormal variety $\Z_k$ of $I_k$ is, by definition, a subvariety of
the cotangent bundle $T^*[G(r,H)\times G(r,H^*)]$. 
The existence, of an embedding of $\Z_k$ in $V$, follows from Claim
\ref{claim-transversality}. 
Let us describe explicitly the isomorphism $\Z_k\cong \Z'_k$, relating 
these two embeddings. 
$V$ embeds as a subbundle of the trivial bundle 
$G(r,H)\times G(r,H^*) \times \End(H)$. The correspondence $\Z$ embeds in 
this trivial bundle as well,
because it is the reduced induced subscheme of the fiber product of
$T^*G(r,H)$ and $T^*G(r,H^*)$ over $\End(H)$. 
The images of $\Z'_k$ and $\Z_k$, 
under these embeddings, are easily seen to be equal. 
Indeed, let $(W_1,W_2,\psi:W_2^*\rightarrow W_1)$ be a point of $\Z'_k$.  
Equations (\ref{eq-g-Psi-is-zero}) and (\ref{eq-Psi-g-is-zero})
imply, that $\ker(\psi)$ contains $[W_1+W_2^\perp]/W_2^\perp$ and 
${\rm Im}(\psi)$ is contained in $W_1\cap W_2^\perp$. 
In particular, $\psi$ lifts to $\varphi_1\in\Hom(H/W_1,W_1)$ as well as to
$\varphi_2\in \Hom(H/W_2^\perp,W_2^\perp)$. The quadruple 
$(W_1,W_2,\varphi_1,\varphi_2)$ is a point in the cotangent bundle 
$T^*[G(r,H)\times G(r,H^*)]$, which belongs to the conormal variety 
$\Z_k$ of $I_k$. 

Summarizing, we have obtained equations 
(\ref{eq-g-Psi-is-zero}) and
(\ref{eq-Psi-g-is-zero})
for $\Z$, as a subvariety of $V$,
which itself is a subbundle of $G(r,H)\times G(r,H^*) \times \End(H)$. 

{\bf Step 2:} (Equations defining $\Y$, set theoretically)
We use equation (\ref{eq-E-is-a-subbundle-of-End-H})
in order to exhibit 
$\Y$ as a subvariety of the smooth variety 
\[
A := G(r,H)\times G(r,H^*)\times \End(H) \times \ComplexNumbers.
\] 
Set $\Xi(H):= G(r,H)\times \End(H)\times \ComplexNumbers$ and
$\Xi(H^*):= G(r,H^*)\times \End(H)\times \ComplexNumbers$. 
Consider the following embedding 
$E(H)\hookrightarrow \Xi(H)$
of $E(H)$, as a subbundle of the trivial vector bundle over $G(r,H)$. 
The embedding is determined by two projections. 
The first projection $E(H)\rightarrow \End(H_{G(r,H)})$, 
into the trivial $\End(H)$ bundle over $G(r,H)$,
is the embedding given by
equation (\ref{eq-E-is-a-subbundle-of-End-H}). 
The second projection $\lambda_H:E(H)\rightarrow \ComplexNumbers$
is defined by the equation
$r\cdot \restricted{\varphi}{W}=\lambda_H(W,\varphi)\cdot I_W$,
where $(W,\varphi)$ is a point in $E(H)$, as in equation
(\ref{eq-E-is-a-subbundle-of-End-H}). Similarly, we have 
the function $\lambda_{H^*}:E(H^*)\rightarrow \ComplexNumbers$ 
and the embedding 
$E(H^*)\hookrightarrow \Xi(H^*)$.
The two projections, of both $E(H)$ and
$E(H^*)$ into $\End(H)\times \ComplexNumbers$, are
birational morphisms onto the same image.
(In case $h\neq 2r$, the same is true for the image in $\End(H)$,
but when $h=2r$, the projection to $\End(H)$ has degree two onto the image). 
The two projections are a dual pair of resolutions of their image,
extending the dual pair (\ref{eq-dual-springer-resolutions}), 
and $\Y$ is the graph of the birational isomorphism. Thus, 
$\Y$ is a closed subvariety of the fiber
product $\Y^{fp}:=E(H)\times_{\End(E)\times \ComplexNumbers} E(H^*)$
(and the two are equal, set theoretically). 
Note, that the variety $A$ is the
fiber product of $\Xi(H)$ and $\Xi(H^*)$ over 
$\End(H)\times \ComplexNumbers$. 
The embeddings $E(H)\subset \Xi(H)$ and $E(H^*)\subset \Xi(H^*)$ 
give rise to a closed imersion $\Y^{fp}\hookrightarrow A$.
We conclude, that $\Y$ is a closed subvariety of $A$.

Set theoretically, $\Y$ is the zero locus, of a section
of the vector bundle 
\[
\Hom(\pi_1^*\tau_{G(r,H)},H) \ \oplus \  
\Hom(\pi_2^*\tau^\perp_{G(r,H^*)},H) \ \oplus \  
\Hom(\pi_2^*\tau_{G(r,H^*)},H^*) \ \oplus \  
\Hom(\pi_1^*\tau^\perp_{G(r,H)},H^*)
\]
over $A$.
Above, $\tau^\perp_{G(r,H^*)}$
is the rank $(h-r)$ subbundle of $H$ annihilating $\tau_{G(r,H^*)}$.
Let $\lambda:A\rightarrow \ComplexNumbers$ be the projection. 
Denote by $\Phi:A\rightarrow \End(H_A)$ the 
tautological section of the trivial $\End(H)$ bundle over $A$. 
$\Y$ is set theoretically cut out by the four equations:
\begin{eqnarray}
\label{eq-tau-H-is-phi-invariant}
\restricted{\Phi}{\pi_1^*\tau_{G(r,H)}} 
& = & 
\left(\frac{\lambda}{r}\right) \cdot I_{\pi_1^*\tau_{G(r,H)}},
\\
\label{eq-tau-perp-H-dual-is-phi-invariant}
\restricted{\Phi}{\pi_2^*\tau^\perp_{G(r,H^*)}} 
& = & - 
\left(\frac{\lambda}{h-r}\right) \cdot I_{\pi_2^*\tau^\perp_{G(r,H^*)}},
\\
\label{eq-tau-H-dual-is-phi-dual-invariant}
\restricted{\Phi^*}{\pi_2^*\tau_{G(r,H^*)}} 
& = & 
\left(\frac{\lambda}{r}\right) \cdot I_{\pi_2^*\tau_{G(r,H^*)}},
\\
\label{eq-tau-perp-H-is-phi-dual-invariant}
\restricted{\Phi^*}{\pi_1^*\tau^\perp_{G(r,H)}} 
& = & - 
\left(\frac{\lambda}{h-r}\right) \cdot I_{\pi_1^*\tau^\perp_{G(r,H)}}.
\end{eqnarray}
Denote by $\Y'$ the subscheme defined by these four equations.

{\bf Step 3:} ($\Z$ is the special fiber of $\Y'$)
We claim, that the subscheme $\Z''$ of $\Y'$, defined by the equation 
$\lambda=0$, is equal to $\Z$. It suffices to prove, that 
$\Z''$ is a closed subscheme of $\Z$. 
Equations (\ref{eq-tau-H-is-phi-invariant}), 
(\ref{eq-tau-perp-H-dual-is-phi-invariant}), 
(\ref{eq-tau-H-dual-is-phi-dual-invariant}), 
(\ref{eq-tau-perp-H-is-phi-dual-invariant}), and 
$\lambda=0$, imply the containments
\begin{eqnarray}
\label{eq-kernel-of-Phi-contains-tau-H}
\ker(\restricted{\Phi}{\Z''}) & \supset & \pi_1^*\tau_{G(r,H)},
\\
\label{eq-kernel-of-Phi-contains-tau-perp-H-dual}
\ker(\restricted{\Phi}{\Z''}) & \supset & \pi_2^*\tau^\perp_{G(r,H^*)},
\\
\label{eq-image-of-Phi-is-in-tau-perp-H-dual}
{\rm Im}(\restricted{\Phi}{\Z''}) & \subset & \pi_2^*\tau^\perp_{G(r,H^*)},
\\
\label{eq-image-of-Phi-is-in-tau-H}
{\rm Im}(\restricted{\Phi}{\Z''}) & \subset & \pi_1^*\tau_{G(r,H)}.
\end{eqnarray}
Equations
(\ref{eq-kernel-of-Phi-contains-tau-perp-H-dual})  and 
(\ref{eq-image-of-Phi-is-in-tau-H}) imply, 
that the restriction $\restricted{\Phi}{\Z''}$ 
belongs to the subbundle $p^*V$ and coincides with the restriction of 
the tautological section $\Psi$. 
Equations (\ref{eq-kernel-of-Phi-contains-tau-H}) and 
(\ref{eq-image-of-Phi-is-in-tau-perp-H-dual}) imply the
equations 
(\ref{eq-g-Psi-is-zero}) and 
(\ref{eq-Psi-g-is-zero}) defining $\Z$ in $V$. 
This completes the proof of Proposition 
\ref{prop-rational-singularities-of-correspondence}
in the case $X=T^*G(r,H)$.
\EndProof


\begin{cor}
\label{cor-Z-is-a-fiber-product}
\begin{enumerate}
\item
\label{cor-item-Z-is-a-fiber-product}
Descriptions 
\ref{correspondence-is-a-fiber-product} and 
\ref{correspondence-fiber-products-of-dual-grassmanian-bundles} 
of $\Z$, in Definition \ref{def-Z-for-cotangent-bundles}, are equivalent.
\item
\label{cor-item-Y-is-a-fiber-product}
$\Y$ is the fiber product, over the affine space
$\End(H)\times \ComplexNumbers$, of $E(H)$ and $E(H^*)$. 
Furthermore, $\Y$ is equal to the subscheme defined by the four equations 
(\ref{eq-tau-H-is-phi-invariant})-(\ref
{eq-tau-perp-H-is-phi-dual-invariant}).
\end{enumerate}
\end{cor}

\noindent
{\bf Proof:}
\ref{cor-item-Z-is-a-fiber-product}) 
Let $\Z^{fp}$ be the fiber product of $T^*G(r,H)$ and $T^*G(r,H^*)$
over $\overline{\N}_r$ given in 
(\ref{eq-dual-springer-resolutions}).
We already know, that $\Z$ is the reduced induced subscheme of
$\Z^{fp}$. Hence, it suffices to prove, that 
$\Z^{fp}$ satisfies the equations 
(\ref{eq-g-Psi-is-zero}) and 
(\ref{eq-Psi-g-is-zero}) defining $\Z$ in $V$. 
Now $\Z^{fp}$ clearly satisfies the four equations 
(\ref{eq-kernel-of-Phi-contains-tau-H})-(\ref
{eq-image-of-Phi-is-in-tau-H}),
as each is the pullback of equations satisfied by
$T^*G(r,H)$ or $T^*G(r,H^*)$. 
These four equations,
combined with the equation $\lambda=0$, imply equations 
(\ref{eq-g-Psi-is-zero}) and 
(\ref{eq-Psi-g-is-zero}) (see the proof 
of part \ref{prop-item-Z-is-a-fiber} of 
Proposition \ref{prop-rational-singularities-of-correspondence}).

\ref{cor-item-Y-is-a-fiber-product}) 
Let $\Y^{fp}$ be the fiber product and $\Y'$ the 
subscheme defined by the four equations 
(\ref{eq-tau-H-is-phi-invariant})-(\ref
{eq-tau-perp-H-is-phi-dual-invariant}).
We continue to use the notation, 
introduced in the proof of part \ref{prop-item-Z-is-a-fiber} of 
Proposition \ref{prop-rational-singularities-of-correspondence}.

{\bf Step 1}
The three subsets $\Y$, $\Y'$, and $\Y^{fp}$, of $A$ 
coincide, set theoretically, and so $\Y$, being an
integral scheme, is the reduced induced subscheme of $\Y'$ and $\Y^{fp}$.
Note, that $\Y^{fp}$ is a closed subscheme of $\Y'$. Indeed,
the two equation (\ref{eq-tau-H-is-phi-invariant}) and 
(\ref{eq-tau-perp-H-is-phi-dual-invariant}) are pullback 
of equations in $\Xi(H)$ satisfied by $E(H)$. 
The two equation (\ref{eq-tau-perp-H-dual-is-phi-invariant}) and 
(\ref{eq-tau-H-dual-is-phi-dual-invariant}) are
the pullback of equations in $\Xi(H^*)$ satisfied by 
$E(H^*)$. Hence, it suffices to prove, the $\Y'$ is reduced. 
We have already shown, in the proof 
of part \ref{prop-item-Z-is-a-fiber} of 
Proposition \ref{prop-rational-singularities-of-correspondence},
that the fiber $\Z'':=\Y'\cap (\lambda=0)$ of $\Y'$
is equal to $\Z$ and is hence reduced.

{\bf Step 2}
Both structure sheaves $\StructureSheaf{\Y}$ and $\StructureSheaf{\Y'}$
are sheaves on the topological space $\Y$. 
Let $\rho : \StructureSheaf{\Y'} \rightarrow \StructureSheaf{\Y}$
be the restriction homomorphism. The kernel $N$ of $\rho$ is a sheaf of 
nilpotent ideals. The open subset $\lambda\neq 0$, of the scheme $\Y'$, is 
easily seen to be reduced. Hence, $N$ is set theoretically supported on the
fiber $\Z$. Denote by $\lambda$ both the function on $A$ and its  
restriction to $\Y'$. Consider the commutative diagram 
\[
\begin{array}{ccccccccc}
0 & \rightarrow & N & \rightarrow & \StructureSheaf{\Y'} & 
\LongRightArrowOf{\rho} & \StructureSheaf{\Y} & \rightarrow & 0
\\
& & (\cdot \lambda\restricted{)}{N} \ \downarrow \ \hspace{6ex} & & 
\cdot \lambda \ \downarrow \ \hspace{3ex} & & 
\hspace{4ex} \ \downarrow \ \cdot\rho(\lambda)
\\
0 & \rightarrow & N & \rightarrow & \StructureSheaf{\Y'} & 
\LongRightArrowOf{\rho} & \StructureSheaf{\Y} & \rightarrow & 0
\end{array}
\]
of multiplication by $\lambda$.
We have shown, that the cokernels, of the middle and right vertical 
homomorphisms, are both isomorphic to $\StructureSheaf{\Z}$. 
Furthermore, the kernel of $\rho(\lambda)$ is trivial, as $\Y$ is integral. 
By the Snake Lemma, the left vertical homomorphism is surjective. Hence, 
the multiplication homomorphism 
$\lambda^k : N \rightarrow N$ is surjective, for all positive powers $k$. 
But $N$ is supported on the set $\lambda=0$. Hence, the latter homomorphisms 
vanish, for all $k$ sufficiently large. It follows, that the sheaf $N$
vanishes and $\Y'$ is reduced. 
\EndProof

\subsection{The induced action on the cohomology}
\label{sec-induced-action-on-cohomology}

Assume now, that $X$ is 
a symplectic irreducible projective variety, or at least,
a compact and hyperk\"{a}ler variety. 
The correspondence $\Z$ in $X\times X'$ induces an isomorphism 
of cohomology rings $\Z_* : H^*(X,\Integers) \rightarrow H^*(X',\Integers)$,
because it is the limit of the family $\Y\rightarrow B$ of isomorphisms 
(Proposition \ref{prop-rational-singularities-of-correspondence}). 
In particular, it preserves the Poincare pairing, and 
$(\Z_*)^*=(\Z_*)^{-1}$. 
The transpose $(\Z_*)^*$ is represented by
the image of $\Z$, via the transposition of the factors of $X\times X'$.
We expect $\Z$ to induce an isomorphism of the Grothendieck K-rings,
preserving the multiplication. This is known when $X$ is the 
cotangent bundle of a Grassmannian \cite{namikawa-stratified}.

The variety $X$ is endowed with a collection of natural
cohomology classes. These include the classes $[X^t]$ Poincare-dual
to the strata. Additional classes are the images, via the Gysin map
$H^*(B^{[t+1]}X^t)\rightarrow H^*(X)$, of the characteristic classes
of the tautological (projective) sub and quotient bundles,
over the Grassmannian bundle 
$B^{[t+1]}X^t\rightarrow B^{[1]}X(t)$, given in (\ref{eq-f-k-t}). 

A somewhat indirect approach, to the action of $\Z_*$ 
on the cohomology and $K_0$-group of $X$, is taken 
in Theorems \ref{thm-class-of-correspondence-in-stratified-elementary-trans} 
and \ref{thm-reflection-sigma-satisfies-main-conj}.
$X$ is a moduli space of coherent sheaves on a K3 surface $S$. 
Generators for $H^*(X,\RationalNumbers)$ are encoded in
the universal sheaf $\E$ over $X\times S$ (see \cite{markman-diagonal}). 
We will calculate the K-theoretic image, 
via $\Z\times \one_S$, of the universal sheaf. 
It would be interesting to use the latter calculation, in order to 
work out the $\Z_*$-image of the 
natural classes listed above. We ilustrate this method next, only for the
classes of the strata. (See also Section 
\ref{sec-relationship-between-two-pullcacks-of-U_t-to-Delta-t}).

Clearly, $\Z_*$ maps the class $[X^t]$, Poincare dual to $X^t$, 
to a multiple of $[(X')^t]$. The coefficient $c_t$ must be $1$ or $-1$, 
since $[X^t]=(\Z_*)^*\circ \Z_*[X^t]=(c_t)^2[X^t]$. 
Note, that the codimension of $X^t$ is $t(t+n-1)$, 
where $n$ is the codimension of $X^1$ in $X$. 
We expect the coefficient to be given by the parity of the codimension, 
\[
\Z_*[X^t] \ \ \ = \ \ \ (-1)^{t(t+n-1)}[(X')^t].
\]
We will prove this equality in two sequences of examples of self-dual 
moduli spaces, where $n=1$ and $n=3$
(Corollaries \ref{cor-action-of-Z-on-individual-strata} 
and \ref{cor-action-of-Z-on-individual-strata-plus-2-case}). 
The proof of Corollary \ref{cor-action-of-Z-on-individual-strata} 
would apply for all $n$, provided the K-theoretic image 
of the universal sheaf, via $\Z\times \one_S$, 
is given by the conjectural equality
(\ref{eq-conjectural-K-theoretic-equality}).

\section{Reflections of moduli spaces by line bundles}
\label{sec-self-dual-stratified-elemetary-transformations}

We recall in section \ref{sec-moduli-spaces-with-stratified-a-1-sing} 
the relationship between 
1) the autoequivalence $\Phi^{\StructureSheaf{S}}$ of $D(S)$, associated
to the reflection with respect to the trivial line bundle 
$\StructureSheaf{S}$, 
and 
2) stratified elementary transformations, 
of moduli spaces of sheaves on the K3 surface $S$. 

Theorem \ref{thm-class-of-correspondence-in-stratified-elementary-trans}, 
which 
is the main result of this section, is stated in section 
\ref{sec-the-class-of-the-correspondence-in-terms-of-E-v}. 
Let $\tau$ be the isometry, of the Mukai lattice, induced by 
the autoequivalence $\Phi^{\StructureSheaf{S}}$. 
The Theorem verifies Theorem \ref{thm-Gamma-v-acts-motivicly}, for 
the isometry $\tau$ and certain moduli spaces 
$\M(v)$, whose Mukai vector is $\tau$-invariant (but may have rank larger 
than $1$). 
Recall the Steinberg correspondence $\Z\subset \M(v)\times \M(v)$, 
associated to the self-dual stratified elementary transformation of $\M(v)$
(Definition \ref{def-Z-for-cotangent-bundles}).
The Theorem expresses, in addition, the cohomology class Poincare-dual to 
$\Z$ in terms of 1) the Chern classes of the universal sheaf $\E_v$ over
$\M(v)$ and 2) the autoequivalence $\Phi^{\StructureSheaf{S}}$.

Theorem \ref{thm-class-of-correspondence-in-stratified-elementary-trans} 
is proven in section 
\ref{sec-proof-of-theorem-on-class-of-correspondence-in-minus-2-reflection}.
The proof boils down to the computation of the 
K-theoretic push-forward of $\tau(\E_v)$ via the correspondence $\Z$. 
Recall, that $\tau(\E_v)$, given in (\ref{eq-tau-of-universal-sheaf}), 
is the K-theoretic representative of the 
relative Fourier-Mukai transform $\Phi^{\StructureSheaf{S}}(\E_v)$ of 
the universal sheaf $\E_v$.
The push-forward of $\tau(\E_v)$ via $\Z$
agrees, in the K-group of $\M(v)\times S$, 
with a universal sheaf over $\M(v)\times S$.

The notation for sections 
\ref{sec-self-dual-stratified-elemetary-transformations} and 
\ref{sec-stratified-elementary-trans-plus-2-vectors} 
is summerized in section \ref{sec-notation-minus-2-vector}.

\subsection{Brill-Noether stratifications}
\label{sec-moduli-spaces-with-stratified-a-1-sing}
We constructed in \cite{markman-reflections} a symplectic birational 
isomorphism between $\M_H(v)$ and $\M_H(\tau(v))$, 
for the reflection $\tau$ in
the $-2$ vector $v_0=(1,0,1)$, 
under the assumption that the determinant line-bundle $\LB$ of $v$ satisfies:

\begin{condition}
\label{condition-minimality}
\begin{enumerate}
\item 
\label{cond-part-minimality}
$\LB$ is an effective cartier divisor with {\em minimal}
degree in the sense, that the subgroup $H\cdot \Pic(S)$ of $H^4(S,\Integers)$
is generated by $H\cdot c_1(\LB)$. 
In particular, all curves in the linear system
$\linsys{\LB}$ are reduced and irreducible. 
\item 
The base locus of 
$\linsys{\LB}$ is either empty or zero-dimensional. 
\item \label{cond-part-generic-curve-is-smooth}
The generic curve in $\linsys{\LB}$ is smooth. 
\item
$H^1(S,\LB)=0$.
\end{enumerate}
\end{condition}


Each moduli space admits a Brill-Noether stratification 
\begin{equation}
\label{eq-brill-noether-stratification}
\M(v) \ = \ \M(v)^0 \ \supset \ \M(v)^1 \ \supset \ \cdots \ \supset \ 
\M(v)^\mu
\end{equation}
If $\chi(v)\geq 0$, 
the stratum $\M(v)^t\setminus\M(v)^{t+1}$ parametrizes stable sheaves $F$
with $h^1(F)=t$. When $\chi(v)<0$, we use $h^0(F)$ instead. We denote by
$\mu$ the length of the stratification. 
The smallest stratum $\M(v)^\mu$ is smooth. We denote by 
$B^{[k]}\M(v)$ the itereted blow-up starting at $\M(v)^\mu$
and proceeding, in decreasing order $t=\mu, \mu\!-\!1, \cdots, k\!+\!1, k$, 
along the strict transforms 
$B^{[t+1]}\M(v)^t$ of $\M(v)^t$ in $B^{[t+1]}\M(v)$. 
Then the top iterated blow-ups $B^{[1]}\M(v)$ and $B^{[1]}\M(\tau(v))$
are isomorphic. 
More precisely, Theorem 3.2 in \cite{markman-reflections}
states, that the birational isomorphism 
between $\M(v)$ and $\M(\tau(v))$ is a 
{\em stratified elementary transformation} (see section
\ref{sec-stratified-elemetary-transformations} or Definition
2.3 in \cite{markman-reflections}). 
Furthermore, the symplectic varieties $X(t)$, in the diagonal entries of 
(\ref{eq-diagram-of-X-v}), are moduli spaces as well. 
If $\chi(v)\geq 0$, then $X(0)=\M(v)$ and $X(t)=\M(v+\vec{t})$, where 
$\vec{t}$ is the Mukai vector $(t,0,t)$ of the trivial rank $t$ vector 
bundle. If $\chi(v)\leq 0$, take $X(t)=\M(v-\vec{t})$ (or
$\M(\sigma\tau[v-\vec{t}])$, whichever has positive rank,
where $\sigma$ is given in (\ref{eq-sigma})). 
The Grassmannian fibration
$B^{[t+1]}\M(v)^t\rightarrow B^{[1]}\M(v+\vec{t})$, mentioned in
(\ref{eq-f-k-t}) for a general stratified elementary transformation,
is described in part 
\ref{prop-item-morphism-is-a-classifying-morphism-of-F-t} 
of Proposition \ref{prop-decomposition-of-complex-over-interated-blow-up}
below. 

The correspondence
\begin{equation}
\label{eq-total-brill-noether-correspondence}
\Z \ \ \subset \ \ \M(v) \times \M(\tau(v)),
\end{equation}
defined in part 
\ref{correspondence-fiber-products-of-dual-grassmanian-bundles}
of definition \ref{def-Z-for-cotangent-bundles}, was shown 
to induce an isomorphism of cohomology rings
\[
H^*(\M(v),\Integers) \ \ \IsomRightArrow \ \ 
H^*(\M(\tau(v)),\Integers) 
\] 
(Theorem 1.2 in \cite{markman-reflections}).
A point in $\Z$ consists of a pair $(E_1,E_2)$ of stable sheaves, such that
the canonical extensions $F_i$, of $E_i$ by the trivial vector bundle
$H^1(E_i)\otimes\StructureSheaf{S}$, are isomorphic $F_1\cong F_2\cong F$. 
The sheaf $F$ is necessarily stable. 
Notice the equalities $h^1(E_1)+\chi(v)=h^1(E_2)$ and
$h^0(F)=h^1(E_1)+h^1(E_2)$. 
These equalities follow from the equalities 
$\chi(E_1)=-\chi(E_2)$, 
$\chi(E_1)+2h^1(E_1)=\chi(F)=\chi(E_2)+2h^1(E_2)$, the vanishing of 
$H^2(E_i)$, and the vanishing of $H^i(F)$, for $i\geq 1$. 
The triple $(F,E_1,E_2)$ is equivalent to the data of two subspaces of
$H^0(F)$ of complementary dimensions $t:=h^1(E_1)$ and 
$t+\chi(v)$. 
The correspondence $\Z$ is reduced, but reducible in general. 
One component $\Z_0$ is the graph of the birational isomorphism. 
When $\chi(v)\geq 0$, the 
other components consist of (the closure of) 
fiber products of Grassmannian bundles 
\begin{equation}
\label{eq-closure-of-fiber-product-of-grassmannian-bundles}
\Z_t := \mbox{closure of} \ \ 
[\M(v)^t\setminus \M(v)^{t+1}] \ \ \times_{\M(v+ \vec{t})}\ \ 
[\M(\tau(v))^{t+\chi(v)}\setminus \M(\tau(v))^{t+\chi(v)+1}].
\end{equation}
When $\chi(v)< 0$, take the fiber product over
$\M(\tau(v)+ \vec{t})$.

\begin{example}
\label{example-geometric-meaning-of-reflections-by-minus-2-vectors}
{\rm
The vector $v=(1,\LB,-1)$ is fixed under the reflection by the 
vector $v_0=(1,0,1)$. $\M(v)$ is isomorphic to the Hilbert scheme $S^{[n]}$
when $\deg{\LB}=2n-4$. Tensoring by $\LB^{-1}$, we may replace the pair $v$
and $v_0$ by $(1,0,1-n)$ and the $-2$ vector $(1,-\LB,n-1)$. 
\begin{enumerate}
\item
\label{example-item-geometric-meaning-of-reflections-rank-1-case}
The reflection of the Hilbert scheme $S^{[n]}\cong \M(1,0,1-n)$ 
in a $-2$ vector of the form $(1,-\LB,n-1)$ is a self dual 
stratified elementary transformation, provided the line bundle 
$\LB$ satisfies Condition \ref{condition-minimality}
(\cite{markman-reflections} Section 3.2). 
The center $(S^{[n]})^1$ of the operation is, 
birationally, a $\PP^1$-bundle over
$\M(0,\LB,-2)$. $(S^{[n]})^1$ consists of length $n$ subschemes, 
whose image in
$\linsys{\LB}^*\cong\PP^{n-1}$ is contained in some hyperplane. 
When the hyperplane is unique, 
the twist $\Ideal{}\otimes\LB$, of the ideal sheaf by $\LB$, 
has a one-dimensional space of global sections. The quotient
$F:=\Ideal{}\otimes\LB/[H^{0}(\Ideal{}\otimes\LB)\otimes\StructureSheaf{S}]$
is a sheaf $F\in \M(0,\LB,-2)$, 
of Euler characteristic $-2$, supported on a curve $C$ in
the linear system $\linsys{\LB}$. 
\item
The case $n=1$ can be considered as a special case, if $S$ contains an
irreducible $-2$ curve $\Sigma\cong\PP^1$ and the line bundle
$\LB$ is taken to be $\StructureSheaf{S}(\Sigma)$. In that case, 
$(S^{[1]})^1=\Sigma$. 
\item
\label{example-item-geometric-meaning-of-reflections-rank-1-n-2}
Specialize to the case 
$n=2$ and let $\LB$ be the pull-back of $\StructureSheaf{\PP^1}(1)$ 
via an elliptic fibration $\pi:S\rightarrow \PP^1$. Then 
$\M(0,\LB,-2)\rightarrow \PP^1$ is the relative Jacobian of degree $-2$.
$(S^{[2]})^1\rightarrow \M(0,\LB,-2)$ is the $\PP^1$-bundle 
obtained by projectivizing the ``vector bundle'' of relative first cohomologies
of a ``universal sheaf'' over $\M(0,\LB,-2)\times S$ (existence of a global 
universal sheaf is not needed as we projectivize). 
\end{enumerate}
}
\end{example}

\subsection{The class of the correspondence in terms of universal sheaves}
\label{sec-the-class-of-the-correspondence-in-terms-of-E-v}
Next, we provide an identification of the class of 
the correspondence (\ref{eq-total-brill-noether-correspondence})
in terms of the Chern classes of the universal family $\E_v$. 
Such an identification provides a formula for an automorphism of the 
cohomology ring of moduli spaces of sheaves on a more general K3 surface 
(removing Condition \ref{condition-minimality} on the Picard group of $S$). 
This enables us to degenerate the K3 surface to one with a large Picard group 
and obtain many such reflections.

\bigskip
Assume, that $v=(r,\LB,-r)$, $r\geq 1$, 
$\LB$ satisfies condition
\ref{condition-minimality} with respect to the ample line-bundle $H$, 
and $\mbox{gcd}(r,\deg{\LB})=1$
(so that there is a universal family $\E_v$ over $\M(v)\times S$). 
Note, that $v=\tau(v)$ and the Brill-Noether locus $\M(v)^1$ is a divisor. 

\begin{thm}
\label{thm-class-of-correspondence-in-stratified-elementary-trans}
\begin{enumerate}
\item
The class of the correspondence 
(\ref{eq-total-brill-noether-correspondence}), in the cohomology ring 
$H^*(\M(v) \times \M(v),\Integers)_{\rm free}$, is equal to the class 
$\gamma_{\tau}(\E_v,\E_v)$, given in 
(\ref{eq-Chow-theoretic-formula-for-gamma-tau}).
\item
Parts \ref{thm-item-gamma-g-is-a-ring-isomorphism} and 
\ref{thm-item-normalized-universal-chern-character-is-invariant} of Theorem 
\ref{thm-Gamma-v-acts-motivicly}
hold for the reflection $\tau$ and the moduli space $\M_H(v)$. 
\end{enumerate}
\end{thm}

Let $p:\M(v)\times S \rightarrow \M(v)$  be the projection.
Lemma 
4.4,
in part one of the paper \cite{markman-part-one}, 
implies the second equality below
\begin{eqnarray}
\nonumber
\gamma_{\tau}(\E_v,\E_{\tau(v)}) & := & 
c_m\left\{- \ 
\pi_{13_!}\left(
\left[\pi^*_{12}\E_v \ - \ \pi_1^! p_!\E_v
\right]^\vee
\otimes \pi_{23}^*\E_{\tau(v)}
\right)
\right\} \ \ \ =  
\\
\label{eq-two-expressions-for-gamma-tau-E-v-E-tau-v}
& = & 
c_m\left\{- \ 
\pi_{13_!}\left(
\pi_{12}^*(\E_v)^\vee
\otimes
\left[\pi^*_{23}\E_{\tau(v)} \ - \ \pi_3^! p_!\E_{\tau(v)}
\right]
\right)
\right\}.
\end{eqnarray}
The first equality above is the definition 
(\ref{eq-Chow-theoretic-formula-for-gamma-tau}).

\begin{rem}
\label{rem-K-theory-level-minus-2-case}
{\rm
Our proof of Theorem 
\ref{thm-class-of-correspondence-in-stratified-elementary-trans}
implies a K-theoretic version of part 
\ref{thm-item-normalized-universal-chern-character-is-invariant}
of Theorem  \ref{thm-Gamma-v-acts-motivicly}. 
We prove, that the correspondence
(\ref{eq-total-brill-noether-correspondence}) 
maps the class of the Fourier-Mukai transform of the 
universal sheaf $\E_v$, appropriately chosen,
to the class, in the Grothendieck K-group of $S\times \M(v)$, of another 
universal sheaf $\E'_v$ (see (\ref
{eq-correspondence-takes-k-theoretic-class-of-univ-sheaf-to-its-twist})). 
We expect that, even if we drop the assumption $\chi(v)=0$ of Theorem 
\ref{thm-class-of-correspondence-in-stratified-elementary-trans}, 
the K-theoretic equivalence 
\begin{equation}
\label{eq-conjectural-K-theoretic-equality}
(\Z\times \one_S)_![\tau(\E_v)] \ \ \ \equiv \ \ \
\E_{\tau(v)},
\end{equation}
holds for some universal sheaf $\E_{\tau(v)}$ over $\M(\tau(v))\times S$. 
Above, $\tau(\E_v)$ is the K-theoretic representative, given in
(\ref{eq-tau-of-universal-sheaf}),
of the reflection of $\E_v$ with respect to $\StructureSheaf{S}$.
On the other hand, we do not prove the equality of the 
class of the correspondence (\ref{eq-total-brill-noether-correspondence}) 
and of the Chow-theoretic representation 
(\ref{eq-Chow-theoretic-formula-for-gamma-tau})
of the class $\gamma_\tau(\E_v,\E_v)$ in the Chow ring of 
$\M(v)\times \M(v)$. 
One problem is, that we use the uniqueness statement in lemma
\ref{lemma-recovering-f},
 which is proven only on the level of singular cohomology. 
}
\end{rem}

A major role, in the proof of Theorem 
\ref{thm-class-of-correspondence-in-stratified-elementary-trans},
will be played by the following proposition. 
Let $v$ be a primitive Mukai vector with $\chi(v)\geq 0$.
Assume, that $c_1(v)=\LB$ satisfies condition
\ref{condition-minimality} and that a universal family $\E_v$
exists over $\M(v)$.

\begin{prop} 
\label{prop-decomposition-of-complex-over-interated-blow-up}
(Proposition 3.18 in \cite{markman-reflections})
For $1\leq t\leq \mu+1$, 
there exists a complex of vector bundles 
\begin{equation}
\label{eq-decomposable-complex}
B^{[t]}\rho:B^{[t]}V_0 \rightarrow B^{[t]}V_1
\end{equation}
over $B^{[t]}\M(v)$ with the following properties:
\begin{enumerate}
\item 
\label{property-equivalent-to-push-forward-of-universal-sheaf}
Over $\M(v)$, the complex 
$p_!\E_v$ is equivalent in the K-group to the complex of
$B^{[\mu+1]}V_\bullet$.
\item
\label{property-constant-rank}
The homomorphism $B^{[t+1]}\rho$ has constant co-rank $t$ along
$B^{[t+1]}\M(v)^t$. Denote the kernel and co-kernel of its restriction by
$W_t$ and $U_t$ respectively. 
These are vector bundles  over $B^{[t+1]}\M(v)^t$
of ranks $\chi(v)+t$ and $t$. 
\item
\label{prop-item-U-and-W-are-direct-images}
$W_t$ is isomorphic to the direct image of
the restriction of the universal family $\beta^*\E_v$ to 
the stratum $B^{[t+1]}\M(v)^t\times S$. The vector bundle 
$U_t(\sum_{i=t+1}^\mu B^{[t+1]}\M(v)^i)$ is isomorphic to 
the first higher direct image of 
the restriction of the universal family $\beta^*\E_v$ to 
the stratum $B^{[t+1]}\M(v)^t\times S$.
\item
\label{property-U-t}
$\beta^*B^{[t+1]}V_\bullet$ is equivalent in the K-group of $B^{[t]}\M(v)$ 
to $B^{[t]}V_\bullet- e_{t,!}\phi_t^*U_t$. 

\noindent
Above, $e_t:B^{[t]}\M(v)^t\hookrightarrow B^{[t]}\M(v)$ 
is the closed imersion and $\phi_t:B^{[t]}\M(v)^t\rightarrow B^{[t+1]}\M(v)^t$
is the natural projection from the exceptional divisor onto the center of the 
blow-up. 
\item
There exists an extension of sheaves over $B^{[t+1]}\M(v)^t\times S$ 
\begin{equation}
\label{eq-the-tautological-extension-of-families-over-BN-locus}
0\rightarrow U_t \rightarrow \F_t \rightarrow
\beta^*(\E_{\restricted{v}{\M(v)^t}})
\left(-\sum_{i=t+1}^\mu B^{[t+1]}\M(v)^i\right)
\rightarrow 0.
\end{equation}
\item
\label{prop-item-morphism-is-a-classifying-morphism-of-F-t}
$\F_t$ is a family of stable sheaves with Mukai vector $v+\vec{t}$,
where $\vec{t}$ denotes the Mukai vector $(t,0,t)$ of
the trivial rank $t$ bundle.
The morphism $B^{[t+1]}\M(v)^t\rightarrow \M(v+\vec{t})$
is induced by the classifying morphism associated to $\F_t$. 
It factors through $B^{[1]}\M(v+\vec{t})$ via a grassmannian fibration
$B^{[t+1]}\M(v)^t\rightarrow B^{[1]}\M(v+\vec{t})$.
\end{enumerate}
\end{prop}

A modified version of the proposition holds, even if a
universal family does not exist over $\M(v)$ 
(Remark 3.20 in \cite{markman-reflections}).
In that case, the vector bundles $U_t$ and $W_t$ need not exist,
but the projective bundles $\PP{W}_t$ and $\PP{U}_t$ do exist.
When several Mukai vectors are considered, we denote these bundles
by $\PP{W}_{v,t}$ and $\PP{U}_{v,t}$. The grassmannian bundle
$B^{[t+1]}\M(v)^t\rightarrow B^{[1]}\M(v+\vec{t})$ 
is isomorphic to $G(t,W_{v+\vec{t},0})$ 
(\cite{markman-reflections} Proposition 3.18 part 6). Compare with equation 
(\ref{eq-relation-between-various-grassmannian-bundles}) above. 

Let $W_{v+\vec{t}}$ be the direct image, 
of the sheaf $\F_t$ in proposition 
\ref{prop-decomposition-of-complex-over-interated-blow-up},
via the projection $B^{[t+1]}\M(v)^t\times S\rightarrow B^{[t+1]}\M(v)^t$.
Caution: $\F_0$ is the tensor product of $\E_v$ with a non-trivial 
line-bundle. 
Consequently, $W_{v+\vec{0}}$ is the tensor product of $W_{v,0}$ 
with that line bundle.

\begin{new-lemma}
\label{lemma-W-v-plus-t-is-a-direct-image}
$W_{v+\vec{t}}$ 
is a locally free sheaf, which fits in the short exact sequence
\begin{equation}
\label{eq-W-v-plus-t-as-an-extension}
0 \rightarrow U_t \rightarrow W_{v+\vec{t}} \rightarrow 
W_t\left(-\sum_{i=t+1}^\mu B^{[t+1]}\M(v)^{t\cap i}\right) 
\rightarrow 0. 
\end{equation}
The projectivization $\PP{W}_{v+\vec{t}}$ is
the pullback of the bundle
$\PP{W}_{v+\vec{t},0}$, defined over $B^{[1]}\M(v+\vec{t})$ in part
\ref{property-constant-rank} of Proposition 
\ref{prop-decomposition-of-complex-over-interated-blow-up}. 
There is a natural isomorphism of the two grassmannian bundles 
$G(t,W_{v+\vec{t},0})$ and 
$B^{[t+1]}\M(v)^t\rightarrow B^{[1]}\M(v+\vec{t})$. This isomorphism 
identifies the subbundle $U_t$ 
with the tautological subbundle of $W_{v+\vec{t}}$.
\end{new-lemma}

\noindent
{\bf Proof:}
Part 
\ref{prop-item-U-and-W-are-direct-images}
of proposition \ref{prop-decomposition-of-complex-over-interated-blow-up}
implies, that $W_t$ is the direct image of the universal family. 
We conclude, that $W_{v+\vec{t}}$ is locally free, from 
the long exact sequence of higher direct images, associated to the
short exact sequence 
(\ref{eq-the-tautological-extension-of-families-over-BN-locus}). 
The equality $\PP{W}_{v+\vec{t}}=\PP{W}_{v+\vec{t},0}$ was proven
in part 6 of Proposition 3.18 in \cite{markman-reflections}. 
The identification of $U_t$, as the tautological subbundle, 
is included in part 8 of Proposition 3.18 in \cite{markman-reflections}. 
\EndProof

\medskip
Properties \ref{property-constant-rank} and  \ref{property-U-t} in Proposition
\ref{prop-decomposition-of-complex-over-interated-blow-up} 
imply the following K-theoretic decomposition over $B^{[1]}\M(v)$,
of the K-theoretic pullback of $p_!\E_v$. 
\begin{equation}
\label{eq-decomposition-of-K-pushforward-of-universal-sheaf-over-B-1-M-v}
\beta^!p_!\E_v \ \equiv \ 
W_0
\ - \ 
\sum_{t=1}^\mu e_{t,!}\tilde{\phi}_t^*U_t,
\end{equation}
where $\tilde{\phi}_t$ is the composition 
$\phi_t\circ\beta : B^{[1]}\M(v)^t\rightarrow B^{[t+1]}\M(v)^t$
and $e_t:B^{[1]}\M(v)^t\hookrightarrow B^{[1]}\M(v)$ is the closed imersion
(we used also lemma
\ref{lemma-K-theoretic-pullback-of-a-sheaf-supported-on-a-divisor}
for the equality $\beta^!e_{t_!}\phi_t^*U_t=e_{t_!}\tilde{\phi}_t^*U_t$).
The vector bundle $W_0$ vanishes, if $\chi(v)\leq 0$, and is of rank $\chi(v)$
otherwise. 

\subsection{Notation}
\label{sec-notation-minus-2-vector}

We summerize here some of the notation, that will be 
repeatedly used throughout the rest of sections
\ref{sec-self-dual-stratified-elemetary-transformations} and 
\ref{sec-stratified-elementary-trans-plus-2-vectors}. 
This summary is intended for reference; the notation is introduced 
again when first mentioned. Additional notation is
summerized in section \ref{sec-notation}.
Given a proper morphism $f:X\rightarrow Y$, 
we denote by $f_*$ the sheaf theoretic push-forward and by $f_!$ the
pushforward on the level of Grothendieck
$K_0$-groups of coherent sheaves. 
Similarly, $f^*$ denotes the sheaf theoretic 
pullback, while $f^!$ the K-theoretic one. 
We denote the pullback, in singular cohomology, 
and the Gysin map by $f^*$ and $f_*$ as well. 

Let $\tau_{v_0}$ be the reflection
(\ref{eq-tau}), of the Mukai lattice, with respect to
a $-2$-vector $v_0$. We set $\tau:=\tau_{v_0}$, where $v_0=(1,0,1)$.
We denote by $\tau$ also the related K-theoretic operation, 
on the universal family of sheaves, 
defined in equation (\ref{eq-tau-of-universal-sheaf}).

The correspondence $\Z$ is the one in Definition
\ref{def-Z-for-cotangent-bundles}
part \ref{correspondence-fiber-products-of-dual-grassmanian-bundles}. 
The notation of the normal crossing model $\Z^{[1]}$ of $\Z$, as well as the 
components $\Z_t$, $\Z^{[1]}_t$, and the birational models $\Delta_t^{[1]}$ 
and $\Delta^{[t+1]}_t$, is discussed in Remark 
\ref{rem-notation-Z-vs-Delta}. 
$\Delta^{[1]}_{\leq t}$ is the union of components of $\Z^{[1]}$,
with index $\leq t$. 
$I_t$ denotes the incidence divisor in the fiber product $\Delta_t^{[t+1]}$,
of dual Grassmannian bundles. 
We set $I^{[1]}_t:=\Delta_t^{[1]}\cap \Delta^{[1]}_{\leq (t-1)}$. 
In section 
\ref{sec-self-dual-stratified-elemetary-transformations}, 
the correspondence $\Z$ and its components $\Z_t$ 
are described in equations
(\ref{eq-total-brill-noether-correspondence}) and 
(\ref{eq-closure-of-fiber-product-of-grassmannian-bundles}), 
incorporating the identification, of the stratified
elementary transform of the moduli space $\M(v)$, as a moduli space as well. 

The morphism $p:\M(v)\times S \rightarrow \M(v)$  is the projection. 
The projection $\pi_{ij}$ is 
from $\M(v)\times S\times \M(v)$ on the product of
the $i$-th and $j$-th factors. We will denote by 
$\pi_{ij}$ also the restriction of the latter projection to $\Z\times S$. 
The two projections from $\Z^{[1]}$ to $\M(v)$ are denoted
by $\tilde{\pi}_i$, $i=1,2$.
The morphism
$\tilde{\pi}_{ij}$, from $\Z^{[1]}\times S$ to the product of
the $i$-th and $j$-th factors of $\M(v)\times S\times \M(v)$, 
will denote the composition of $\Z^{[1]}\times S\rightarrow \Z\times S$
with $\pi_{ij}$. 

$\M(v)^t$ denotes a Brill-Neother stratum, introduced in 
(\ref{eq-brill-noether-stratification}). 
We denote blow-up morphisms, as well as their compositions,
by $\beta$. Examples include
$\beta:B^{[i]}\M(v)^t\rightarrow \M(v)^t$, for $i<t$, and
$\beta:\Delta_t^{[1]}\rightarrow \Delta_t^{[t+1]}$. 
The projective bundle $\phi_t: B^{[t]}\M(v)^t\rightarrow B^{[t+1]}\M(v)^t$ is
the natural morphism from the exceptional divisor to the center of the 
blow-up. The morphism 
$\tilde{\phi}_t: B^{[1]}\M(v)^t\rightarrow B^{[t+1]}\M(v)^t$
is the composition of $\phi_t$ with the blowing-down morphism
$B^{[1]}\M(v)^t\rightarrow B^{[t]}\M(v)^t$. 
We denote by $B^{[1]}\M(v)^{t\cap k}$ the intersection 
$B^{[1]}\M(v)^t\cap B^{[1]}\M(v)^k$.

The projection $\delta_{ij}$, $1\leq i,j\leq 3$, 
is from $\Delta_t^{[t+1]}\times S$ to the product of the $i$-th and $j$-th 
factor of $B^{[t+1]}\M(v)^t\times S\times B^{[t+1]}\M(v)^t$.
The morphism $\tilde{\delta}_{ij}$, from $\Delta_t^{[1]}\times S$, 
is the composition of $\delta_{ij}$ with
$\Delta_t^{[1]}\times S\rightarrow \Delta_t^{[t+1]}\times S$.

We will use $e_t$ to denote various 
embeddings, including the following
$B^{[i]}\M(v)^t\hookrightarrow B^{[i]}\M(v)$,
$\Delta_t^{[1]}\hookrightarrow \Z^{[1]}$, and 
$\Delta_t^{[1]}\times S\hookrightarrow \Z^{[1]}\times S$.

We denote by $\E_v$ the universal sheaf over $\M(v)\times S$. 
The vector bundles $U_t$ and $W_t$, over $B^{[t+1]}\M(v)^t$, are 
introduced in proposition
\ref{prop-decomposition-of-complex-over-interated-blow-up}.
The vector $\vec{t}$ is the Mukai vector $(t,0,t)$ 
of the trivial rank $t$ vector bundle.
The vector bundle $W_{v+\vec{t}}$, over $B^{[t+1]}\M(v)^t$, is introduced
in (\ref{eq-W-v-plus-t-as-an-extension}).
The family $\F_t$, over $B^{[t+1]}\M(v)^t\times S$, is the one in
extension 
(\ref{eq-the-tautological-extension-of-families-over-BN-locus}).


\subsection{Proof of Theorem 
\ref{thm-class-of-correspondence-in-stratified-elementary-trans}}
\label{sec-proof-of-theorem-on-class-of-correspondence-in-minus-2-reflection}

Lemma \ref{lemma-recovering-f} 
and equation (\ref{eq-chern-character-of-tau-E-v}) 
reduce the proof to the verification of the 
relation
\begin{equation}
\label{eq-correspondence-takes-normalizes-chern-char-to-same}
\Z_*\otimes 1_S\left(ch(\eta_v)\cdot ch(\E_v-p^!p_!\E_v)\cdot \sqrt{td_S}
\right) \ = \
ch(\eta_v)\cdot ch(\E_v)\cdot \sqrt{td_S},
\end{equation}
where $\eta_v$ is the class in 
(\ref{eq-invariant-normalized-class-of-chern-character-of-universal-sheaf}).


\smallskip
\noindent
We perform a sequence of reductions of the equality 
(\ref{eq-correspondence-takes-normalizes-chern-char-to-same})
to simpler ones.

{\bf Step I}:
We would like to carry out the calculation on the level 
of K-groups of complexes of sheaves on moduli. 
For this purpose, we will apply the push-forward operation $\Z_*$ 
to $\E_v$ rather than its normalization $\eta_v\otimes\E_v$ 
(because $\eta_v$ is only $\RationalNumbers$-Cartier). 
We need to calculate the
difference between $\eta_v$ and its image under $\Z_*$. 
Denote by $D_1$ the divisor $\M(v)^1$. 


Choose a universal family $\E_v$ normalized as in Lemma 
\ref{lemma-higher-direct-image-of-E-v-is-canonical-bundle-along-fibers}
part \ref{lemma-item-univ-family-restricts-with-degree-1-r}.
Part 
\ref{lemma-item-degree-theta-v-v} of Lemma 
\ref{lemma-higher-direct-image-of-E-v-is-canonical-bundle-along-fibers}
implies the equality 
\begin{equation}
\label{eq-push-forward-by-Z-of-eta-v}
\Z_*(\eta_v^{-1}) \ \ = \ \ \eta_v^{-1}-D_1
\end{equation}
We know, that $\Z_*$ is a ring automorphism
(\cite{markman-reflections} Theorem 1.2).
Hence, if $\Z_*$ further satisfies
(\ref{eq-correspondence-takes-normalizes-chern-char-to-same}), 
we get  
\[
\begin{array}{lccc}
\Z_*\otimes 1_S\left(ch[\tau(\E_v)]\sqrt{td_S}\right) & = & 
\Z_*\otimes 1_S\left(ch[\eta_v^{-1}\otimes\eta_v\otimes\tau(\E_v)]
\sqrt{td_S}
\right) 
&
\stackrel{(\ref{eq-correspondence-takes-normalizes-chern-char-to-same})}{=} 
\\
\Z_*\left(ch[\eta_v^{-1}]\right)\cdot 
ch[\eta_v\otimes\E_v]\sqrt{td_S} 
& \stackrel{(\ref{eq-push-forward-by-Z-of-eta-v})}{=} & 
ch[\E_v]\cdot\sqrt{td_S}\cdot p^*ch[\StructureSheaf{\M(v)}(-D_1)]
&  = 
\\
ch[\E_v(-D_1)]\cdot\sqrt{td_S}
\end{array}
\]
We see, that 
(\ref{eq-correspondence-takes-normalizes-chern-char-to-same})
implies the equality 
\begin{equation}
\label{eq-correspondence-takes-non-normalizes-chern-char-to-its-twist}
\Z_*\otimes 1_S\left(ch[\tau(\E_v)]\sqrt{td_S}\right) \ = \ 
ch[\E_v(-D_1)]\cdot\sqrt{td_S}.
\end{equation}
The same argument shows, furthermore, that 
(\ref{eq-correspondence-takes-normalizes-chern-char-to-same})
is equivalent to 
(\ref{eq-correspondence-takes-non-normalizes-chern-char-to-its-twist}).

{\bf Step II}:
Next, we translate the equality 
(\ref{eq-correspondence-takes-non-normalizes-chern-char-to-its-twist})
to a K-theoretic equivalence. 

\begin{new-lemma}
\label{lemma-ch-intertwines-cohomological-and-K-theoretic-pushforward}
We have the relations
\begin{eqnarray}
\label{eq-Z-takes-td-to-td}
\Z_*(td_{\M(v)}) & = & td_{\M(v)} \ \ \mbox{and}
\\
\label{eq-ch-character-intertwines-cohomological-and-K-theoretic-pushforward}
(\Z\times 1_S)_*(ch(\alpha)) & = &
ch\left(\pi_{12_!}\left[\pi_{23}^!\alpha\stackrel{L}{\otimes}
\StructureSheaf{\Z\times S}
\right]\right)
\end{eqnarray}
for every class $\alpha$ in the derived category of $\M(v)\times S$. 
In other words, the Chern character intertwines the 
operations of pushforward by $\Z\times S$ 
in cohomology and K-theory. 
\end{new-lemma}

\noindent
{\bf Proof:}
Proposition \ref{prop-rational-singularities-of-correspondence} 
part \ref{prop-item-Z-is-a-fiber} exhibits 
the class $\Z$ as a special fiber over $0\in B$, in a flat 
family $\Y\rightarrow B$. 
The variety $\Y$ is embedded in the fiber product $\X\times_B\X'$
of two flat families $\X\rightarrow B$ and 
$\X'\rightarrow B$.
Both fibers of $\X$ and $\X'$ over $0\in B$ are isomorphic to
$\M(v)$. The complement $\Y\setminus \Z$, of the special fiber,
is the graph of an isomorphism of the
restriction of the two families to $B\setminus \{0\}$. 
Equation (\ref{eq-Z-takes-td-to-td}) follows, since it holds
for a generic fiber of $\Y$. 

Let $\iota : \Z\times S \hookrightarrow \M(v)\times S\times \M(v)$
be the closed imersion. 
The right hand side of equation
(\ref{eq-ch-character-intertwines-cohomological-and-K-theoretic-pushforward})
is equal to the right hand side of 
\begin{equation}
\label{eq-ch-character-intertwines-cohomological-and-K-theoretic-pushforward2}
(\Z\times 1_S)_*(ch(\alpha)) \ \ = \ \ 
\pi_{12_*}\left[
\pi_{23}^*ch(\alpha)\cdot 
ch\left(\iota_!\StructureSheaf{\Z\times S}\right)\cdot 
\pi_3^* td_{\M(v)}
\right]
\end{equation}
(Grothendieck-Riemann-Roch \cite{fulton} Theorem 15.2 page 286). 
We will now prove equation (\ref
{eq-ch-character-intertwines-cohomological-and-K-theoretic-pushforward2}).
The Chern character 
of the structure sheaf of $\Z\times S$ is a pullback of that of $\Z$ in 
$\M(v)\times \M(v)$. The latter varies in a flat family in the above 
mentioned deformation. 
Denote by 
$(\ref{eq-ch-character-intertwines-cohomological-and-K-theoretic-pushforward2}
)_b$ the equation obtained when we replace $\Z$ in 
(\ref{eq-ch-character-intertwines-cohomological-and-K-theoretic-pushforward2})
by the graph $\Y_b$,  of the isomorphism $\X_b\cong \X'_b$, and 
we replace $ch(\alpha)$ by an arbitrary cohomology class
of $\X_b\times S$. Equation 
$(\ref{eq-ch-character-intertwines-cohomological-and-K-theoretic-pushforward2}
)_b$ follows by calculating $ch(\iota_!\StructureSheaf{\Y_b\times S})$ 
in terms of the Todd class of $\Y_b$, via Grothendieck-Riemann-Roch, 
and using the projection formula. Equality 
(\ref{eq-ch-character-intertwines-cohomological-and-K-theoretic-pushforward2})
follows by taking the limit.
\EndProof

%

\medskip
The equality 
(\ref{eq-correspondence-takes-non-normalizes-chern-char-to-its-twist})
would follow from 
(\ref{eq-ch-character-intertwines-cohomological-and-K-theoretic-pushforward})
and the equivalence
\begin{equation}
\label{eq-correspondence-takes-k-theoretic-class-of-univ-sheaf-to-its-twist}
\pi_{12_!}\left(
\pi_{23}^!\left[\E_v-p^*p_!\E_v\restricted{\right]}{\Z\times S}\right)
\ \equiv \ 
\E_v(-D_1)
\end{equation}
in the K-group. 
Indeed, we have
\begin{eqnarray*}
(\Z\times S)_* ch(\tau(\E_v)) & 
\stackrel{
(\ref{eq-ch-character-intertwines-cohomological-and-K-theoretic-pushforward})
+{\rm G.R.R}
}{=} &
\pi_{12_*}\left[
\pi_3^*td_{\M(v)}\cdot
\left(
\pi_{23}^*ch(\tau(\E_v))\cdot ch(\iota_!\StructureSheaf{\Z\times S})
\right)
\right]
\\
& \stackrel{{\rm G.R.R}}{=}& 
ch\left[
\pi_{12_!}\left(
\pi_{23}^!
[
\tau(\E_v)
\restricted{]}{\Z\times S}
\right)
\right]
\\
& \stackrel{
(\ref{eq-correspondence-takes-k-theoretic-class-of-univ-sheaf-to-its-twist})}
{=} &
ch(\E_v(-D_1)).
\end{eqnarray*}

\medskip
{\bf Step III}:
Next we translate 
(\ref{eq-correspondence-takes-k-theoretic-class-of-univ-sheaf-to-its-twist})
to an equivalence over the top iterated blow-up $\Z^{[1]}$ of $\Z$. 
The blow-up $\Z^{[1]}$ is defined in equation 
(\ref{eq-definition-of-top-iterated-blow-up-of-Z}).

The K-theoretic push-forward, via 
the morphism $\beta: \Z^{[1]}\rightarrow \Z$, satisfies the equivalence
\[
\beta_!\beta^!\left(
\pi_{23}^!\left[\E_v-p^*p_!\E_v\restricted{\right]}{\Z\times S}\right)
\ \equiv \ 
\pi_{23}^!\left[\E_v-p^*p_!\E_v\restricted{\right]}{\Z\times S},
\]
by Proposition 
\ref{prop-rational-singularities-of-correspondence}
part \ref{prop-item-K-theoretic-push-forward-is-O-Z}. 
Denote by $\tilde{\pi}_i$ the projection from
$\Z^{[1]}\times S$ onto the $i$-th factor of
$\M(v)\times S \times \M(v)$.
Then 
(\ref{eq-correspondence-takes-k-theoretic-class-of-univ-sheaf-to-its-twist})
is equivalent to
\begin{equation}
\label{eq-blown-up-correspondence-takes-k-theoretic-class-of-univ-sheaf-to}
(\tilde{\pi}_{12})_!\left(
(\tilde{\pi}_{23})^!\left[\E_v-p^*p_!\E_v\right]
\right)
\ = \ 
\E_v(-D_1).
\end{equation}


{\bf Step IV}:
Next, we carry out a recursive decomposition of the class 
$\tilde{\pi}_{23}^!\left[\E_v-p^*p_!\E_v\right]$. 
Since $\chi(v)=0$, the push-forward sheaf $p_*\E_v$ vanishes 
(\cite{markman-reflections} Corollary 3.16).
The second higher direct image $R^2p_*\E_v$ vanishes as well
because $c_1(v)=\LB$ is assumed to satisfy Condition 
\ref{condition-minimality}. Thus, $p_!\E_v=R^1p_*\E_v[1]$. 

Each component $\Delta_t^{[1]}$ of $\Z^{[1]}$ 
is the iterated blow-up, of the fiber product
$\Delta_t^{[t+1]}$, of dual grassmannian bundles
(\ref{eq-Delta-t-is-a-fiber-product}).
Denote by $\delta_{1,2}$ and $\delta_{2,3}$ the two projections from
$\Delta_t^{[t+1]}\times S$ to  $B^{[t+1]}\M(v)^t\times S$.
Recall that $\Delta_0^{[1]}$ is the diagonal in 
$B^{[1]}\M(v)\times B^{[1]}\M(v)$. 
Denote by $e:\Delta_{\leq (t-1)}^{[1]}\hookrightarrow \Z^{[1]}$ 
the subscheme of $\Z^{[1]}$ corresponding to the union of components 
$\Delta_{i}^{[1]}$ with index $i\leq (t-1)$. 
Let $I^{[1]}_t$ be the intersection 
$\Delta_{\leq (t-1)}^{[1]}\cap \Delta_t^{[1]}$. 
Let us derive the following  decomposition:
\begin{eqnarray}
\label{eq-decomposition-of-universal-sheaf-on-blown-up-correspondence}
\tilde{\pi}_{23}^![\E_v-p^*p_!\E_v]
& = &
\tilde{\pi}_{23}^*
\left(\E_v
\restricted{\right)}{\Delta_0^{[1]}\times S}
\left(-\sum_{i=1}^\mu\Delta_0^{[1]}\cap\Delta_i^{[1]}\right) 
\\ 
\nonumber
& & +
\sum_{t=1}^\mu 
e_{t,!}\left[
\beta^*\delta_{23}^*\F_t-(\beta^*\delta_{23}^*W_t)(-I^{[1]}_t)
\right].
\end{eqnarray}
Above, 
$e_{t}: \Delta_t^{[1]}\times S \hookrightarrow \Z^{[1]}\times S$
is the closed imersion. 
Let
$\Delta_{\geq t}^{[1]}$ be the union of components $\Delta_i^{[1]}$, $i\geq t$.
First observe, that we have the decomposition
\[
\begin{array}{lcl}
\tilde{\pi}_{23}^*(\E_v) & = &
(\tilde{\pi}_{23}^*\E_v\restricted{)}{\Delta_0^{[1]}\times S}\left(
-\sum_{t=1}^\mu\Delta_0^{[1]}\cap \Delta_t^{[1]}
\right) +
(\tilde{\pi}_{23}^*\E_v\restricted{)}{\Delta^{[1]}_{\geq 1}\times S}
\\
&=&
\sum_{t=0}^\mu
(\tilde{\pi}_{23}^*\E_v\restricted{)}{\Delta_t^{[1]}\times S}\left(
-\sum_{i=t+1}^\mu\Delta_t^{[1]}\cap \Delta_i^{[1]}
\right) 
\\
&\stackrel{(\ref{eq-the-tautological-extension-of-families-over-BN-locus})}
{=}&
(\tilde{\pi}_{23}^*\E_v\restricted{)}{\Delta_0^{[1]}\times S}\left(
-\sum_{t=1}^\mu\Delta_0^{[1]}\cap \Delta_t^{[1]}
\right) +
\sum_{t=1}^\mu 
e_{t,!}\beta^*\left[
\delta_{23}^*(\F_t-U_t)
\right].
\end{array}
\]
In the last equality above we used also Lemma
\ref{lemma-structure-of-correspondence} part 
\ref{lemma-item-pull-back-of-BN-divisors-to-correspondence}.

The decomposition 
(\ref{eq-decomposition-of-universal-sheaf-on-blown-up-correspondence})
follows from the Claim \ref{claim-decomposition-of-p-lower-shriek-E-v}, 
which is an analogue of 
(\ref{eq-decomposition-of-K-pushforward-of-universal-sheaf-over-B-1-M-v}).
Let $\tilde{\pi}_i$, $i=1,2$, be the two projections from $\Z^{[1]}$ to
$\M(v)$. Denote by $\delta_i$, $i=1,2$,  the two projections from
$\Delta_t^{[t+1]}$ to $B^{[t+1]}\M(v)^t$. 

\begin{claim}
\label{claim-decomposition-of-p-lower-shriek-E-v}
\[
\tilde{\pi}_{2}^!p_!\E_v
\ = \ 
\sum_{t=1}^\mu 
e_{t,!} [(\beta^*\delta_{2}^*W_t)(-I^{[1]}_t)-\beta^*\delta_{2}^*U_t],
\]
where $e_t:\Delta_t^{[1]}\hookrightarrow \Z^{[1]}$ is the closed imersion.
\end{claim}

\noindent
{\bf Proof:}
The morphism $\tilde{\pi}_2:\Z^{[1]}\rightarrow \M(v)$
does not factor through $B^{[t]}\M(v)$. Nevertheless, 
the restriction of $\tilde{\pi}_2$ to $\Delta_{\leq (t-1)}^{[1]}$
does factor through $B^{[t]}\M(v)$. 
Let $B^{[t]}(\rho\restricted{)}{\Delta_{\leq (t-1)}^{[1]}}$ 
be the pullback of the complex (\ref{eq-decomposable-complex}), 
from $B^{[t]}\M(v)$ to $\Delta_{\leq (t-1)}^{[1]}$, 
via the factorization of $\tilde{\pi}_2$.
We prove, by descending induction on $t$, the decomposition over $\Z^{[1]}$
\begin{equation}
\label{eq-decomposition-induction-statement}
\tilde{\pi}_2^!p_!\E_v \ \equiv \
\sum_{i=t}^\mu 
e_{i,!}\left[(\beta^*\delta_2^*W_i)(-I^{[1]}_i)-\beta^*\delta_2^*U_i\right]
+ e_![B^{[t]}(\rho\restricted{)}{\Delta_{\leq (t-1)}^{[1]}}].
\end{equation}
The claim follows from the case $t=1$, since $B^{[1]}(\rho)$
is an isomorphism over $\Delta_0^{[1]}$ (and the class of the corresponding 
complex is trivial in $K$-Theory).
The case $t=\mu+1$ is part 
\ref{property-equivalent-to-push-forward-of-universal-sheaf} of Proposition 
\ref{prop-decomposition-of-complex-over-interated-blow-up}.
Assume, that the equivalence holds for $t+1$.
The equivalence (\ref{eq-decomposition-induction-statement})
in the $t+1$ case can be rewritten in the form
\begin{eqnarray*}
\tilde{\pi}_2^!p_!\E_v & \equiv &
e_![B^{[t+1]}(\rho\restricted{)}{\Delta_{\leq (t-1)}^{[1]}}]
+ e_{t,!}[B^{[t+1]}(\rho\restricted{)}{\Delta_t^{[1]}}(-I^{[1]}_t)] +
\\& & +
\sum_{i=t+1}^\mu 
e_{i,!}\left[(\beta^*\delta_2^*W_i)(-I^{[1]}_i)-\beta^*\delta_2^*U_i\right].
\end{eqnarray*}
Let $\iota: I^{[1]}_t\hookrightarrow \Delta^{[1]}_{\leq(t-1)}$ be the 
inclusion map. The claim follows from the equivalences
\begin{eqnarray}
\label{eq-B-tplus-1-rho-differs-from-B-t-rho-along-incidence-divisor}
B^{[t+1]}(\rho\restricted{)}{\Delta_{\leq (t-1)}^{[1]}}
& \equiv &
B^{[t]}(\rho\restricted{)}{\Delta_{\leq (t-1)}^{[1]}}
-\iota_!(\delta_2^*U_t\restricted{)}{I^{[1]}_t}, 
\\
\label{eq-B-t-plus-1-rho-is-difference-between-kernel-and-cokernel}
B^{[t+1]}(\rho\restricted{)}{\Delta_t^{[1]}}
& \equiv & 
\beta^*\delta_2^*[W_t-U_t], \ \ \ \mbox{and}
\\
\nonumber
\beta^*\delta_2^*U_t &\equiv& 
\beta^*\delta_2^*U_t(-I^{[1]}_t)+
\iota_!(\delta_2^*U_t\restricted{)}{I^{[1]}_t}.
\end{eqnarray}
The equivalence 
(\ref{eq-B-tplus-1-rho-differs-from-B-t-rho-along-incidence-divisor})
follows from part \ref{property-U-t}
of proposition \ref{prop-decomposition-of-complex-over-interated-blow-up}
and Lemma \ref{lemma-structure-of-correspondence} part 
\ref{lemma-item-pull-back-of-BN-divisors-to-correspondence}
(the lemma is applied to each of the components 
$\Delta^{[1]}_j$ of $\Delta^{[1]}_{\leq(t-1)}$, taking the indices 
$i$ and $t$ of Lemma \ref{lemma-structure-of-correspondence} part 
\ref{lemma-item-pull-back-of-BN-divisors-to-correspondence} to be 
$t$ and $j$ in 
(\ref{eq-B-tplus-1-rho-differs-from-B-t-rho-along-incidence-divisor})
respectively).
The equivalence 
(\ref{eq-B-t-plus-1-rho-is-difference-between-kernel-and-cokernel})
follows from part \ref{property-constant-rank}
of proposition \ref{prop-decomposition-of-complex-over-interated-blow-up}.
\EndProof

\medskip
Equality 
(\ref{eq-blown-up-correspondence-takes-k-theoretic-class-of-univ-sheaf-to})
would follow from the decomposition 
(\ref{eq-decomposition-of-universal-sheaf-on-blown-up-correspondence}),
the vanishing
\begin{eqnarray}
\label{eq-vanishing-of-k-theory-pushforward-of-pullback-of-F-t}
\tilde{\pi}_{12_!}\left(
e_{t,!}\beta^*[\delta_{23}^*(\F_t)]
\right) & = & 0, \ \ \ t\geq 1
\\
\label{eq-vanishing-of-k-theory-pushforward-of-pullback-of-W-t-minus-I}
\tilde{\pi}_{12_!}\left(
e_{t,!}[(\beta^*\delta_{23}^*W_t)(-I^{[1]}_t)]
\right) & = & 0, \ \ \ t\geq 1
\end{eqnarray}
(proven below in steps V and VI respectively), 
and the following claim

\begin{claim}
\label{claim-push-forward-of-ideal-sheaf-is-ideal-sheaf}
Let $\beta:B^{[1]}\M(v)\rightarrow \M(v)$ be the iterated blow-down morphism.
Denote by $\widetilde{D}_t$ the divisor $B^{[1]}\M(v)^t$. 
Then
\[
\beta_{!}\left[
(\beta^*\E_v)\left(-\sum_{t=1}^\mu \widetilde{D}_t\right)
\right] \ = \ 
\beta_{!}\beta^*\left[\E_v(-\M(v)^1)
\right] \ = \ 
\E_v(-\M(v)^1). 
\]
\end{claim}

\noindent
{\bf Proof:} Note, that the sheaf theoretic pullback 
$\beta^*[\E_v(-\M(v)^1)]$ represents 
the K-theoretic pullback, because $\E_v$ is flat over $\M(v)$ and $\M(v)^1$ 
is a Cartier divisor. 
The second equality in the claim follows from the 
projection formula and the equivalence
$\beta_!\StructureSheaf{B^{[1]}\M(v)}\equiv \StructureSheaf{\M(v)}$. 
The latter equivalence is due to the fact, that 
$\beta$ is the composition of blow-ups with smooth centers. 
Note, next, that we have the equality
\[
\beta^*\left[\E_v(-\M(v)^1)
\right] \ = \ 
(\beta^*\E_v)\left(-\sum_{t=1}^\mu t\cdot \widetilde{D}_t\right)
\]
because the divisor $\M(v)^1$ has multiplicity $t$ along the locus
$\M(v)^t$. The determination of the multiplicity is 
an analogue of the Riemann-Kempf theorem, which 
expresses the fact, that the Brill-Noether stratification
is defined by a section of a homomorphism bundle and the section
is transversal to the relative determinantal stratification. 
The transversality follows from the results of 
\cite{markman-reflections} Proposition 3.18 and the general, well known, 
criterion mentioned in part 2 of 
Lemma 2.10 of \cite{markman-reflections}). 

The proof of the first equality
requires the following decomposition:
\begin{eqnarray*}
\beta^*\E_v\left(-\sum_{t=1}^\mu t\cdot \widetilde{D}_t\right)
& = & 
\beta^*\E_v\left(-\sum_{t=1}^\mu \widetilde{D}_t\right)
\ - \ \beta^*\E_v\left(-\sum_{t=1}^\mu \widetilde{D}_t
\restricted{\right)}{\widetilde{D}_2} - \ \cdots \ -
\\
& & -
\left[
\sum_{k=1}^{t-1} \beta^*\E_v
\left(-\sum_{i=1}^{t-1}i\cdot \widetilde{D}_i \ -k\cdot \widetilde{D}_t \ 
-\sum_{j=t+1}^\mu\widetilde{D}_j
\right)
\restricted{\right]}{\widetilde{D}_t} - \ \cdots \ -
\\
& & -
\left[
\sum_{k=1}^{\mu-1} \beta^*\E_v
\left(-\sum_{i=1}^{\mu-1}i\cdot \widetilde{D}_i \ -k\cdot \widetilde{D}_\mu \ 
\right)
\restricted{\right]}{\widetilde{D}_\mu}
\end{eqnarray*}

The claim would follow from the projection formula and the vanishing
\begin{equation}
\label{eq-vanishing-of-pushforward-of-line-bundle-on-D-t}
\beta_!
\StructureSheaf{\widetilde{D}_t}
\left(-\sum_{i=1}^{t-1}i\cdot \widetilde{D}_i \ -k\cdot \widetilde{D}_t \ 
-\sum_{j=t+1}^\mu\widetilde{D}_j
\right) \ \ \equiv \ \ 0, \ \ \mbox{for} 
\ \ 2\leq t \leq \mu, \ \ \mbox{and} \ \ 1 \leq k \leq t-1.
\end{equation}
$\widetilde{D}_t$ is the iterated blow-up of $B^{[t]}\M(v)^t$. 
Let $B^{[t]}\M(v)^{t\cap 1}$ be the intersection 
$B^{[t]}\M(v)^t \cap B^{[t]}\M(v)^1$. 
The pullback to $\widetilde{D}_t$ 
of the line-bundle 
$\StructureSheaf{B^{[t]}\M(v)^t}(-B^{[t]}\M(v)^{t\cap 1})$ 
is $\StructureSheaf{\widetilde{D}_t}(-
\sum_{i=1}^{t-1}i\cdot \widetilde{D}_t\cap\widetilde{D}_i)$. 
Consequently, the K-theoretic push-forward of 
$\StructureSheaf{\widetilde{D}_t}
(-\sum_{i=1}^{t-1}i\cdot \widetilde{D}_t\cap\widetilde{D}_i)$ to 
$B^{[t]}\M(v)^t$ is
$\StructureSheaf{B^{[t]}\M(v)^t}(-B^{[t]}\M(v)^{t\cap 1})$. 
The vanishing (\ref{eq-vanishing-of-pushforward-of-line-bundle-on-D-t}) 
reduces to 
\begin{equation}
\label{eq-push-forward-of-line-bundle-of-relative-degree-k-t}
\beta_!
\left[
\StructureSheaf{B^{[t]}\M(v)^t}
(-B^{[t]}\M(v)^{t\cap 1}) \ \otimes \ N_{B^{[t]}\M(v)^t/B^{[t]}\M(v)}^{-k}
 \ \otimes \ 
\StructureSheaf{B^{[t]}\M(v)^t}\left(-\sum_{j=t+1}^\mu\widetilde{D}_j
\right)
\right] \ \ \equiv \ \ 0, 
\end{equation}
for $2 \leq t \leq \mu$ and $1 \leq k \leq t-1$.

$B^{[t]}\M(v)^t$ is the projectivised relative cotangent bundle
$\PP{T}^*_{B^{[t+1]}\M(v)^t/B^{[1]}\M(v+\vec{t})}$
of the grassmannian bundle 
$B^{[t+1]}\M(v)^t \rightarrow  B^{[1]}\M(v+\vec{t})$. 
The projectivised cotangent bundle $\PP{T}^*G(t,2t)$,
of a grassmannian fiber, is the projectivization of the 
homomorphism bundle
$\Hom(q_{G(t,2t)},\tau_{G(t,2t)})$ from the tautological quotient to 
the tautological sub-bundle. The degeneracy divisor in 
$\PP\Hom(q_{G(t,2t)},\tau_{G(t,2t)})$ 
has degree $t$ in each $\PP^{t^2-1}$-fiber over $G(t,2t)$. 
We conclude, that the degeneracy locus 
$B^{[t]}\M(v)^{t\cap 1}$,
in the exceptional divisor $B^{[t]}\M(v)^t$, 
has degree $t$ in each $\PP^{t^2-1}$-fiber over $B^{[t+1]}\M(v)^t$. 
The line-bundle 
$\StructureSheaf{B^{[t]}\M(v)^t}\left(-\sum_{j=t+1}^\mu\widetilde{D}_j
\right)$ is pulled-back from $B^{[t+1]}\M(v)^1$
(see \cite{markman-reflections} Corollary 3.19,
or equation (\ref{eq-conditioned-equality-of-catier-divisor}) above).
Hence, the line bundle in 
(\ref{eq-push-forward-of-line-bundle-of-relative-degree-k-t})
has degree $k-t$ along each $\PP^{t^2-1}$-fiber. The vanishing
(\ref{eq-push-forward-of-line-bundle-of-relative-degree-k-t})
follows.
This completes the proof of the claim.
\EndProof

\medskip
{\bf Step V:}
Next, we prove the vanishing 
(\ref{eq-vanishing-of-k-theory-pushforward-of-pullback-of-F-t}). 
The composition of the pullback $\beta^!$, from the K-group of 
$\Delta^{[t+1]}_t$ to  the K-group of $\Delta^{[1]}_t$, 
followed by the pushforward $\beta_!$, 
is the identity, since $\beta:\Delta^{[1]}_t\rightarrow \Delta^{[t+1]}_t$ 
is the composition of a sequence of blow-ups with smooth centers. 
The equality $\tilde{\pi}_{12_!}\circ e_{t,!}=\delta_{12_!}\circ \beta_!$
reduces the vanishing 
(\ref{eq-vanishing-of-k-theory-pushforward-of-pullback-of-F-t})
to the vanishing
\begin{eqnarray}
\label{eq-vanishing-of-k-theory-pushforward-of-pullback-of-F-t-by-delta}
\delta_{12_!}\left(\delta_{23}^*(\F_t)
\right) & = & 0, \ \ \ t\geq 1.
\end{eqnarray}

Let $W_t$ be the vector bundle in Proposition
\ref{prop-decomposition-of-complex-over-interated-blow-up}. 
The direct image of $\F_t$, via the projection
$B^{[t+1]}\M(v)^t\times S\rightarrow B^{[t+1]}\M(v)^t$,
is the rank $2t$ vector bundle, denote by  $W_{v+\vec{t}}$ 
(Lemma \ref{lemma-W-v-plus-t-is-a-direct-image}). 
These vector bundles are related by the 
short exact sequence (\ref{eq-W-v-plus-t-as-an-extension}).
Let 
\begin{equation}
\label{eq-pi-from-grass-bundle-to-BN-locus}
\pi \ : \ G(t,W_{v+\vec{t}}) \ \rightarrow \ B^{[t+1]}\M(v)^t
\end{equation}
be the Grassmannian bundle of $t$-dimensional subspaces.
Lemma \ref{lemma-W-v-plus-t-is-a-direct-image} implies, that 
$G(t,W_{v+\vec{t}})$ is isomorphic to 
the fiber product $\Delta^{[t+1]}_t$ of $B^{[t+1]}\M(v)^t$ with itself
over $B^{[1]}\M(v+\vec{t})$. 
We choose the isomorphism, so that the projection $\pi$
in (\ref{eq-pi-from-grass-bundle-to-BN-locus}) corresponds to the projection
$\delta_1:\Delta^{[t+1]}_t\rightarrow B^{[t+1]}\M(v)^t$.
Let 
\[
\kappa \ : \ G(t,W_{v+\vec{t}}) \ \rightarrow \ B^{[t+1]}\M(v)^t
\]
be the second projection. Then equation 
(\ref{eq-vanishing-of-k-theory-pushforward-of-pullback-of-F-t-by-delta}) 
is equivalent to the following equality.
\begin{equation}
\label{eq-reformulation-vanishing-of-k-theory-pushforward}
(\pi\times 1_S)_!(\kappa\times 1_S)^*\F_t \ \ \  \equiv \ \ \ 0, \ \ \
t\geq 1.
\end{equation}

Denote by $\tau_{G(t,W_{v+\vec{t}})}$ the tautological subbundle. 
We get the short exact sequence
\begin{equation}
\label{eq-tautological-extension-of-families-over-grass}
0\rightarrow \tau_{G(t,W_{v+\vec{t}})} \rightarrow 
(\pi\times 1_S)^*\F_t \rightarrow Q
\rightarrow 0. 
\end{equation}
The quotient $Q$ is a family of stable sheaves flat over $G(t,W_{v+\vec{t}})$.
Moreover, the family $Q$ defines a classifying morphism 
\[
\kappa' \ : \ G(t,W_{v+\vec{t}}) \ \rightarrow \ 
\M(v)^t
\]
(see \cite{markman-reflections} Theorem 3.15). 
The morphism $\kappa'$ factors through 
$B^{[t+1]}\M(v)^t$ via the projection $\kappa$
(use Proposition 3.18 part 5 in \cite{markman-reflections}).
Hence, 
there is a line bundle $K$ on $G(t,W_{v+\vec{t}})$, 
such that 
$(\ref{eq-tautological-extension-of-families-over-grass})\otimes K$ 
is isomorphic to 
$\kappa^*(\ref{eq-the-tautological-extension-of-families-over-BN-locus})$.
In particular, we get the isomorphisms 
\begin{eqnarray}
\label{eq-pullback-of-F-t-by-kappa-is-K-tensor-pullback-by-pi}
(\kappa\times 1_S)^*\F_t & \cong & [(\pi\times 1_S)^*\F_t]\otimes K 
\ \ \mbox{and}
\\
\nonumber
\kappa^*U_t & \cong & \tau_{G(t,W_{v+\vec{t}})}\otimes K.
\end{eqnarray}
The projection formula and equation 
(\ref{eq-pullback-of-F-t-by-kappa-is-K-tensor-pullback-by-pi})
imply the equality 
\[
(\pi\times 1_S)_!(\kappa\times 1_S)^*\F_t 
\ = \ \F_t\otimes (\pi_! K). 
\]
Lemma \ref{lemma-higher-direct-image-of-E-v-is-canonical-bundle-along-fibers}
implies, that $K$ restricts as
$\StructureSheaf{G(t,2t)}(-1)$ along the fibers of 
$\pi$ in (\ref{eq-pi-from-grass-bundle-to-BN-locus}). 
Hence, the pushforward $\pi_! K$ vanishes and 
equation 
(\ref{eq-reformulation-vanishing-of-k-theory-pushforward}) follows.

\medskip
{\bf Step VI:}
Next, we prove the vanishing 
(\ref{eq-vanishing-of-k-theory-pushforward-of-pullback-of-W-t-minus-I}).
Part 
\ref{lemma-item-pushforward-of-ideal-sheaf-is-ideal-sheaf-of-incidence-div}
of Lemma \ref{lemma-structure-of-correspondence} implies the equality
\[
\beta_!\left[\left(\beta^*\delta_{23}^*W_t\right)(-I^{[1]}_t)\right]
\ \ \ = \ \ \ 
(\delta_{23}^*W_t)(-I_t),
\]
where $I_t$ is the incidence divisor in $\Delta_t^{[t+1]}$. 
The vanishing 
(\ref{eq-vanishing-of-k-theory-pushforward-of-pullback-of-W-t-minus-I})
reduces to
\[
\delta_{12_!}\left[(\delta_{23}^*W_t)(-I_t)\right] \ \ = \ \ 0, 
\ \ \ \mbox{for} \ \ t\geq 1.
\]
$W_t$ is defined over $B^{[t+1]}\M(v)^t$, so we may work over $\Delta_t^{[1]}$
rather than $\Delta_t^{[1]}\times S$.
Using the notation of Step V, 
the last equation translates to
\begin{equation}
\label{eq-pushforward-by-pi-kills-I-t-twist-of-pullback-of-W-t}
\pi_!\left\{
\left[\kappa^*W_t\right](-I_t)
\right\} \ \ = \ \ 0, 
\ \ \ \mbox{for} \ \ t\geq 1.
\end{equation}

The two isomorphisms 
(\ref{eq-pullback-of-F-t-by-kappa-is-K-tensor-pullback-by-pi}) 
and the sequence (\ref{eq-W-v-plus-t-as-an-extension})
imply the isomorphisms 
\begin{eqnarray*}
\kappa^*W_{v+\vec{t}} & \cong & (\pi^*W_{v+\vec{t}})\otimes K \ \ \ \mbox{and}
\\
\kappa^*\left[W_t\left(-\sum_{i=t+1}^\mu B^{[t+1]}\M(v)^{t\cap i}\right)
\right] & \cong & q_{G(t,W_{v+\vec{t}})}\otimes K,
\end{eqnarray*}
where $q_{G(t,W_{v+\vec{t}})}$ is the tautological quotient bundle over
$G(t,W_{v+\vec{t}})$.
Lemma \ref{lemma-higher-direct-image-of-E-v-is-canonical-bundle-along-fibers}
implies, that $K$ restricts as
$\StructureSheaf{G(t,2t)}(-1)$ along the fibers of 
$\pi$ in (\ref{eq-pi-from-grass-bundle-to-BN-locus}). 
The line bundles 
$\StructureSheaf{B^{[t+1]}\M(v)^t}(B^{[t+1]}\M(v)^{t\cap i})$
are pullbacks of line bundles over $B^{[1]}\M(v+\vec{t})$, when $i>t$. 
We conclude, that $\kappa^*W_t$ restricts to the Grassmannian fibers 
of $\pi$ as $q(-1)$. 
The fiber of $\Delta_t^{[t+1]}$ over $B^{[1]}\M(v+\vec{t})$
is a product $G(t,U)\times G(t,U)$, where $U$ is a $2t$-dimensional subspace 
of $H^0(F)$, for a sheaf $F$ in $\M(v+\vec{t})$. 
The incidence divisor $I_t$ restricts to the product
$G(t,U)\times G(t,U)$ as the zero locus of a section
of $\det\SheafHom(p_1^*\tau,p_2^*q)$. Hence, $I_t$ 
restricts to the Grassmannian fibers 
of $\pi$ as $\StructureSheaf{G(t,2t)}(t)$.
Lemma \ref{lemma-vanishing-of-cohomology-of-vb-on-grassmannians}
implies, that $H^i(G(t,2t),q(-t-1))$ vanishes, for all $i$. 
The vanishing 
(\ref{eq-pushforward-by-pi-kills-I-t-twist-of-pullback-of-W-t}) follows. 
This completes the proof of Theorem 
\ref{thm-class-of-correspondence-in-stratified-elementary-trans}. 
\EndProof

\medskip
Let $v$ be as in Theorem 
\ref{thm-class-of-correspondence-in-stratified-elementary-trans}. 
The following is an easy corollary of the theorem.

\begin{cor}
\label{cor-action-of-Z-on-individual-strata}
The correspondence $\Z$, given in 
(\ref{eq-total-brill-noether-correspondence}), 
acts via multiplication by $(-1)^t$ on the class of
$\M(v)^t$ in the Chow group of $\M(v)$. 
\end{cor}

Note: We present a more general proof of the equality
\[
\Z_*[\M(v)^t]=(-1)^{t(t+\chi(v))}[\M(\tau(v)^t],
\]
without the assumption $\chi(v)=0$
of Theorem \ref{thm-class-of-correspondence-in-stratified-elementary-trans}.
We assume, instead, the conjectural 
equality (\ref{eq-conjectural-K-theoretic-equality}). 
When $\chi(v)=0$, equality 
(\ref{eq-conjectural-K-theoretic-equality}) is verified in 
equality 
(\ref{eq-correspondence-takes-k-theoretic-class-of-univ-sheaf-to-its-twist}). 

\noindent
{\bf Proof of corollary \ref{cor-action-of-Z-on-individual-strata}:}
The varieties $\M(v)^t$ are irreducible, since each iterated blow-up
$B^{[t+1]}\M(v)^t$ is a grassmannian bundle over $B^{[1]}\M(v+\vec{t})$,
and $\M(v+\vec{t})$ is irreducible. It is clear, that 
$\Z_*(\M(v)^t)$ is a multiple of the class of $\M(\tau(v))^t$. 
The computation of the coefficient can thus be carried out on the level
of cohomology. We will need to use the fact, that $\Z_*$ 
preserves the multiplicative structure of the  cohomology. 
The argument given below would 
prove the statement directly, for the classes in the Chow ring,
provided we knew, that $\Z_*$ preserves the ring srtucture of
the Chow ring as well. 

$\M(v)^t$ is the degeneracy locus, of expected dimension, of a homomorphism
$\rho:V_0\rightarrow V_1$ of vector bundles. 
The homomorphism $\rho$ is the case  $B^{[\mu+1]}\rho$ of the
sequence given in equation (\ref{eq-decomposable-complex}). 
The ranks $r_0$ and $r_1$, of the vector bundles $V_0$ and $V_1$, 
satisfy $r_0-r_1=\chi(v)$. 
Proposition \ref{prop-decomposition-of-complex-over-interated-blow-up}
part \ref{property-equivalent-to-push-forward-of-universal-sheaf}
identifies the class of $V_0- V_1$, 
in the K-group of $\M(v)$, 
with the class of the pushforward $p_!\E_v$, of the universal sheaf,
via the projection $p:\M(v)\times S\rightarrow \M(v)$. 
Note the equality $\chi(\tau(v))=-\chi(v)$. 
We may assume, that $\chi(v)\geq 0$, 
possibly after interchanging $v$ with $\tau(v)$. 
Then the classes of $\M(v)^t$ and $\M(\tau(v))$ 
are given by the Porteous formulas
\begin{eqnarray*}
\M(v)^t & \equiv & \Delta_t^{(t+\chi(v))}\left(c[-p_!\E_v]\right)
\\ 
 \M(\tau(v))^t & \equiv & 
\Delta^{(t)}_{t-\chi(\tau(v))}\left(c[-p_!\E_{\tau(v)}]\right),
\end{eqnarray*}
where $c(-p_!\E_v)$ is Chern polynomial, 
and $\Delta_t^{(s)}$ is the 
determinant of the $s\times s$ matrix, $(c_{t+j-i})_{1\leq i,j\leq s}$
(Theorem 14.4 in \cite{fulton}). 

The commutativity 
$\Z_!\circ p_! = p_!\circ (\Z\times \one_S)_!$
and the equality
\[
p_!\left[\E_v-p^*p_!\E_v\right]  \ \ \ \equiv \ \ \ 
p_!\E_v -\chi(\StructureSheaf{S})\cdot p_!\E_v   \ \ \ \equiv \ \ \ -p_!\E_v,
\]
translate the equivalence (\ref{eq-conjectural-K-theoretic-equality}) to 
\begin{equation}
\label{eq-correspondence-multiplies-K-theoretic-image-of-E-v-by-minus-1}
\Z_!(-p_!\E_v) \ \ \ \equiv \ \ \ p_!\E_{\tau(v)}.
\end{equation}
$\Z$ induces an isomorphism of
cohomology {\em rings}, and thus commutes with $\Delta^{(t+\chi(v))}_t$
\[
\Z_*\left(\Delta^{(t+\chi(v))}_t(c[-p_!\E_v])\right)  \ \ = \ \ 
\Delta^{(t+\chi(v))}_t\left(c[p_!\E_{\tau(v)}]\right).
\]
The corollary follows from the equality 
\[
\Delta^{(t+\chi(v))}_t(c[p_!\E_{\tau(v)}]) \ = \ 
(-1)^{t(t+\chi(v))}\Delta_{t+\chi(v)}^{(t)}(c[-p_!\E_{\tau(v)}]) 
\]
(see Example 14.4.9 in \cite{fulton}).
\EndProof

\smallskip
The Lemmas 
below were used in the proof of Theorem 
\ref{thm-class-of-correspondence-in-stratified-elementary-trans}.

\begin{new-lemma}
\label{lemma-degree-of-determinantal-line-bundles-on-fibers-of-M-1}
\begin{enumerate}
\item
\label{lemma-item-degree-of-determinantal-line-bundles-on-fibers-of-M-1}
Let $v_0=(1,0,1)$, $v=(r,\LB,s)$ with $\mbox{gcd}(r,\deg(\LB),s)=1$. 
Assume, that $\LB$ satisfies Condition \ref{condition-minimality}. 
The class $[\PP^1]\in H^{2\dim(v)-2}(\M(v),\Integers)$, of 
a line in the $\PP^{\Abs{\chi(v)}+1}$-fiber of 
$[\M(v)^1\setminus\M(v)^2]\rightarrow \M(v+v_0)$, 
is taken to the element $(-v_0,\bullet)$ of $(v^\perp)^*$, 
under the composition
\[
H^{2\dim(v)-2}(\M(v),\Integers)
\ \cong \ H^2(\M(v),\Integers)^* \ \cong \ (v^\perp)^*
\]
of Poincare-Duality with the isomorphism $\theta_v$ 
in (\ref{eq-mukai-homomorphism}). Equivalently, 
if $w\in v^\perp$, then $(w,-v_0)$ is the degree of the restriction of the 
class $\theta_v(w)$ in (\ref{eq-mukai-homomorphism}), 
to a $\PP^{\Abs{\chi(v)}+1}$-fiber of 
$[\M(v)^1\setminus\M(v)^2]\rightarrow \M(v+v_0)$. 
\item
\label{lemma-item-brill-noether-divisor-goes-to-its-fiber}
Assume, in addition, that $\chi(v)=0$ (i.e., $r=-s$) and 
consider the composition 
\[
H^2(\M(v),\Integers) \LongIsomRightArrowOf{\theta_v^{-1}}
v^\perp \rightarrow (v^\perp)^* \LongIsomRightArrow 
H^2(\M(v),\Integers)^*
\LongIsomRightArrow 
H^{2\dim(v)-2}(\M(v),\Integers),
\]
of $\theta_v^{-1}$, the Mukai pairing, $(\theta_v^{-1})^*$, 
and Poincare-Duality. 
The class of $\M(v)^1$ is taken by the above composition to the class of a
$\PP^1$-fiber. 
\end{enumerate}
\end{new-lemma}

\noindent
{\bf Proof:}
\ref{lemma-item-degree-of-determinantal-line-bundles-on-fibers-of-M-1})
We prove this part under the assumption $\chi(v)\geq 0$. 
The proof  is similar in case $\chi(v)$ is negative.
Let $F$ be a stable sheaf in $\M(v+v_0)\setminus\M(v+v_0)^1$ and $E_0$
a sheaf in $\M(v)^1$ over $F$. 
$\PP{H}^0(F)$ is naturally identified with 
the fiber over $F$ of the $\PP^{\chi(v)+1}$-bundle 
$[\M(v)^1\setminus \M(v)^2]\rightarrow [\M(v+v_0)\setminus\M(v+v_0)^1]$.
We have the short exact sequence 
\[
0\rightarrow \StructureSheaf{S}\rightarrow F\rightarrow E_0\rightarrow 0.
\]
Equation (\ref{eq-the-tautological-extension-of-families-over-BN-locus})
implies, that the restriction of the universal family $\E_v$ to 
$\PP{H}^0(F)\times S$ fits in the short exact sequence 
\[
0 \rightarrow
\pi_1^*\StructureSheaf{\PP{H}^0(F)}(\ell-1) \rightarrow 
\pi_2^*F\otimes\pi_1^*\StructureSheaf{\PP{H}^0(F)}(\ell) \rightarrow 
\E_{\restricted{v}{\PP{H}^0(F)\times S}} \rightarrow 0
\]
for some integer $\ell$.
Assume first, that $w=v(A)$ for some complex of sheaves $A$. 
The line bundle $\theta_v(w\restricted{)}{\PP{H}^0(F)}$ is  
\begin{eqnarray*}
\det\pi_{1_!}(\pi_2^*A^\vee\otimes\E_{\restricted{v}{\PP{H}^0(F)\times S}}) 
& = &
\det\pi_{1_!}
(\pi_2^*A^\vee\otimes\pi_2^*F\otimes\pi_1^*\StructureSheaf{\PP{H}^0(F)}(\ell))
-
\det\pi_{1_!}(\pi_2^*A^\vee\otimes\pi_1^*\StructureSheaf{}(\ell-1)) 
\\
& = & 
\StructureSheaf{\PP{H}^0(F)}(\ell\cdot\chi(A,F))
-\StructureSheaf{\PP{H}^0(F)}((\ell-1)\chi(A^\vee))
\\
&=&
\StructureSheaf{\PP{H}^0(F)}((-\ell)\langle w,v+v_0\rangle \ - \ 
(1-\ell)\langle w,v_0\rangle)
\ = \ \StructureSheaf{\PP{H}^0(F)}(-\langle w,v_0\rangle).
\end{eqnarray*}
The topological translation of the above argument, 
via Grothendieck-Riemann-Roch, establishes the equality
\[
\theta_v(w)\cap [\PP^1] \ = \ -\langle w,v_0\rangle
\]
also for a non-algebraic vector $w\in v^\perp$. 

\ref{lemma-item-brill-noether-divisor-goes-to-its-fiber}) 
We know, that $\theta_v(-v_0)$ is isomorphic to
$\StructureSheaf{\M(v)}(\M(v)^1)$. 
Part \ref{lemma-item-brill-noether-divisor-goes-to-its-fiber} is an immediate 
consequence of Theorem 
\ref{thm-irreducibility} and part
\ref{lemma-item-degree-of-determinantal-line-bundles-on-fibers-of-M-1} 
of the Lemma.
\EndProof

\begin{new-lemma}
\label{lemma-higher-direct-image-of-E-v-is-canonical-bundle-along-fibers}
Let $v_0=(1,0,1)$, $v=(r,\LB,s)$ with $\mbox{gcd}(r,\deg(\LB),s)=1$ and
assume that $\LB$ satisfies Condition \ref{condition-minimality}. 
Assume, further, that either $r=1$, or $r=-s$. 
Denote by $\tau$ the rank $t$ tautological subbundle over the Grassmannian
$G(t,H)$ of $t$-dimensional subspaces of a vector space $H$. 
\begin{enumerate}
\item
\label{lemma-item-univ-family-restricts-with-degree-1-r}
There exists a universal family $\E_v$ over
$\M(v)\times S$, such that the restriction of $\det\E_v$ 
to a grassmannian fiber, of every Brill-Noether stratum, is 
$\StructureSheaf{}(1-r)$. 
\item
\label{lemma-item-univ-U-t-vs-tau}
Let $\E_v$ be as in part
\ref{lemma-item-univ-family-restricts-with-degree-1-r}. 
The vector bundle $U_t$ over $B^{[t+1]}\M(v)^t$ 
(given in (\ref{eq-the-tautological-extension-of-families-over-BN-locus}))
restricts as $\tau(-1)$ along the fibers of
$f:B^{[t+1]}\M(v)^t\rightarrow B^{[1]}\M(v+\vec{t})$.
\item
\label{lemma-item-degree-theta-v-v}
Let $\E_v$ be as in part
\ref{lemma-item-univ-family-restricts-with-degree-1-r}. 
The degree, of the restriction of the line-bundle 
$\theta_v(v)$ to a fiber of 
$[\M(v)^1\setminus\M(v)^2]\rightarrow \M(v+v_0)$,
is $\langle v,v\rangle+\chi(v)$.
\end{enumerate}
\end{new-lemma}

\noindent
{\bf Proof:} 
\ref{lemma-item-univ-family-restricts-with-degree-1-r} and 
\ref{lemma-item-univ-U-t-vs-tau})
We prove the Lemma under the assumption that $\chi(v)\geq 0$ (the case
$\chi(v)<0$ is similar).
Choose a universal family $\E_v$ (possible because
$\mbox{gcd}(r,\deg(\LB),s)=1$). 
Fix a sheaf $F_0\in \M(v+\vec{t})\setminus \M(v+\vec{t})^1$.
The family $\F_t$ restricts to a trivial family on each 
Grassmannian fiber $G(t,H^0(F_0))$ of 
$B^{[t+1]}\M(v)^t \rightarrow B^{[1]}\M(v+\vec{t})$,  
up to a twist by a line-bundle on $G(t,H^0(F_0))$. 
Say $(\F_t\restricted{)}{G(t,H^0(F_0))\times S}\cong 
F_0\otimes \StructureSheaf{G(t,H^0(F_0))}(n)$.
Then (\ref{eq-the-tautological-extension-of-families-over-BN-locus}) 
implies, that $U_t$ restricts to $\tau(n)$. 
Moreover, (\ref{eq-the-tautological-extension-of-families-over-BN-locus}) 
implies, that the restrictions of 
$R^1p_*\E_v$ and $U_t$ to $G(t,H^0(F_0))$ are isomorphic. 
Furthermore, the restriction of $\det(\E_v)$  to 
$G(t,H^0(F_0))\times \{pt\}$ is equal 
to that of $\det(\F_t)\otimes \det[R^1p_*\E_v]^*$.
The former is 
$\StructureSheaf{G(t,H^0(F_0))}(k_t)$, for some integer $k_t$, 
while the latter is 
$\StructureSheaf{G(t,H^0(F_0))}((t+r)n-(tn-1))=
\StructureSheaf{G(t,H^0(F_0))}((rn+1))$. We conclude the equality
\[
k_t \ \ = \ \ rn+1.
\]

Lemma
\ref{lemma-grassmannian-fibers-are-homologous-to-sub-grassmannians} implies,
that each fiber of 
$B^{[2]}\M(v)^1\rightarrow B^{[1]}\M(v+\vec{1})$,
is homologous to a $\PP^{\chi(v)+1}$, embedded ``linearly'' in a
grassmannian fiber of $B^{[t+1]}\M(v)^t \rightarrow B^{[1]}\M(v+\vec{t})$.
Consequently, $k_t$ is independent of $t$.

For $r=1$, we can choose $\E_v$ such that $\det(\E_v)$ restricts to the
trivial line-bundle on $\M(v)\times \{pt\}$ and $k_t=0$ for all $t$. 
If $r\geq 1$ and $k_t=1-r$, 
then $n=-1$ and $R^1p\E_v$ restricts to $\tau(-1)$.

It remains to show, that we can choose $\E_v$ with $k_t=1-r$,
for all grassmannian fibers (all $t$), provided $\chi(v)=0$.
Lemma 
\ref{lemma-degree-of-determinantal-line-bundles-on-fibers-of-M-1} 
implies, that the line bundle $\theta_v(v_0)$ restricts as 
$\StructureSheaf{\PP^{\chi(v)+1}}(2)$ on the $\PP^{\chi(v)+1}$ fibers of 
$\M(v)^1\setminus \M(v)^2$ (here we used the assumption that $\chi(v)=0$. 
Without this assumption, we can only find $w\in v^\perp$, such that
$(w,v_0)=d$, where $d=\gcd(c_1(\LB)^2,r-s)$). 
If $n$ is odd, we can replace $\E_v$ by 
$
\E_v\otimes 
p^*\theta_v\left(-\left(\frac{n+1}{2}\right)\cdot v_0\right)
$ 
to obtain a suitably normalized universal family. 

It remains to show, that $n$ can not be even.
If $n$ was even, there would be a line bundle on $\M(v)$, which restricts
with degree $1$ to a $\PP^{\chi(v)+1}$ fiber. Theorem 
\ref{thm-irreducibility} implies, that $\Pic_{\M(v)}$ is isomorphic
to the lattice $\theta_v(v^\perp_{Alg})$. The algebraic 
Mukai vectors in $v^\perp$ are spanned by $v_0$, and algebraic vectors in 
$\{v,v_0\}^\perp$. The line bundle $\theta(v_0)$ restricts with degree $2$
to a fiber. If $u\in \{v,v_0\}^\perp$,
the line bundle $\theta_v(u)$ restricts to the trivial bundle on
a fiber (Lemma
\ref{lemma-degree-of-determinantal-line-bundles-on-fibers-of-M-1}). 
Hence, $n$ can not be even. 

\ref{lemma-item-degree-theta-v-v})
Repeat the calculation in Lemma 
\ref{lemma-degree-of-determinantal-line-bundles-on-fibers-of-M-1}
with arbitrary $w$ (without the assumption that $w\in v^\perp$) and 
$\ell=-1$. We get the equality
\[
\theta_v(w)\cap [\PP^1] \ = \ \langle w,v-v_0\rangle
\]
In partcular, $\theta_v(v)\cap [\PP^1]=\langle v,v\rangle+\chi(v)$. 
\EndProof


\begin{new-lemma}
\label{lemma-vanishing-of-cohomology-of-vb-on-grassmannians}
Let $G(t,H)$ be the grassmannian of $t$-dimensional subspaces of an
$h$-dimensional vector space, $1\leq t \leq h-1$. 
Denote by $\tau$ the rank $t$ universal 
sub-bundle and let $q$ be the universal quotient bundle. 
\begin{enumerate}
\item
\label{lemma-item-kodaira-vanishing-ober-grassmannians}
If $1\leq j\leq h-1$, then $H^i(G(t,H),\StructureSheaf{G(t,H)}(-j))$
vanishes for all $i$.
\item
\label{lemma-item-vanishing-of-cohomology-of-twists-of-tau}
$H^i(G(t,H),\tau(-j))$ vanishes, for all $i$,
provided $0\leq j \leq h-1$, $(t,j)\neq (1,h-1)$, and $(t,j)\neq (h-1,1)$.
\item
\label{lemma-item-vanishing-of-cohomology-of-twists-of-tau-dual}
$H^i(G(t,H),\tau^*(-k))$ vanishes, for all $i$,
provided $1\leq k \leq h$, $(t,k)\neq (1,1)$, and 
$(t,k)\neq (h-1,h-1)$.
\item
\label{lemma-item-vanishing-of-cohomology-of-twists-of-q}
$H^i(G(t,H),q(-k))$ vanishes, for all $i$,
provided $1\leq k \leq h$, $(t,k)\neq (h-1,1)$, and $(t,k)\neq (1,h-1)$. 
\end{enumerate}
\end{new-lemma}

\noindent
{\bf Proof:}
\ref{lemma-item-kodaira-vanishing-ober-grassmannians})
Follows from Kodaira's Vanishing Theorem and the fact, that 
the canonical line bundle of $G(t,H)$ is $\StructureSheaf{G(t,H)}(-h)$. 

\ref{lemma-item-vanishing-of-cohomology-of-twists-of-tau})
Follows from part 
\ref{lemma-item-vanishing-of-cohomology-of-twists-of-tau-dual}
and Serre's Duality. 

\ref{lemma-item-vanishing-of-cohomology-of-twists-of-tau-dual})
The case $t=1$ reduces to part 
\ref{lemma-item-kodaira-vanishing-ober-grassmannians}, 
so we assume $t\geq 2$. 
We first reduce the statement to the vanishing of 
the cohomologies of a line bundle on the full flag variety.
Let $\pi_i$, $i=1,t$,  be the two projections from the partial 
flag variety $Flag(1,t,H)$.
\[
\PP{H} \ \ \LongLeftArrowOf{\pi_1} \ \ Flag(1,t,H) \ \ 
\LongRightArrowOf{\pi_t} \ \ G(t,H).
\]
Then $\tau^*_{G(t,H)}$ is isomorphic to 
$\pi_{t_*}\pi_1^*\StructureSheaf{\PP{H}}(1)$, because 
$Flag(1,t,H)$ is isomorphic to $\PP\tau_{G(t,H)}$.
Moreover, the higher direct images
$R^i\pi_{t_*}\pi_1^*\StructureSheaf{\PP{H}}(1)$ vanish, for $i>0$.
We can, in fact, carry the whole construction over the full flag variety
$\B:= SL(h)/B$, parametrizing Borel subgroups of $SL(h)$. 
The morphism $SL(h)\rightarrow SL(h)/B$ is a principal $B$-bundle over $\B$.
We get an associated $T$-bundle, via the natural homomorphism from $B$ to
the maximal torus $T\subset B$. 
Given a character $\lambda$ of $T$,
we denote by $L_\lambda$ the line bundle, associated to the $T$-bundle, via
the character $-\lambda$. 
Let $\{\epsilon_i-\epsilon_j\}_{i,j=1}^h$, $i\neq j$, denote the roots of
$\LieAlg{sl}_h$ and
$\lambda_i:=
\sum_{j=1}^i\epsilon_j-\frac{i}{h}(\epsilon_1+\cdots +\epsilon_h)$, 
$1\leq i \leq h-1$, the fundamental weights of $SL(h)$.
Then $L_{\lambda_i}$ is isomorphic to
$\pi^*_i\StructureSheaf{G(i,H)}(1)$, where
$\pi_i:\B\rightarrow G(i,H)$ is the natural projection. 
We get the isomorphism
$
\tau^*_{G(t,H)}(-k) \ \ \cong \ \
\pi_{t_*}\left[
L_{\lambda_1-k\lambda_t}
\right].
$
Moreover, the higher direct images 
$R^i\pi_{t_*}[L_{\lambda_1-k\lambda_t}]$ vanish, for $i>0$.
The isomorphism
\[
H^i(G(t,H),\tau^*(-k)) \ \ \cong \ \ 
H^i(\B,L_{\lambda_1-k\lambda_t})
\]
follows for all $i$. 

Set $\rho:=\sum_{i=1}^{h-1}\lambda_i$ and 
$\delta:=\rho+\lambda_1-k\lambda_t\equiv
\sum_{i=1}^h
\frac{h+1-2i}{2}\epsilon_i+\epsilon_1-k(\epsilon_1+\cdots + \epsilon_t)$, 
where the second equivalence is taken modulo 
$\epsilon_1+\cdots + \epsilon_h$. 
The Borel-Weil-Bott Theorem implies, that
$H^i(\B,L_\lambda)$ vanishes for all $i$, if and only if
$(\lambda+\rho,\alpha)=0$, for some root $\alpha$
(see \cite{bott}). 
We have
\[
(\delta,\epsilon_i) \ = \ 
\left\{
\begin{array}{lcc}
(h+1)/2 - k & \mbox{if} & i=1,
\\
(h+1)/2 - i - k & \mbox{if} & 2\leq i \leq t,
\\
(h+1)/2 - i & \mbox{if} & t < i.
\end{array}
\right.
\]
For $i<j$ we get,
\[
(\delta,\epsilon_i-\epsilon_j) \ = \ 
\left\{
\begin{array}{lcccc}
j & \mbox{if} & i=1 & \mbox{and} & 2\leq j \leq t,
\\
j-k & \mbox{if} & i=1 & \mbox{and} & t<j,
\\
j-i & \mbox{if} & 2\leq i < j \leq t,
\\
j-i-k & \mbox{if} & 2\leq i \leq t & \mbox{and} & t < j,
\\
j-i & \mbox{if} & t < i < j.
\end{array}
\right.
\]

For $k$ in the range $t < k \leq h$, we get
$(\delta,\epsilon_1-\epsilon_k)=0$. 
For $k$ in the range $1\leq k \leq h-2$, 
set $i:=\max\{2,t-k+1\}$. Then 
$2\leq i \leq t$, \ \ $t < i+k\leq h$, and
$(\delta,\epsilon_i-\epsilon_{i+k})=0$. 
Part 
\ref{lemma-item-vanishing-of-cohomology-of-twists-of-tau-dual} follows 
from the Borel-Weil-Bott Theorem.

\ref{lemma-item-vanishing-of-cohomology-of-twists-of-q})
The identification, of $G(t,H)$ with
$G(h-t,H^*)$, identifies the vector bundle $q$ with $\tau^*$. 
The vanishing follows from part 
\ref{lemma-item-vanishing-of-cohomology-of-twists-of-tau-dual}.
\EndProof

\begin{new-lemma}
\label{lemma-K-theoretic-pullback-of-a-sheaf-supported-on-a-divisor}
Let $f:X\rightarrow Y$ be a morphism of smooth and irreducible varieties, 
$D$ a divisor on $Y$, and $\tilde{D}$ its pullback to $X$. Denote by 
$e:D\hookrightarrow Y$ and $\tilde{e}:\tilde{D}\hookrightarrow X$
the closed imersions. Let $U$ be a locally free sheaf over $D$.
Then the sheaf theoretic pullback 
$f^*(e_*U)$ is equal to $\tilde{e}_*(\restricted{f}{\tilde{D}}^*U)$ and 
$f^*(e_*U)$ represents the $K$-theoretic pullback as well
(the higher torsion sheaves vanish). 
\end{new-lemma}

\section{Stratified reflections of Hilbert schemes 
with respect to $+2$ vectors}
\label{sec-stratified-elementary-trans-plus-2-vectors}

We recall in section \ref{sec-hilbert-scheme-with-a-birational-involution} 
the relationship between 
\begin{enumerate}
\item an equivalence $\Phi:D(S)\rightarrow D(S)^{op}$, obtained as the 
composition of 
a) the reflection with respect to the trivial line bundle, 
b) the duality functor, and c) the shift autoequivalence, 
and 
\item
stratified Mukai elementary transformations, 
of moduli spaces of sheaves on the K3 surface $S$. 
\end{enumerate}

Let $\sigma$ denote the isometry (\ref{eq-sigma})
of the Mukai lattice, induced by 
the functor $\Phi$. The main result of this section, 
Theorem \ref{thm-reflection-sigma-satisfies-main-conj}, 
verifies Theorem \ref{thm-Gamma-v-acts-motivicly}, for 
the isometry $\sigma$ and the Hilbert scheme $S^{[n]}$. 
Recall the Steinberg correspondence $\Z\subset S^{[n]}\times S^{[n]}$, 
associated to the self-dual stratified elementary transformation of $S^{[n]}$
(Definition \ref{def-Z-for-cotangent-bundles}).
The Theorem expresses, in addition, the cohomology class Poincare-dual to 
$\Z$ in terms of 1) the Chern classes of the universal sheaf over
$S^{[n]}$ and 2) the functor $\Phi$. 

In section \ref{sec-notation} we summarize the notation used throughout
section \ref{sec-stratified-elementary-trans-plus-2-vectors}.
In sections \ref{sec-the-structure-of-the-correspondence-in-the-plus-2-case} 
and \ref{sec-more-on-the-structure-of-the-correspondences} we describe
properties of the Steinberg correspondence $\Z$ in $S^{[n]}\times S^{[n]}$.
In sections 
\ref{sec-relating-the-two-pullbacks-of-universal-sheaf-to-Delta-t} 
and \ref{sec-the-line-bundle-A-t} we relate the two pullbacks 
of the universal sheaf to $\Z^{[1]}_t\times S$, 
via the two projections to $S^{[n]}\times S$. Above,
$\Z^{[1]}_t$ is the resolution of a component of the correspondence
$\Z$. 

Theorem \ref{thm-reflection-sigma-satisfies-main-conj} 
is proven in section 
\ref{sec-proof-of-thm-reflection-sigma}.
The proof boils down to the computation of the 
K-theoretic push-forward of $\sigma(\E_v)$, 
via the correspondence $\Z$. 
Recall, that  $\sigma(\E_v)$ is given in 
(\ref{eq-sigma-of-universal-sheaf}). The class $\sigma(\E_v)$  is 
the K-theoretic representative, of the 
transform $\Phi(\E_v)$ of the universal sheaf $\E_v$, 
via the functor $\Phi$.
The push-forward of $\sigma(\E_v)$ via $\Z$
agrees, in the K-group of $S^{[n]}\times S$, 
with the class of a universal sheaf.

The proof of Theorem \ref{thm-reflection-sigma-satisfies-main-conj} is 
more involved, than the proof of Theorem 
\ref{thm-class-of-correspondence-in-stratified-elementary-trans},
due to the following two reasons. First,
the construction of the K-theoretic representative $\sigma(\E_v)$
of $\Phi(\E_v)$,
which preserves stability, is more involved
(Section \ref{sec-relating-the-two-pullbacks-of-universal-sheaf-to-Delta-t}).
We decompose the pullback of $\sigma(\E_v)$ to 
$\Z^{[1]}\times S$, as in the proof of Theorem
\ref{thm-class-of-correspondence-in-stratified-elementary-trans}. 
Again, the contribution of most of the terms is annihilated.
In the proof of Theorem 
\ref{thm-class-of-correspondence-in-stratified-elementary-trans},
each of the irrelavent terms was annihilated
by the second pushforward to $S^{[n]}\times S$ (see
equations (\ref{eq-vanishing-of-k-theory-pushforward-of-pullback-of-F-t}) 
and (\ref{eq-vanishing-of-k-theory-pushforward-of-pullback-of-W-t-minus-I})). 
Unfortunately, this is no longer the case in the 
proof of Theorem \ref{thm-reflection-sigma-satisfies-main-conj}.
The proof, of the cancelation of the irrelavent terms, is acheived using 
the finer analysis carried out in section
\ref{sec-more-on-the-structure-of-the-correspondences} and 
\ref{sec-the-line-bundle-A-t}. 

\subsection{Hilbert schemes with a birational involution}
\label{sec-hilbert-scheme-with-a-birational-involution}

Let $u_0=(1,0,-1)$ be the $+2$ Mukai vector of the ideal sheaf of
two points. 
Denote by $\sigma := \sigma_{u_0}$ the reflection 
of the Mukai lattice with respect to $u_0$, given in (\ref{eq-sigma}).
Let $v_0$ be the $-2$ vector $(1,0,1)$ of the trivial line-bundle. 
The corresponding reflections, $\sigma_{u_0}$ and $\tau_{v_0}$,
commute and satisfy the relation (\ref{eq-relation-between-sigma-and-tau}).


\begin{example}
\label{example-geometric-reflections-by-plus-2-vectors}
{\rm 
The reflection by $u_0:=(1,0,-1)$ was studied in detail in
\cite{markman-reflections}, for Mukai vectors with
first Chern class satisfying Condition \ref{condition-minimality}.
If $\sigma_{u_0}(v)=v$, the resulting reflection is often related
to an anti-symplectic birational automorphism. 
Let us describe a few cases: 
\begin{enumerate}
\item
Let $v=(0,\LB,0)$, where $\LB$ is an effective  line-bundle of degree
$2g-2$, $g\geq 0$. Assume Condition \ref{condition-minimality}. 
Then the relative compactified Jacobian $\M(v)$ is smooth and compact. 
A point in $\M(v)$  parametrizes a stable sheaf $F$, 
with pure one-dimensional support $C_F$ in the linear system $\linsys{\LB}$. 
The generic point corresponds to the pushforward of a 
line-bundle, of zero Euler-characteristic,  
on a curve in the linear system $\linsys{\LB}$. 
The Mukai vector $v$ is fixed by both $\sigma_{u_0}$ and 
$\tau_{v_0}$. The composition $\sigma_{u_0}\tau_{v_0}$ is $-D$, where $D$
is the duality operator  $D(w):=w^\vee$. 
The reflection $\gamma_{\sigma_{u_0}}$, given in 
(\ref{eq-action-of-stabilizer}), 
is induced by the composition $\gamma_{\tau_{v_0}}\gamma_{-D}$.
The class $\gamma_{\tau_{v_0}}$ is described in Theorem
\ref{thm-class-of-correspondence-in-stratified-elementary-trans}, 
under the assumption, that the rank of $v$ is positive. 
We expect  Theorem 
\ref{thm-class-of-correspondence-in-stratified-elementary-trans} 
to extend to the rank $0$ case 
(see the conjecture in part one of the paper
\cite{markman-part-one}). The class
$\gamma_{-D}$ is induced by the composition of 

\hspace{2ex}
a) the regular automorphism $\iota$, sending a sheaf $F$ to 
   $\SheafHom(F,\omega_{C_F})$, and 

\hspace{2ex}
b) the duality involution, acting by $(-1)^i$ on $H^{2i}(\M(v))$.

\noindent
The sheaf $\SheafHom(F,\omega_{C_F})$ 
is isomorphic to $\SheafExt^1_{S}(F,\StructureSheaf{S})$ 
(\cite{markman-reflections} Lemma 3.23). 
The regularity, of the involution $\iota$ of $\M(v)$, was proven 
by Le Potier (\cite{le-potier} Theorem 5.7).  
\item
\label{example-item-birational-involution-of-hilbert-scheme}
Let $\LB$ be a very ample line-bundle of degree
$2n$ embedding $S$ in $\PP^{n+1}$. 
If $v=(1,\LB,1)$, the moduli space $\M(v)$ is isomorphic to the Hilbert scheme
$S^{[n]}$. 
The Mukai vector $v$
is fixed under $\sigma_{u_0}$ and the involution is related to 
the following birational transformation. A generic length $n$ subscheme $A$ 
spans a $\PP^{n-1}$ in $\PP^{n+1}$ intersecting $S$ along a length $2n$
subscheme $B$. The birational transformation sends $A$ to the complementary 
subscheme $B\setminus A$. When $n=1$ and $\LB$ is ample and base-point free,
$S$ is a double cover of $\PP^2$ and the involution is regular. 
The involution is regular also for $n=2$, provided 
$\LB$ generates the Picard group of $S$
\cite{tyurin-cycles-curves-surfaces,markman-reflections}. 
For $3\leq n\leq 7$, the involution is a Mukai elementary transformation along
a $\PP^3$-bundle over the $(2n-6)$-dimensional moduli space 
$\M_\LB(2,\LB,2)$. 
The $\PP^3$-bundle parametrizes length $n$ subschemes, which span a 
$\PP^{n-2}$.
For $n\geq 8$, the involution is a stratified Mukai elementary transformation
\cite{markman-reflections}. 
The smallest stratum $(S^{[n]})^\mu$, $\mu=\lfloor\sqrt{n+1}\rfloor-1$,
is a $G(\mu,2\mu+2)$ bundle over $\M(1\!+\!\mu,\LB,1\!+\!\mu)$. 
\item
\label{example-item-reflection-by-plus-2-vector-in-higher-rank}
Each of the moduli spaces $\M(r,\LB,r)$ admits an anti-symplectic birational 
involution $\iota$ (see \cite{markman-reflections} Theorem 3.21). 
In the special case, where $F\in \M(r,\LB,r)$ is locally
free and $h^1(F)=0$, 
then $\iota(F)$ is defined by
\[
0\rightarrow F^* \rightarrow 
H^0(F)^*\otimes_{\ComplexNumbers}\StructureSheaf{S} \rightarrow
\iota(F) \rightarrow 0.
\]
Consequently, we have an isomorphism 
$H^0(F)^*\rightarrow H^0(\iota(F))$, canonical up to a scalar factor. 
In particular, if $F=\iota(F)$, then $H^0(F)$ is {\em self dual}.
\end{enumerate}
}
\end{example}

Let $\LB$ a line bundle of degree $2g-2$, satisfying Condition
\ref{condition-minimality}, and
$v=(1,\LB,1)$ the Mukai vector of the $\LB$-twisted ideal
sheaf of $g-1$ points. 
We consider again the Brill-Noether stratification
(\ref{eq-brill-noether-stratification}) of $\M(v)$.
Theorem 3.2 in \cite{markman-reflections}
states, that this stratification is the top row, 
of a stratified dualizable collection (\ref{eq-diagram-of-X-v}). 
Furthermore, the symplectic varieties $X(t)$, in the diagonal entries of 
(\ref{eq-diagram-of-X-v}), are $X(0)=\M(v)$ and $X(k)=\M(v+\vec{k})$,
where $\vec{k}$ is the Mukai vector $(k,0,k)$ of the trivial rank $k$ vector
bundle. The collection is self-dual;
the dual collection $X'(k)^t$ is isomorphic to $X(k)^t$
(see Remark \ref{rem-self-duality-in-the-plus-2-case}). 
We have the grassmannian bundles
\begin{equation}
\label{eq-grassmannian-bundle}
f_t \ : \ [\M(v)^t\setminus\M(v)^{t+1}] \ \rightarrow \ 
[\M(v+\vec{t})\setminus \M(v+\vec{t})^1].
\end{equation}
They extend to the grassmannian fibrations
$B^{[t+1]}\M(v)^t\rightarrow B^{[1]}\M(v+\vec{t})$, mentioned in
(\ref{eq-f-k-t}) for a general stratified elementary transformation,
and described in part 
\ref{prop-item-morphism-is-a-classifying-morphism-of-F-t} 
of Proposition \ref{prop-decomposition-of-complex-over-interated-blow-up}.

Using the self-duality, the correspondence 
\begin{eqnarray}
\label{eq-correspondence-inducing-involution-of-hilbert-scheme}
\Z & \subset & \M(v) \times \M(v)
\\
\nonumber
\Z & := & \sum_{i=0}^\mu \Z_i,
\end{eqnarray}
introduced in Definition \ref{def-Z-for-cotangent-bundles}
part 
\ref{correspondence-fiber-products-of-dual-grassmanian-bundles},
can be described as follows. 
$\Z_0$ is the closure of the graph of the involution 
\begin{equation}
\label{eq-involution-iota-on-open-stratum}
\iota \ : \ [\M(v)\setminus\M(v)^1] \ \LongIsomRightArrow \ 
[\M(v)\setminus\M(v)^1],
\end{equation} 
described in Example \ref{example-geometric-reflections-by-plus-2-vectors} 
part \ref{example-item-birational-involution-of-hilbert-scheme}.
$\Z_t$ is the closure of the $\iota$-twist of the 
fiber product of the grassmannian bundles 
(\ref{eq-closure-of-fiber-product-of-grassmannian-bundles}). 
By {\em $\iota$-twist} we mean, that we consider the inverse image in
$\M(v)^t\times\M(v)^t$, via $f_t\times f_t$, 
of the graph of the regular involution
\begin{equation}
\label{eq-sigma-involution-on-lower-strata}
\iota \ : \ [\M(v+\vec{t})\setminus \M(v+\vec{t})^1] \ \ 
\IsomRightArrow \ \ [\M(v+\vec{t})\setminus \M(v+\vec{t})^1] 
\end{equation}
described in Example 
\ref{example-geometric-reflections-by-plus-2-vectors} part 
\ref{example-item-reflection-by-plus-2-vector-in-higher-rank}.
We proved in \cite{markman-reflections}, that $\Z$ induces an involution
of the cohomology ring of $\M(v)$. 
Our goal in this section is to prove:

\begin{thm}
\label{thm-reflection-sigma-satisfies-main-conj}
\begin{enumerate}
\item
\label{thm-item-class-of-the-correspondence-Z-of-sigma}
The class of the correspondence
(\ref{eq-correspondence-inducing-involution-of-hilbert-scheme}),
in the cohomology ring  $H^*(\M(v)\times \M(v),\Integers)_{\rm free}$, 
is equal to the class 
of the composition $D_\M\circ\gamma_\sigma(\E_v,\E_v)$.
\item
Parts \ref{thm-item-gamma-g-is-a-ring-isomorphism} and 
\ref{thm-item-normalized-universal-chern-character-is-invariant}
of Theorem \ref{thm-Gamma-v-acts-motivicly} hold for the reflection $\sigma$
and the moduli space $\M_H(v)$, 
where the polarization $H$ and
the line-bundle $\LB=c_1(v)$ satisfy condition 
\ref{condition-minimality}. 
\end{enumerate}
\end{thm}

\begin{rem}
\label{rem-K-theory-level-plus-2-case}
{\rm
The analogue of Remark 
\ref{rem-K-theory-level-minus-2-case}
holds for Theorem \ref{thm-reflection-sigma-satisfies-main-conj} as well.
}
\end{rem}

\begin{rem}
\label{rem-self-duality-in-the-plus-2-case}
{\rm
Let $\W$ be the stratified transform of $\M(v)$, with respect to the 
Brill-Neother stratification, in the sense of \cite{markman-reflections}
Theorem 2.4 (and section \ref{sec-stratified-elemetary-transformations}). 
The involution of $H^*(\M(v),\Integers)$ arises from {\em two}
isomorphisms of the cohomology rings of 
$\M(v)$ and $\W$. One is the correspondence in $\M(v)\times \W$,
described in part 
\ref{correspondence-fiber-products-of-dual-grassmanian-bundles} of 
definition \ref{def-Z-for-cotangent-bundles}.
The other is a regular isomorphism $\eta:\M(v)\IsomRightArrow \W$. 
The above correspondence $\Z$, given in 
(\ref{eq-correspondence-inducing-involution-of-hilbert-scheme}),
is the image, via $id\times \eta$, of the standard correspondence given 
in part \ref{correspondence-fiber-products-of-dual-grassmanian-bundles} of 
definition \ref{def-Z-for-cotangent-bundles}.
Let us briefly describe the isomorphism $\eta$. 
The transform $\W$ is endowed with a stratification and the open stratum
$\W\setminus \W^1$ is naturally equal to 
$\M(v)\setminus \M(v)^1$, being the complement of the center of the 
birational isomorphism. 
The lower strata are the  Grassmannian bundles dual to
(\ref{eq-grassmannian-bundle}).
The identification 
$[\W\setminus \W^1]=[\M(v)\setminus \M(v)^1]$ does {\em not}
extend to an isomorphism, but its composition with the involution
(\ref{eq-involution-iota-on-open-stratum}) does extend
(\cite{markman-reflections} Proposition 3.6). 
The point is that (\ref{eq-involution-iota-on-open-stratum})
is compatible with the involutions
(\ref{eq-sigma-involution-on-lower-strata}), 
which pulls back the Grassmannian 
bundle (\ref{eq-grassmannian-bundle}) to its {\em dual} bundle
(for a precise description of this compatibility, see the diagram
(\ref{eq-commutative-diagram-of-two-involutions-g-iota})). 
}
\end{rem}

\subsection{Notation}
\label{sec-notation}

We summerize here some of the notation, that will be 
repeatedly used throughout the rest of section 
\ref{sec-stratified-elementary-trans-plus-2-vectors}. 
This summary is intended for reference; the notation is introduced 
again when first mentioned. 
We continue to use the notation summerized in
section \ref{sec-notation-minus-2-vector}.
Let $\sigma_{u_0}$ be the reflection
(\ref{eq-sigma}), of the Mukai lattice, with respect to
a $+2$-vector $u_0$. We set $\sigma:=\sigma_{u_0}$, where $u_0=(1,0,-1)$.
We denote by $\sigma$ also the related operation, on families of sheaves, 
defined in lemma 
\ref{lemma-lazarsfeld-reflection-of-a-family-of-sheaves}.

The birational involution $\iota:\M(r,\LB,r)\rightarrow \M(r,\LB,r)$
is introduced in Example
\ref{example-geometric-reflections-by-plus-2-vectors} part 
\ref{example-item-reflection-by-plus-2-vector-in-higher-rank}. 
The regular involution $g:B^{[1]}\M(v)\rightarrow B^{[1]}\M(v)$ 
is the resolution of $\iota$, introduced in diagram
(\ref{eq-commutative-diagram-of-two-involutions-g-iota}).
The morphism $\pi_1:B^{[1]}\M(v)\times S \rightarrow B^{[1]}\M(v)$
is the projection. 

Each moduli space $\M(v)$, considered in this section, 
is part of a {\em self-dual} startified collection. 
The correspondence $\Z$, in Definition
\ref{def-Z-for-cotangent-bundles}
part \ref{correspondence-fiber-products-of-dual-grassmanian-bundles},
is thus contained in the cartesian square $\M(v)\times \M(v)$. 
The description of $\Z$, incorporating the self-duality, is 
given in (\ref{eq-correspondence-inducing-involution-of-hilbert-scheme}). 
The normal crossing model $\Z^{[1]}$ of $\Z$ is described in
section \ref{sec-the-structure-of-the-correspondence-in-the-plus-2-case}.
$\Delta_t^{[1]}$ denotes an irreducible component of $\Z^{[1]}$
as in Lemma \ref{lemma-structure-of-correspondence-plus-2-case}. 
The twisted fiber product $\Delta_t^{[t+1]}$ is defined in 
Lemma \ref{lemma-structure-of-correspondence-plus-2-case} as well.
$I_t$ denotes the incidence divisor in $\Delta_t^{[t+1]}$, defined in
section \ref{sec-incidence-divisor}.
The line bundle $A_t$ over $\Delta_t^{[t+1]}$ is introduced in
Corollary \ref{cor-lazarsfeld-reflection-of-universal-sheaf-over-Delta-t}.

We will denote by $\tilde{\delta}_{t,ij}$ the morphism from
$\Delta_t^{[1]}\times S$ to the product of factors of 
$\M(v)\times S\times \M(v)$. 
The morphism 
$\tilde{\delta}_{t\cap i,3}:
[\Delta_t^{[1]}\cap \Delta_i^{[1]}]\times S \rightarrow 
B^{[i+1]}\M(v)^{t\cap i}$
is the composition of the inclusion in $\Delta_i^{[1]}\times S$, the
projection to $\Delta_i^{[i+1]}$, 
followed by the second projection from
$\Delta_i^{[i+1]}$ to $B^{[i+1]}\M(v)^{i}$. 
For $i\geq t$, $\Delta_{t,i}^{[t+1]}$ is the subvariety of 
$\Delta_t^{[t+1]}$ defined in corollary 
\ref{cor-lazarsfeld-reflection-of-universal-sheaf-over-Delta-t}.
The morphism 
$\delta_{t\cap i,1}:\Delta_{t,i}^{[t+1]}\rightarrow B^{[i+1]}\M(v)^{t\cap i}$
is defined in corollary 
\ref{cor-lazarsfeld-reflection-of-universal-sheaf-over-Delta-t}
as well.

\subsection{The structure of the correspondence in the $+2$ vector case}
\label{sec-the-structure-of-the-correspondence-in-the-plus-2-case}

We define 
the top iterated blow-up $\Z^{[1]}$, of the correspondence $\Z$, in
a way analogous to the definition 
in Section \ref{sec-resolution-of-correspondences} equation
(\ref{eq-definition-of-top-iterated-blow-up-of-Z}). 
We observed, in Remark \ref{rem-self-duality-in-the-plus-2-case},
that $\M(v)$ is isomorphic to its stratified elementary transform. 
Hence, components of the correspondence $\Z^{[1]}$ can be described, in
terms of $\M(v)$ directly. Due to the subtlety of the self-duality
in Remark \ref{rem-self-duality-in-the-plus-2-case}, we include 
this explicit description in the following lemma. 

\begin{new-lemma}
\label{lemma-structure-of-correspondence-plus-2-case}
\begin{enumerate}
\item
\label{lemma-item-iota-extends-to-a-regular-involution}
The birational isomorphism $\iota:\M(v+\vec{t})\rightarrow \M(v+\vec{t})$,
given in (\ref{eq-sigma-involution-on-lower-strata}), 
extends to a regular involution 
$\iota: B^{[1]}\M(v+\vec{t})\rightarrow B^{[1]}\M(v+\vec{t})$.
\item
\label{lemma-item-component-is-a-fiber-product-of-bundle-with-its-pullback}
The component $\Delta^{[1]}_t$ of $\Z^{[1]}$ is the top iterated blow-up 
of the fiber product $\Delta^{[t+1]}_t$ of the Grassmannian bundle
$B^{[t+1]}\M(v)^t$ with its pullback $\iota^*B^{[t+1]}\M(v)^t$
\[
\Delta^{[t+1]}_t \ := \ B^{[t+1]}\M(v)^t\times_{B^{[1]}\M(v+\vec{t})} 
\iota^*B^{[t+1]}\M(v)^t,
\]
where $\iota$ is the involution of $B^{[1]}\M(v+\vec{t})$ 
in part \ref{lemma-item-iota-extends-to-a-regular-involution}. 
\item
\label{lemma-item-components-of-correspondence-intersect-as-flag-variety}
The Grassmannian bundles 
$B^{[t+1]}\M(v)^t\rightarrow B^{[1]}\M(v+\vec{t})$ and 
$\iota^*B^{[t+1]}\M(v)^t\rightarrow B^{[1]}\M(v+\vec{t})$ in part 
\ref{lemma-item-component-is-a-fiber-product-of-bundle-with-its-pullback} 
are the dual Grassmannian bundles 
$G(t,W)$ and $G(t,W^*)$, where 
$\PP{W}$ is the $\PP^{\chi(v+\vec{t})-1}$-bundle over 
$\M(v+\vec{t})$ in Proposition 
\ref{prop-decomposition-of-complex-over-interated-blow-up}
($\PP{W}$ comes from a vector bundle $W$, if there is a universal sheaf over
$\M(v+\vec{t})\times S$).
The intersection $\Delta^{[1]}_t\cap\Delta^{[1]}_0$ maps 
onto the relative flag subvariety $Flag(t,\chi(v)+t,W)$ 
of $\Delta^{[t+1]}_t$, under the identification of the latter with
the (untwisted) fiber product 
$G(t,W)\times_{B^{[1]}\M(v+\vec{t})}G(t,W^*)$.
\item
\label{lemma-item-pull-back-of-BN-divisors-to-correspondence-in-plus-2-case}
The composition
$\Delta_t^{[1]} \ \rightarrow \ \Delta_t^{[t+1]} \ \rightarrow \ 
B^{[t+1]}\M(v)^t \ \hookrightarrow \ B^{[t+1]}\M(v)$
pulls back (scheme theoretically) 
the divisor $B^{[t+1]}\M(v)^i$, \ $i\geq t+1$, 
to the divisor $\Delta_t^{[1]}\cap\Delta_i^{[1]}$.
\item
\label{lemma-item-normal-bundle-of-Delta-t}
The line-bundle 
$
\StructureSheaf{\Z^{[1]}}(\sum_{t=0}^\mu\Delta_t^{[1]}) 
$
is the trivial line bundle over $\Z^{[1]}$.
\end{enumerate}
\end{new-lemma}

{\bf Proof:} Part \ref{lemma-item-iota-extends-to-a-regular-involution}
is proven in \cite{markman-reflections} Proposition 3.26.
Part 
\ref{lemma-item-component-is-a-fiber-product-of-bundle-with-its-pullback} 
is the analogue, for $+2$ reflections, of lemma 
\ref{lemma-two-descriptions-of-the-components-of-the-correspondence}. 
Its proof is identical. 
Part \ref{lemma-item-components-of-correspondence-intersect-as-flag-variety} 
follows from Theorem 3.2 in \cite{markman-reflections}. 
The proof of
part \ref{lemma-item-pull-back-of-BN-divisors-to-correspondence-in-plus-2-case}
is similar to that of part
\ref{lemma-item-pull-back-of-BN-divisors-to-correspondence} 
of lemma \ref{lemma-structure-of-correspondence}. 
Part 
\ref{lemma-item-normal-bundle-of-Delta-t}
is a consequence of the construction of $\Z^{[1]}$. 
We have two non-isomorphic families $\X\rightarrow B$ and $\X'\rightarrow B$ 
of deformations of $\M(v)$, over a one-parameter base $B$. 
The iterated blow-ups $B^{[1]}\X$ and $B^{[1]}\X'$ are isomorphic and
$\Z^{[1]}$ is the special fiber of $B^{[1]}\Y\rightarrow B$, in 
the graph $B^{[1]}\Y$,
of the isomorphism (\ref{eq-top-iterated-blow-up-of-Y}). 
The line bundle 
$
\StructureSheaf{\Z^{[1]}}(\sum_{t=0}^\mu\Delta_t^{[1]}) 
$
is the normal line-bundle to the fiber of
$B^{[1]}\Y$.
Now, $B^{[1]}\Y$ is isomorphic to $B^{[1]}\X$ and 
the normal bundle to the special fiber of $B^{[1]}\X$ is trivial. 
\EndProof


\subsection{Reflection by a $+2$ vector: a relative construction}
\label{sec-relating-the-two-pullbacks-of-universal-sheaf-to-Delta-t}

Let $F$ be a stable sheaf on $S$, with $\chi(F)\geq 0$. Given 
a $\chi(F)$-dimensional subspace $W\subset H^0(F)$,
we get the complex $W\otimes \StructureSheaf{S}\rightarrow F$,
of the evaluation homomorphism, with $F$ in degree $1$.
The Mukai vector of the complex is equal to $-\tau(v(F))$. 
Recall, that $\sigma(v)=-\tau(v)^\vee$. 
Theorem 3.21 in \cite{markman-reflections}
involves the construction of a complex 
$W^*\otimes \StructureSheaf{S}\rightarrow \sigma(F,W)$, 
with $\sigma(F,W)$ a stable sheaf with Mukai vector
$\sigma(v(F))$. The sheaf $\sigma(F,W)$ is equivalent to
$[W\otimes \StructureSheaf{S}\rightarrow F]^\vee$
in the derived category of $S$.
The proof of Theorem \ref{thm-reflection-sigma-satisfies-main-conj}
requires the following relative version of the construction. 

\begin{new-lemma}
\label{lemma-lazarsfeld-reflection-of-a-family-of-sheaves}
Let $T$ be a smooth quasi-projective variety and $\F$ a 
sheaf over $T\times S$, flat over $T$,  representing a family
of $H$-stable sheaves over $S$, with Mukai vector $v'=(r,\LB,s)$. 
Assume, that $r$ and $s$ are $\geq 0$ 
and $\LB$ satisfies Condition \ref{condition-minimality}.
Assume, that a vector bundle $W$ of rank $\chi(v')$
is given over $T$, together with a surjective homomorphism
\begin{equation}
\label{eq-a-surjective-homomorphism-ev-dual}
ev^* \ : \ \RelExt^2_{\pi_{T,*}}(\F,\StructureSheaf{T\times S}) 
\ \ \longrightarrow \ \ W^*.
\end{equation}
Then there exists a natural sheaf $\sigma(\F,W)$, flat over $T$,  
representing a family of $H$-stable sheaves over $S$, 
with Mukai vector $\sigma(v')=(s,\LB,r)$. 
In addition, there is a natural surjective homomorphism
\[
\sigma(ev)^* \ : \ 
\RelExt^2_{\pi_{T,*}}(\sigma(\F,W),\StructureSheaf{T\times S}) 
\ \ \longrightarrow \ \ W.
\]
Furthermore, $\sigma(\F,W)$ fits in the exact sequence
\begin{equation}
\label{eq-relation-between-F-and-its-reflection-wrt-sigma}
0\rightarrow
\SheafHom(\F,\StructureSheaf{}) \LongRightArrowOf{ev^*}
\pi_T^*W^*\LongRightArrowOf{\sigma(ev)}
\sigma(\F,W) \longrightarrow
\SheafExt^1_{T\times S}(\F,\StructureSheaf{})\rightarrow 0
\end{equation}
and $\SheafExt^i_{T\times S}(\F,\StructureSheaf{})$ vanishes for $i\geq 2$.
Finally, the construction is involutive; the complex
(\ref{eq-a-surjective-homomorphism-ev-dual}) is naturally isomorphic to 
\[
\sigma(\sigma(ev))^* \ : \ 
\RelExt^2_{\pi_{T,*}}(\sigma(\sigma(\F,W),W^*),\StructureSheaf{T\times S}) 
\ \ \longrightarrow \ \ W^*.
\]
\end{new-lemma}

{\em Note: The exact sequence 
(\ref{eq-relation-between-F-and-its-reflection-wrt-sigma}) implies,
that $\sigma(\F,W)$ is equivalent to 
the dual of complex $\pi_T^*W\rightarrow \F$, 
in the derived category of $T\times S$. 
The sequence (\ref{eq-relation-between-F-and-its-reflection-wrt-sigma}) 
is a relative 
version of the one given in (103) in \cite{markman-reflections}. 
}

\medskip
\noindent
{\bf Proof:} 
The general idea of the proof was communicated to the author by K. Yoshioka.
It leads to a simplification of the proof of Theorem 3.21 in
\cite{markman-reflections}.

Assume first, that $r\geq 1$. 
We prove first, that $\F$ has projective dimension $\leq 1$. 
This is equivalent to the vanishing $\SheafExt^i_{T\times S}(\F,E)=0$, 
for every sheaf $E$ over $T\times S$ and for $i\geq 2$.
Choose, locally over $T$, a rank $r-1$ subbundle $V$ of $W$. 
The evaluation homomorphism $\pi_T^*V\rightarrow \F$ is injective 
(\cite{markman-reflections} Lemma 3.7 part 3).
Moreover, the quotient $Q$ is a family, flat over $T$, of stable sheaves on $S$
with Mukai vector $(1,\LB,s+1-r)$ (\cite{markman-reflections} Theorem 3.15). 
$Q$ is the tensor product of the ideal sheaf 
of a codimension $2$ subscheme $\Y\subset [T\times S]$, 
with a line bundle on $T\times S$. 
Since $T\times S$ is smooth, the ideal sheaf of $\Y$ has projective dimension 
$1$ (or zero if $\Y$ is empty). The vanishing of the sheaf 
$\SheafExt^i_{T\times S}(Q,E)$, $i\geq 2$, follows. 
We have a short exact sequence 
\[
0\rightarrow \pi_T^*V \rightarrow \F \rightarrow Q \rightarrow 0.
\]
The vanishing
$\SheafExt^i_{T\times S}(\F,E)=0$, for $i\geq 2$, 
follows from the corresponding long exact sequence of extension sheaves. 

Let 
\[
0\rightarrow K_1'\rightarrow W_1' \rightarrow \F \rightarrow 0
\]
be a locally free resolution of $\F$. 
Set $W_1:=W_1'\oplus \pi_T^*W$ and let $K_1$ be the kernel
of the natural homomorphism from $W_1$ to $\F$
\[
0\rightarrow K_1 \rightarrow W_1\rightarrow \F \rightarrow 0.
\]

Let $\widetilde{\F}$ be the image of the evaluation homomorphism
$ev : \pi_T^*W\rightarrow\F$. 
We get the following commutative diagram, with short exact rows and 
columns, defining $K_0$, $K_2$, and $W_2$.
\[
\begin{array}{ccccc}
K_0 & \rightarrow & \pi_T^*W & \rightarrow & \widetilde{\F}
\\
\downarrow & & \downarrow & & \downarrow 
\\
K_1 & \rightarrow & W_1 & \rightarrow & \F
\\
\downarrow & & \downarrow & & \downarrow 
\\
K_2 & \rightarrow & W_2 & \rightarrow & \F/\widetilde{\F}.
\end{array}
\]
The sheaves $K_1$,$\pi_T^*W$, $W_1$, and $W_2$ are locally free, the 
latter being $W_1'$. 
We know that $\F/\widetilde{\F}$ 
has cokernel with support of codimension $\geq 2$ (\cite{markman-reflections}
Lemma 3.7 part 4). 
Thus, the natural homomorphism
$W_2^*\rightarrow \SheafHom(K_2,\StructureSheaf{T\times S})$
is an isomorphism and 
$\SheafExt^1(\F/\widetilde{\F},\StructureSheaf{T\times S})$ vanishes.
It follows, that 
$\SheafHom(\F,\StructureSheaf{T\times S})\rightarrow
\SheafHom(\widetilde{\F},\StructureSheaf{T\times S})$
is an isomorphism and 
$e:\SheafExt^1_{T\times S}(\F,\StructureSheaf{T\times S})\rightarrow
\SheafExt^1_{T\times S}(\widetilde{\F},\StructureSheaf{T\times S})$
is injective.
We get the following commutative diagram with exact rows 
and columns
\[
\begin{array}{ccccccccccc}
& 0 & \rightarrow & W_2^* & \IsomRightArrow &
\SheafHom(K_2,\StructureSheaf{}) & \rightarrow & 0 
\\
 & \downarrow & & \downarrow & & \downarrow & & \downarrow
\\
0  \rightarrow & \SheafHom(\F,\StructureSheaf{}) & \rightarrow &
W_1^* & \rightarrow & K_1^* & \rightarrow & 
\SheafExt^1(\F,\StructureSheaf{}) & \rightarrow 0
\\
 & \downarrow & & \downarrow & & \downarrow & & \hspace{1ex} \ \downarrow \ e
\\
0 \rightarrow & \SheafHom(\F,\StructureSheaf{}) & \rightarrow &
\pi_T^*W^* & \rightarrow & K_0^* & \rightarrow & 
\SheafExt^1(\widetilde{\F},\StructureSheaf{}) & \rightarrow 0
\\
 & \downarrow & & \downarrow & & \downarrow & & \downarrow 
\\
& 0 & & 0 &\rightarrow & \SheafExt^1(K_2,\StructureSheaf{}) &
\IsomRightArrow & 
\SheafExt^2(\F/\widetilde{\F},\StructureSheaf{}).
\end{array}
\]
Consider the third row above as a $2$-extension. We see, that its pullback, 
via the injective homomorphism $e$, is the $2$-extension
\[
0\rightarrow \SheafHom_{T\times S}(\F,\StructureSheaf{}) \rightarrow
\pi_T^*W^* \rightarrow K_1^*/\SheafHom(K_2,\StructureSheaf{}) \rightarrow 
\SheafExt^1_{T\times S}(\F,\StructureSheaf{}) \rightarrow 0
\]
We set $\sigma(\F,W):=K_1^*/\SheafHom(K_2,\StructureSheaf{})$. 
It is independent of the choice of the locally free resolution of $\F$,
being the third term in the pulled back $2$-extension.
It is a quotient of two locally free sheaves. The homomorphism
$\SheafHom(K_2,\StructureSheaf{})\rightarrow K_1^*$
restricts to an injective homomorphism over $t\times S$, for every
closed point $t\in T$. Hence, $\sigma(\F,W)$ is flat over $T$
(see \cite{matsumura} application 2 page 150 
and Lemma 3.11 in \cite{markman-reflections}).
The stability of the sheaves parametrized by $\sigma(\F,W)$ was proven in
the first step of the proof of Theorem 3.21 in \cite{markman-reflections}.
The involutivity of the construction is straightforward to check.

The case  $r=s=0$ of the lemma was proven in \cite{markman-reflections}
Theorem 3.21, following Le Potier. One 
simply sets  
$\sigma(\F,W):=\SheafExt^1_{T\times S}(\F,\StructureSheaf{})$.
When $r=0$ and $s\geq 0$, we construct the sheaf $\sigma(\F,W)$ as an 
extension of 
$\SheafExt^1_{T\times S}(\F,\StructureSheaf{})$ by $\pi_T^*W^*$. 
Indeed, the exact sequence 
(\ref{eq-relation-between-F-and-its-reflection-wrt-sigma})
degenerates to a short exact sequence, since 
$\SheafHom(\F,\StructureSheaf{})=0$. 
In that case, 
$\SheafExt^1_{T\times S}(\F,\StructureSheaf{})$ is a family of
stable sheaves, with Mukai vector $-v^\vee$. 
Note, that $\sigma(v)=\tau(-v^\vee)$. 
The extension $\sigma(\F,W)$, of 
$\SheafExt^1_{T\times S}(\F,\StructureSheaf{})$ by $\pi_T^*W^*$, 
is constructed in the proof of Theorem 
3.15 in \cite{markman-reflections}. 
\EndProof

\medskip
Recall, that the direct image sheaf
$\pi_{1_*}\beta^*\E_v$ of the universal sheaf $p_*\E_v$, 
is a vector bundle $W_0$ on $B^{[1]}\M(v)$ 
(Proposition \ref{prop-decomposition-of-complex-over-interated-blow-up}
part \ref{prop-item-U-and-W-are-direct-images}). 
More generally, we have the family $\F_t$ over $B^{[t+1]}\M(v)^t$, 
fitting in the exact sequence 
(\ref{eq-the-tautological-extension-of-families-over-BN-locus}).
Its direct image is the vector bundle $W_{v+\vec{t}}$ over $B^{[t+1]}\M(v)^t$
in Lemma \ref{lemma-W-v-plus-t-is-a-direct-image}. $W_{v+\vec{t}}$ fits in 
the extension (\ref{eq-W-v-plus-t-as-an-extension}), 
of a twist of the vector bundle $W_t$ by $U_t$. 
The vector bundles $U_t$ and $W_t$ over $B^{[t+1]}\M(v)^t$ were 
introduced in Proposition
\ref{prop-decomposition-of-complex-over-interated-blow-up}. 
The vector bundle $W_{v+\vec{t}}$ parametrizes $\chi(v)+2t$-dimensional 
subspaces of global sections of the sheaves parametrized by $\F_t$.
The pair $(\F_t,W_{v+\vec{t}})$ satisfies the hypothesis of 
Lemma \ref{lemma-lazarsfeld-reflection-of-a-family-of-sheaves}. 
The families $\F_t$ and $\sigma(\F_t):=\sigma(\F_t,W_{v+\vec{t}})$
induce two classifying morphisms from $B^{[t+1]}\M(v)^t$ onto 
$\M(v+\vec{t})$. Both morphisms factor through $B^{[1]}\M(v+\vec{t})$.
Let $f_t: B^{[t+1]}\M(v)^t\rightarrow B^{[1]}\M(v+\vec{t})$ be the
one associated to $\F_t$. Then the morphism associated to $\sigma(\F_t)$ is 
the composition morphism $\iota \circ f_t$,
where $\iota$ is the involution of $B^{[1]}\M(v+\vec{t})$ in
part \ref{lemma-item-iota-extends-to-a-regular-involution}
of lemma \ref{lemma-structure-of-correspondence-plus-2-case} 
(see Corollary 
\ref{cor-lazarsfeld-reflection-of-universal-sheaf-over-Delta-t} below).

Let $\delta_1$ and $\delta_3$ be the two projections 
from $\Delta^{[t+1]}_t\times S$ to $B^{[t+1]}\M(v)^t$ and 
$\delta_2$ the projection to $S$. 
Denote by $\delta_{12}$ and $\delta_{23}$ the two projections
from $\Delta^{[t+1]}_t\times S$ to $B^{[t+1]}\M(v)^t\times S$. 
The following corollary relates the two pullbacks
$\delta_{12}^*\F_t$ and $\delta_{23}^*\F_t$.

\begin{cor}
\label{cor-lazarsfeld-reflection-of-universal-sheaf-over-Delta-t}
\begin{enumerate}
\item
\label{cor-item-A-t-exists}
The sheaf $\sigma(\delta_{12}^*\F_t,\delta_1^*W_{v+\vec{t}})$
is isomorphic to $\delta_{23}^*\F_t\otimes A_t$, 
for some line bundle $A_t$ on $\Delta^{[t+1]}_t$ 
satisfying
\begin{equation}
\label{eq-A-t-is-symmetric}
\delta_3^*W_{v+\vec{t}}\otimes A_t 
\ \ \cong \ \  
\delta_1^*W_{v+\vec{t}}^*.
\end{equation} 
Consequently, we 
have the following exact sequence of sheaves on $\Delta^{[t+1]}_t\times S$
\[
0\rightarrow 
\SheafHom_{\Delta^{[t+1]}_t\times S}(\delta_{12}^*\F_t,\StructureSheaf{}) 
\LongRightArrowOf{ev^*} \delta_1^*W_{v+\vec{t}}^*\rightarrow 
\delta_{23}^*\F_t\otimes A_t \rightarrow 
\SheafExt^1_{\Delta^{[t+1]}_t\times S}(\delta_{12}^*\F_t,\StructureSheaf{})
\rightarrow 0.
\]
Moreover, 
$\SheafExt^i_{\Delta^{[t+1]}_t\times S}(\delta_{12}^*\F_t,\StructureSheaf{})$
vanishes, for $i\geq 2$. 
\item
\label
{cor-item-first-pullback-of-E-v-dual-minus-W-t-dual-is-a-second-pullback}
The following equivalence holds in the K-group of 
$\Delta_t^{[t+1]}\times S$
\begin{eqnarray*}
\delta_{12}^!\beta^!\left[p^!p_!\E_v - \E_v
\restricted{\right]}{\M(v)^t\times S}^\vee 
& \equiv &
\delta_{23}^!\F_t\otimes A_t\otimes 
\StructureSheaf{\Delta_t^{[t+1]}}\left(
-\sum_{i=t+1}^\mu\Delta_{t,i}^{[t+1]}\right) 
\\
& & -\left(
\sum_{i=t}^\mu e_{i,!}\left[
\delta_{t\cap i,1}^* U_{\restricted{i}{B^{[i]}\M(v)^{t\cap i}}}
\right]
\right)^\vee,
\end{eqnarray*}
where, for $i\geq t$, we let 
$\Delta_{t,i}^{[t+1]}$ be the pullback to $\Delta_t^{[t+1]}$
of the divisor $B^{[t+1]}\M(v)^{t\cap i}$ in $B^{[t+1]}\M(v)^t$
(the pullback is the same, regardless which of $\delta_1$ or $\delta_2$ 
we use, by
lemma \ref{lemma-structure-of-correspondence-plus-2-case} part 
\ref{lemma-item-pull-back-of-BN-divisors-to-correspondence-in-plus-2-case}). 
We set $e_i : \Delta_{t,i}^{[t+1]} \hookrightarrow 
\Delta_t^{[t+1]}$ to be the natural embedding and
$\delta_{t\cap i,1}$ is the composition of 
the restriction 
$\Delta_{t,i}^{[t+1]}\rightarrow B^{[t+1]}\M(v)^{t\cap i}$ of $\delta_1$,
the inclusion $B^{[t+1]}\M(v)^{t\cap i}\hookrightarrow B^{[t+1]}\M(v)^i$,
and the blow-down $B^{[t+1]}\M(v)^i\rightarrow B^{[i+1]}\M(v)^i$.
\end{enumerate}
\end{cor}

\medskip
\noindent
{\bf Proof:} 
Part \ref{cor-item-A-t-exists} of the corollary, 
without the isomorphism (\ref{eq-A-t-is-symmetric}), 
follows from 1) 
the definition of $\Delta^{[t+1]}_t$ as an $\iota$-twisted
fiber product over $B^{[1]}\M(v+\vec{t})$, 
2) the construction of the order $2$ automorphism $\iota$ of
$B^{[1]}\M(v+\vec{t})$ in Theorem 3.21 of \cite{markman-reflections}, 
and 3) Lemma \ref{lemma-lazarsfeld-reflection-of-a-family-of-sheaves}. 

The isomorphism (\ref{eq-A-t-is-symmetric}) is proven as follows. 
The vector bundle $W_{v+\vec{t}}$ is characterized as the 
direct image of $\F_t$ via the projection from $B^{[t+1]}\M(v)^t\times S$ to 
$B^{[t+1]}\M(v)^t$ (see Lemma \ref{lemma-W-v-plus-t-is-a-direct-image}). 
Similarly, $W_{v+\vec{t}}^*$ is isomorphic to the 
direct image of $\sigma(\F_t,W_{v+\vec{t}})$. 
The isomorphism 
$\sigma(\delta_{12}^*\F_t,\delta_1^*W_{v+\vec{t}})\cong
\delta_{23}^*\F_t\otimes A_t$ implies, that $\delta_1^*W_{v+\vec{t}}^*$
is isomorphic to 
$\delta_3^*W_{v+\vec{t}}\otimes A_t$.

Part \ref{cor-item-A-t-exists} 
implies the equivalence
\[
\delta_{1}^!W_{v+\vec{t}}^\vee-\delta_{12}^!\F_t^\vee \ \ \equiv \ \
\delta_{23}^!\F_t\otimes A_t. 
\]
The exact sequences 
(\ref{eq-the-tautological-extension-of-families-over-BN-locus}) and
(\ref{eq-W-v-plus-t-as-an-extension})
imply the equivalence over $B^{[t+1]}\M(v)^t\times S$
\[
[\pi_1^!W_t-\beta^*(\E_v\restricted{)}{\M(v)^t\times S}]\left(
-\sum_{i=t+1}^\mu B^{[t+1]}\M(v)^{t\cap i}\right) \ \ \equiv \ \ 
\pi_1^!W_{v+\vec{t}}-\F_t,
\]
where $\pi_1$ is the projection on the first factor.
We get the equivalence
\[
\delta_{12}^![\pi_1^!W_t-\beta^*(\E_v\restricted{)}{\M(v)^t\times S}]^\vee
\ \ \equiv \ \ 
\delta_{23}^!\F_t\otimes A_t\otimes 
\StructureSheaf{\Delta_t^{[t+1]}}\left(
-\sum_{i=t+1}^\mu\Delta_{t,i}^{[t+1]}\right). 
\]

Part
\ref{cor-item-first-pullback-of-E-v-dual-minus-W-t-dual-is-a-second-pullback}
follows from the equivalence, in the K-group of
$\Delta_t^{[t+1]}$: 
\[
\delta_{1}^! \beta^! [p_!\E_v\restricted{]}{\M(v)^t} \ \ \equiv \ \ 
\delta_{1}^*W_t - \left(
\sum_{i=t}^\mu e_{i,!}\left[
\delta_{t\cap i,1}^* U_{\restricted{i}{B^{[i]}\M(v)^{t\cap i}}}
\right]
\right).
\] 
Lemma \ref{lemma-K-theoretic-pullback-of-a-sheaf-supported-on-a-divisor} 
and lemma 
\ref{lemma-structure-of-correspondence-plus-2-case} part 
\ref{lemma-item-pull-back-of-BN-divisors-to-correspondence-in-plus-2-case}
translate the above equivalence to the pullback, via $\delta_1$,
of the following equivalence over $B^{[t+1]}\M(v)^t$, 
generalizing equation 
(\ref{eq-decomposition-of-K-pushforward-of-universal-sheaf-over-B-1-M-v}).
\[
\beta^! [p_!\E_v\restricted{]}{\M(v)^t} \ \ \equiv \ \ W_t - \left(
\sum_{i=t}^\mu e_{i,!}\left[
\tilde{\phi}_{t\cap i}^* U_{\restricted{i}{B^{[i]}\M(v)^{t\cap i}}}
\right]
\right).
\]
In the last equivalence, 
the morphisms $e_i$ and $\tilde{\phi}_{t\cap i}$ are the analogues of 
$e_i$ and $\tilde{\phi}_i$ in 
(\ref{eq-decomposition-of-K-pushforward-of-universal-sheaf-over-B-1-M-v}). 
The proof of the last equivalence is similar to that of 
(\ref{eq-decomposition-of-K-pushforward-of-universal-sheaf-over-B-1-M-v}). 
%
\EndProof

\subsection{More on the structure of the correspondence $\Delta_t^{[1]}$}
\label{sec-more-on-the-structure-of-the-correspondences}
\subsubsection{Relating the two pullbacks of the tautological subbundle}
\label{sec-relationship-between-two-pullcacks-of-U_t-to-Delta-t}
Our next goal is to prove lemma
\ref{lemma-relationship-between-two-pullcacks-of-U_t-to-Delta-t}.
Let $e_{t,k}:\Delta_t^{[1]}\cap \Delta_k^{[1]}\hookrightarrow \Delta_t^{[1]}$
be the natural embedding.
Let $\tilde{\delta}_i:\Delta^{[1]}_t\rightarrow B^{[t+1]}\M(v)^t$
be the composition of the blow-down 
$\Delta_t^{[1]}\rightarrow \Delta_t^{[t+1]}$, 
followed by the projection $\delta_i$. 
Denote by $\tilde{\delta}_{t\cap k,2}: \Delta_t^{[1]}\cap \Delta_k^{[1]}
\rightarrow B^{[k+1]}\M(v)^k$ the composition of
the inclusion 
$\Delta_t^{[1]}\cap \Delta_k^{[1]}\hookrightarrow \Delta_k^{[1]}$ and 
the projection 
$\tilde{\delta}_2:\Delta_k^{[1]}\rightarrow B^{[k+1]}\M(v)^k$.

\begin{new-lemma}
\label{lemma-relationship-between-two-pullcacks-of-U_t-to-Delta-t}
The following K-theoretic equivalence holds over $\Delta_t^{[1]}$
\[
(\tilde{\delta}_1^*U_t)^*
\otimes A_t^{-1} \ \ \equiv \ \ 
\tilde{\delta}_2^*U_t
\left(
\sum_{i=0}^{t-1}\Delta_i^{[1]}
\right)
 - 
\left\{
\StructureSheaf{\Delta_t^{[1]}}
\left(
-\sum_{j=t+1}^\mu \Delta_j^{[1]}
\right)\otimes
\sum_{k=1}^{t-1}e_{t,k_!}\tilde{\delta}_{t\cap k,2}^*U_k
\right\},
\]
where $A_t$ is the line bundle introduced in Corollary 
\ref{cor-lazarsfeld-reflection-of-universal-sheaf-over-Delta-t}.
\end{new-lemma}

We will need some preliminary work, before we prove the lemma. 
We construct in lemma \ref{lemma-petri-section-over-B-t-X-t}
a complex $\varphi$, roughly extending the complex
$B^{[t]}\rho$ in (\ref{eq-decomposable-complex}), 
from $B^{[t]}\M(v)^t$ to $B^{[t]}\X^t$ 
(see section \ref{sec-resolution-of-correspondences}, 
for the definition of $B^{[t]}\X^t$). 
The precise relation, between the two complexes, 
is stated in Claim \ref{claim-U-k-1-is-cokernel-of-varphi}. 

Let $W_{v+\vec{t}}$ be the rank $2t+2$ vector bundle over $B^{[t+1]}\M(v)^t$, 
defined in (\ref{eq-W-v-plus-t-as-an-extension}),
$W_t$ the vector bundle in proposition 
\ref{prop-decomposition-of-complex-over-interated-blow-up}, and set
\[
Q_t:=W_t\left(-\sum_{i=t+1}^\mu B^{[t+1]}\M(v)^{t\cap i}\right).
\] 
Recall, that $B^{[t+1]}\M(v)^t\rightarrow B^{[1]}\M(v+\vec{t})$
is a grassmannian bundle. The vector bundles $U_t$ and $Q_t$ 
over $B^{[t+1]}\M(v)^t$ are the tautological sub and quotient bundles
(Lemma \ref{lemma-W-v-plus-t-is-a-direct-image}). 

Set
\begin{equation}
\label{eq-E-W-v+t}
E(W_{v+\vec{t}}) \ \ := \ \ 
N_{B^{[t+1]}\M(v)^t/B^{[t+1]}\X}
\left(
\sum_{j=t+1}^\mu B^{[t+1]}\M(v)^{j\cap t}
\right) 
\end{equation}
to be the twist of the normal bundle of $B^{[t+1]}\M(v)^t$ in $B^{[t+1]}\X$.
Recall, that $B^{[t+1]}\X^t = B^{[t+1]}\M(v)^t$, but
$B^{[t]}\X^t$ contains $B^{[t]}\M(v)^t$ as a divisor 
(section \ref{sec-resolution-of-correspondences}). 
Let $\beta: B^{[t]}\X^t\rightarrow B^{[t+1]}\M(v)^t$
be the restriction of the blow-down morphism 
$B^{[t]}\X\rightarrow B^{[t+1]}\X$. 
Consider the tensor product, of 
the natural short exact sequence
\begin{equation}
\label{eq-exact-sequence-of-normal-bundles-over-B-t+1-M-v-t}
0\rightarrow
N_{B^{[t+1]}\M(v)^t/B^{[t+1]}\M(v)}
\rightarrow 
N_{B^{[t+1]}\M(v)^t/B^{[t+1]}\X}
\rightarrow 
\left(N_{B^{[t+1]}\M(v)/B^{[t+1]}\X}\restricted{\right)}{B^{[t+1]}\M(v)^t}
\rightarrow 0,
\end{equation}
with the line bundle 
$\StructureSheaf{B^{[t+1]}\M(v)^t}\left(
\sum_{j=t+1}^\mu B^{[t+1]}\M(v)^{j\cap t}
\right). 
$
We claim, that the resulting sequence is isomorphic to 
\begin{equation}
\label{eq-short-exact-seq-of-E-defining-lambda}
0\rightarrow \Hom(Q_t,U_t) \rightarrow 
E(W_{v+\vec{t}}) \RightArrowOf{\lambda} 
\StructureSheaf{B^{[t+1]}\M(v)^t} \rightarrow 0.
\end{equation}
The identification of left terms, in the two sequences, via the petri-map, 
follows from Proposition 3.18 part 9 in \cite{markman-reflections}. 
Over the special fiber of $B^{[t]}\X$ (and thus over each of its irreducible 
components), we have the rational equivalence
\begin{equation}
\label{eq-rational-equivalence-over-B-t-X}
\left(B^{[t]}\M(v) + \sum_{j=t}^\mu B^{[t]}\X^j\right)
\ \ \equiv \ \  0.
\end{equation}
This equivalence is 
proven via an argument similar to the one used in 
Lemma \ref{lemma-structure-of-correspondence-plus-2-case} 
part \ref{lemma-item-normal-bundle-of-Delta-t}. 
The triviality of the quotient line-bundle in  
(\ref{eq-short-exact-seq-of-E-defining-lambda}) follows from this rational
equivalence. 

Let $\End_0(W_{v+\vec{t}})$ be the bundle of traceless endomorphisms. 
Note, that $\Hom(Q_t,U_t)$ is the relative cotangent bundle $T^*f_t$, 
of the grassmannian bundle
$f_t:B^{[t+1]}\M(v)^t\rightarrow B^{[1]}\M(v+\vec{t})$. 
We identified in equation 
(\ref{eq-E-is-a-subbundle-of-End-H})
the non-trivial extension 
(\ref{eq-non-trivial-extension-is-deformation-of-cotangent-bundle}), 
of the trivial line-bundle, by the cotangent bundle of a Grassmannian. 
Denote by $E(W_{v+\vec{t}})'\subset \End_0(W_{v+\vec{t}})$
the relative analogue, of the subbundle (\ref{eq-E-is-a-subbundle-of-End-H}). 
We claim, that $E(W_{v+\vec{t}})$ and $E(W_{v+\vec{t}})'$ are isomorphic. 
Consequently, the vector bundle $E(W_{v+\vec{t}})$
embeds naturally as a subbundle of 
traceless endomorphisms of $W_{v+\vec{t}}$
\begin{equation}
\label{eq-embedding-of-E-W-v+t-in-traceless-endomorphisms}
E(W_{v+\vec{t}}) \ \ \hookrightarrow \ \ \End_0(W_{v+\vec{t}}).
\end{equation}
The above embedding satisfies the following properties. 
An element $\eta_x$, in the fiber of $E(W_{v+\vec{t}})$ over 
$x$, is nilpotent, if and only if $\lambda(\eta_x)=0$. 
If $\lambda(\eta_x)\neq 0$, then $\eta_x$ is semisimple
with eigenvalues $\lambda(\eta_x)$
and $-\frac{t\lambda(\eta_x)}{(t+2)}$,
of multiplicities $t$ and $t+2$ respectively.

Let us construct the isomorphism $E(W_{v+\vec{t}})\cong E(W_{v+\vec{t}})'$.
The extension (\ref{eq-exact-sequence-of-normal-bundles-over-B-t+1-M-v-t}) 
restricts, as a non-trivial extension, over each grassmannian fiber of $f_t$
(due to our choice, of the one-parameter deformation $\X$ of $\M(v)$; 
See the proof of Theorem 1.2 in \cite{markman-reflections}). 
The extension $E(W_{v+\vec{t}})'$ has the same property.
It suffices to prove, that the extension group 
$H^1(B^{[t+1]}\M(v)^t,T^*_{f_t})$ is one-dimensional. 
Using the Leray spectral sequence, we get the isomorphism
$H^1(B^{[t+1]}\M(v)^t,T^*_{f_t})\cong 
H^0(B^{[1]}\M(v+\vec{t}),R^1_{f_{t,*}}T^*_{f_t})$. 
Now, $R^1_{f_{t,*}}T^*_{f_t}$ is a trivial line-bundle over 
$B^{[1]}\M(v+\vec{t})$, because a K\"{a}hler class on $B^{[t+1]}\M(v)^t$ 
projects onto a generator, of the second singular cohomology, of
each grassmannian fiber of $f_t$. 

\begin{new-lemma}
\label{lemma-petri-section-over-B-t-X-t}
The Petri-map defines a canonical nowhere vanishing section $\eta$, 
of the vector bundle
\[
\beta^*E(W_{v+\vec{t}})(B^{[t]}\M(v)^t)
\]
over $B^{[t]}\X^t$. 
The section $\lambda(\eta)$ of 
$\StructureSheaf{B^{[t]}\X^t}(B^{[t]}\M(v)^t)$ vanishes precisely along
the divisor $B^{[t]}\M(v)^t$.
The section
\[
\varphi \ := \ \eta+\left(\frac{t\lambda(\eta)}{(t+2)}\right)I
\] 
of $\beta^*\Hom(W_{v+\vec{t}},W_{v+\vec{t}}(B^{[t]}\M(v)^t))$ 
has image in $U_t(B^{[t]}\M(v)^t)$. It satisfies the equation
\[
\varphi^2 \ \ = \ \ \frac{t+3}{t+2}\lambda(\eta)\cdot \varphi. 
\]
Consequently, the following statements hold:
\begin{enumerate}
\item
The restriction of $\varphi$ to $B^{[t]}\M(v)^t$ 
is nilpotent of square zero (and thus factors through $Q_t$). 
\item
$U_t$ is the $\lambda(\eta)$ eigen-bundle of $\eta$, 
over the complement of $B^{[t]}\M(v)^t$ in $B^{[t]}\X^t$. 
$U_t$ is thus transversal to the $-\frac{t\lambda(\eta)}{(t+2)}$
eigen-bundle of $\eta$, and $W_{v+\vec{t}}$ is the direct sum 
of the two eigen-bundles, over the complement of 
$B^{[t]}\M(v)^t$ in $B^{[t]}\X^t$.  
Over this locus, the homomorphism $\varphi$ is 
$\frac{t+3}{t+2}\lambda(\eta)$ times the projection onto the 
$\lambda(\eta)$ eigen-bundle $U_t$. Consequently, $\varphi$ 
surjects onto $U_t$, over the complement of 
$B^{[t]}\M(v)^t$ in $B^{[t]}\X^t$. 
\end{enumerate}
\end{new-lemma}

\noindent
{\bf Proof:}
By definition of $E(W_{v+\vec{t}})$, we have the natural isomorphism
$B^{[t]}\X^t\cong \PP{E}(W_{v+\vec{t}})$. 
The normal line bundle $N_{B^{[t]}\X^t/B^{[t]}\X}$
is the tautological subbundle of
$\beta^*N_{B^{[t+1]}\X^t/B^{[t+1]}\X}$. The latter is isomorphic to 
$\beta^*E(W_{v+\vec{t}})(-\sum_{j=t+1}^\mu B^{[t+1]}\M(v)^j)$.
We get an injective homomorphism
\[
N_{B^{[t]}\X^t/B^{[t]}\X} \ \ \hookrightarrow \ \ 
\beta^*E(W_{v+\vec{t}})(-\sum_{j=t+1}^\mu B^{[t+1]}\M(v)^j)
\]
corresponding to a non-vanishing section $\eta$ of
$\beta^*E(W_{v+\vec{t}})(-\sum_{j=t}^\mu B^{[t]}\X^{t\cap j})$.
The latter vector bundle is isomorphic to 
$\beta^*E(W_{v+\vec{t}})(B^{[t]}\M(v)^t)$,
via the rational-equivalence (\ref{eq-rational-equivalence-over-B-t-X}). 

The non-vanishing of $\lambda(\eta)$, away from $B^{[t]}\M(v)^t$, 
follows from the equivalence, of the 
exact sequence (\ref{eq-short-exact-seq-of-E-defining-lambda}), 
to the exact sequence of normal bundles. The rest of the statement
follows from the properties of the embedding
(\ref{eq-embedding-of-E-W-v+t-in-traceless-endomorphisms}).
\EndProof

\bigskip
Set $B^{[t]}U_t := (\beta^*U_t)\otimes N^*_{B^{[t]}\X^t/B^{[t]}\X}$. 
It is a vector bundle over $B^{[t]}\X^t$. 
The equivalence (\ref{eq-rational-equivalence-over-B-t-X}) yields
\[
B^{[t]}U_t \ \ = \ \ 
(\beta^*U_t)\left(
B^{[t]}\M(v)^t+\sum_{j=t+1}^\mu B^{[t]}\X^{t\cap j}
\right).
\] 
We construct, by descending recursion on $k$, $1\leq k\leq t$, 
a sequence of vector bundles $B^{[k]}U_t$ over $B^{[k]}\X^t$ and 
homomorphisms
\[
B^{[k]}\varphi \ : \ \beta^*W_{v+\vec{t}}
\left(\sum_{j=t+1}^\mu B^{[k]}\X^j
\right)
\ \ \rightarrow \ \ B^{[k]}U_t,
\]
whose fiberwise rank is $\geq t-k+1$. 
The rank is equal to $t-k+1$ along $B^{[k]}\X^{t\cap (k-1)}$. 
Note, that $B^{[k]}\X^{t\cap (k-1)}=B^{[k]}\M(v)^{t\cap (k-1)}$,
sinse $B^{[k]}\X^{k-1}=B^{[k]}\M(v)^{k-1}$.
We set $B^{[t]}\varphi:=\varphi$,
the homomorphism constructed in lemma
\ref{lemma-petri-section-over-B-t-X-t}. 
Note, that the homomorphism $B^{[k]}\varphi$ factors through $\beta^*W_t$.

Assume, that the homomorphism $B^{[k]}\varphi$ has been constructed.
Denote by $K_{k-1}$ and $C_{k-1}$ 
the kernel and image of
the restriction
$\restricted{B^{[k]}\varphi}{B^{[k]}\X^{t\cap (k-1)}}$.
Let $B^{[k-1]}U_t$ be the subsheaf of $\beta^*B^{[k]}U_t$
over $B^{[k-1]}\X^t$, with values in $C_{k-1}$ along 
$B^{[k-1]}\X^{t\cap (k-1)}$. 
The pullback $\beta^*B^{[k]}\varphi$, of $B^{[k]}\varphi$ to
$B^{[k-1]}\X^t$, has image in the subsheaf $B^{[k-1]}U_t$ of 
$\beta^*B^{[k]}U_t$. We define $B^{[k-1]}\varphi$ to be the homomorphism 
induced by $\beta^*B^{[k]}\varphi$.

\begin{claim}
\label{claim-U-k-1-is-cokernel-of-varphi}
For $1\leq k\leq t$, the cokernel of 
$\restricted{B^{[k]}\varphi}{B^{[k]}\X^{t\cap (k-1)}}$ 
is precisely the restriction of $U_{k-1}$. 
\end{claim}

\noindent
{\bf Proof:}
This follows from a general property of kernels and cokernels of Petri-maps
(Lemma 4.1 in \cite{markman-reflections}).
In our context of Proposition
\ref{prop-decomposition-of-complex-over-interated-blow-up}, 
this property can be stated as follows.
Over $B^{[t]}\M(v)^t$, the Petri-map is a canonical
homomorphism
\[
\dot{\rho} \ : \ \beta^*W_t \ \ \rightarrow \ \ 
\beta^*U_t\otimes N^*_{B^{[t]}\M(v)^t/B^{[t]}\M(v)}.
\]
Moreover, the kernel and co-kernel of the restriction 
$(B^{[t]}\rho\restricted{)}{B^{[t]}\M(v)^t}$, 
of the homomorphism given in (\ref{eq-decomposable-complex}), 
are equal to the kernel and co-kernel of the Petri map $\dot{\rho}$,
respectively. Furthermore, the Petri map $\ddot{\rho}$ of the Petri-map
$\dot{\rho}$ is the restriction of the Petri-map of $B^{[t]}\rho$ to
$B^{[t]}\M(v)^{t\cap (t-1)}$.

By definition, $U_{t-1}$ is the co-kernel of the homomorphism 
$B^{[t]}\rho$ over $B^{[t]}\M(v)^{t-1}$. 
Thus, the restriction of $U_{t-1}$ to $B^{[t]}\M(v)^{t\cap (t-1)}$ 
is also the cokernel of the Petri map $\dot{\rho}$. 
Our definition of $B^{[t]}\varphi$ is such, that it factors through
$\dot{\rho}$ over $B^{[t]}\M(v)^t$. This proves the claim, in the case $k=t$. 

The general case $k\leq t$ is proven by induction, applying 
repeatedly the above mentioned properties of Petri maps. The point is that, 
over $B^{[k]}\X^{t\cap (k-1)}=B^{[k]}\M(v)^{t\cap (k-1)}$, the Petri-map of 
(the Petri-map) $B^{[k]}\varphi$, is equal to the Petri-map of
$B^{[k]}\rho$. 
\EndProof

\medskip
We have the exact sequence over $B^{[k-1]}\X^t$
\[
0\rightarrow 
B^{[k-1]}U_t \rightarrow 
\beta^*B^{[k]}U_t \rightarrow 
e_{k-1_*}U_{k-1} \rightarrow 0.
\]
Over $B^{[1]}\X^t$  we get, for $1\leq k \leq t$, 
\[
0\rightarrow 
\beta^*B^{[k-1]}U_t \rightarrow 
\beta^*B^{[k]}U_t \rightarrow 
e_{k-1_*}\tilde{\phi}_{t,k-1}^*U_{k-1} \rightarrow 0,
\]
where $e_k$ is the closed emersion
$B^{[j]}\M(v)^{t\cap k} \hookrightarrow B^{[j]}\X^t$, for $1\leq j \leq k$, 
and $\tilde{\phi}_{t,k}:B^{[1]}\M(v)^{t\cap k}\rightarrow B^{[k+1]}\M(v)^k$ 
is the natural morphism. Note, that the pullback of $B^{[t]}U_t$ to
$B^{[1]}\X^t$ is
\[
\beta^*(B^{[t]}U_t) \ \ \cong \ \ 
(\beta^*U_t)\left(
B^{[1]}\M(v)^t+
\sum_{i=1}^{t-1}B^{[1]}\X^i+\sum_{j=t+1}^\mu B^{[1]}\X^j
\right) \ \ \cong \ \ 
(\beta^*U_t)\left(-B^{[1]}\X^t\right),
\]
because the center of the blow-up operation of 
$B^{[k]}\X^t\rightarrow B^{[k+1]}\X^t$, for $k<t$, is contained in the
divisor $B^{[k+1]}\M(v)^t$. 
We conclude the K-theoretic equivalence over 
$B^{[1]}\X^t$
\begin{equation}
\label{eq-K-theoretic-decomposition-of-iterated-blow-up-of-U-t}
B^{[1]}U_t \ \ \equiv \ \ 
(\beta^*U_t)
\left(
-B^{[1]}\X^t
\right)
-\sum_{k=1}^{t-1}e_{k_!}\tilde{\phi}_{t,k}^*U_k. 
\end{equation}
Proposition \ref{prop-decomposition-of-complex-over-interated-blow-up}
part \ref{property-constant-rank}
and Claim \ref{claim-U-k-1-is-cokernel-of-varphi} imply, that
$B^{[1]}\varphi$ is surjective. Hence, the kernel of $B^{[1]}\varphi$
is a rank $t+2$ subbundle of 
$\beta^*W_{v+\vec{t}}(\sum_{j=t+1}^\mu B^{[1]}\X^j)$. 
The kernel determines a morphism 
\[
\tilde{\delta}_1\times \tilde{\delta}_2 \ : \ B^{[1]}\X^t  \ \rightarrow 
 \ G(t,W^*_{v+\vec{t}}).
\]
Let us justify the notation. 
$G(t,W^*_{v+\vec{t}})$ is isomorphic to $\Delta_t^{[t+1]}$, 
by Lemma \ref{lemma-W-v-plus-t-is-a-direct-image}. 
Composing the last two morphisms, we get the morphism
$B^{[1]}\X^t\rightarrow \Delta_t^{[t+1]}$, which factors through
the isomorphism $B^{[1]}\X^t\cong \Delta_t^{[1]}$ of 
Lemma \ref{lemma-two-descriptions-of-the-components-of-the-correspondence}.
The morphism $\tilde{\delta}_i$, $i=1,2$, is identified with 
the composition, of the projection 
$\delta_i: \Delta_t^{[t+1]}\rightarrow B^{[t+1]}\M(v)^t$, 
with the blow-down morphism 
$\beta:\Delta_t^{[1]}\rightarrow \Delta_t^{[t+1]}$. The morphism 
$\tilde{\delta}_1$ above is identified with the restriction to $B^{[1]}\X^t$
of the blow-down morphism $B^{[1]}\X\rightarrow B^{[t+1]}\X$. 

Let $\tau_v^{t}$ and $q_v^{t}$ be the relative tautological 
sub and quotient bundles over $G(t,W^*_{v+\vec{t}})$. 
We conclude the isomorphism, via $B^{[1]}\varphi^*$,
\begin{equation}
\label{eq-tau-is-iterated-blow-up-of-U-t}
\tau_v^{t} \ \ \cong \ \ 
B^{[1]}U_t^*(\sum_{j=t+1}^\mu \Delta_j^{[1]}). 
\end{equation}

We are now able to prove Lemma
\ref{lemma-relationship-between-two-pullcacks-of-U_t-to-Delta-t}. 

\noindent
{\bf Proof of lemma
\ref{lemma-relationship-between-two-pullcacks-of-U_t-to-Delta-t}:}
We have the isomorphism 
$\delta_2^*W_{v+\vec{t}}\otimes A_t \cong \delta_1^* W^*_{v+\vec{t}}$
(see equation (\ref{eq-A-t-is-symmetric}) in Corollary 
\ref{cor-lazarsfeld-reflection-of-universal-sheaf-over-Delta-t}). 
Consequently, $\delta_2$ pulls back 
the short exact sequence (\ref{eq-W-v-plus-t-as-an-extension})
to a twist of the tautological short exact sequence of
$\tau_v^{t}$ by the line bundle $A_t^{-1}$. In particular, we have 
\begin{equation} 
\label{eq-second-pullback-of-U-t-twisted-by-A-t-is-tau}
\delta_2^*U_t  \otimes A_t \ \ \cong \ \ \tau_v^{t} 
\end{equation}
The lemma follows from equations
(\ref{eq-K-theoretic-decomposition-of-iterated-blow-up-of-U-t}) 
and (\ref{eq-tau-is-iterated-blow-up-of-U-t}) and part
\ref{lemma-item-normal-bundle-of-Delta-t} of
Lemma \ref{lemma-structure-of-correspondence-plus-2-case}.
\EndProof

\subsubsection{The incidence divisor $I_t$ in $\Delta_t^{[t+1]}$}
\label{sec-incidence-divisor}
Part 
\ref{lemma-item-component-is-a-fiber-product-of-bundle-with-its-pullback} 
of lemma \ref{lemma-structure-of-correspondence-plus-2-case}
provides two descriptions of the component
$\Delta_t^{[1]}$ of the correspondence $\Z^{[1]}$. Both descriptions realize 
$\Delta_t^{[1]}$ as the end result of a sequence of 
blow-ups, but the starting spaces 
$\Delta_t^{[t+1]}$ and $B^{[t]}\X^t$ 
are non-isomorphic. Our next goal is to relate the exceptional 
divisors, in $\Delta_t^{[1]}$, of the two blow-up sequences. 
We will see, that the exceptional divisors are the same, but the
order of the two blow-down sequence is 
\mbox{\em reversed} in the following sense. 
We will construct below in (\ref{eq-alpha_t}) 
a homomorphism $\alpha_t$, between two vector bundles over 
$\Delta_t^{[t+1]}$, whose determinantal subvarieties are the centers of the
blow-up operations. 
Set $B^{[k]}\X^{t\cap 0}:=B^{[k]}\M(v)^t$ and, for $1\leq k < t$, set
$B^{[k]}\X^{t\cap k}:= B^{[k]}\X^t\cap B^{[k]}\X^k$. 
Then 
\[
B^{[1]}\X^{t\cap k} \ \ = \ \ \Delta^{[1]}_t\cap \Delta^{[1]}_k,
\]
for $0\leq k \leq t-1$, 
by lemma \ref{lemma-structure-of-correspondence-plus-2-case} part 
\ref{lemma-item-pull-back-of-BN-divisors-to-correspondence-in-plus-2-case}.
However, the generic rank of $\beta^*\alpha_t$ along 
$\Delta^{[1]}_t\cap \Delta^{[1]}_k$ is  $k$
(lemma \ref{lemma-two-degeneracy-loci-on-Delta-t-are-identical}), while
the generic rank of $\beta^*\varphi$ along
$B^{[1]}\X^{t\cap k}$ is $(t-k)$,
by lemma \ref{lemma-petri-section-over-B-t-X-t}. 

The bundle $\delta_1^*W_{v+\vec{t}}$ over $\Delta_t^{[t+1]}$
has the tautological rank $t$ subbundle $\delta_1^*U_t$ 
and a tautological rank $t$ quotient bundle, which we denote by $q$. 
We have the homomorphism $\alpha_t$ over $\Delta_t^{[t+1]}$
\begin{equation}
\label{eq-alpha_t}
\delta_1^*U_t \ \hookrightarrow \ 
\delta_1^*W_{v+\vec{t}} \ \RightArrowOf{j} \ q
\end{equation}
relating $\delta_1^*U_t$ and $q$
(compare with 
(\ref{eq-homomorphism-between-tautological-subundle-of-dual-grassmannians})). 
The incidence divisor $I_t\subset \Delta_t^{[t+1]}$ 
is the zero divisor of the 
determinant of the homomorphism $\alpha_t$. 
The homomorphism $\alpha_t$ vanishes along the flag subvariety
$Flag(t,t+2,W_{v+\vec{t}}) \subset \Delta_t^{[t+1]}$. 

\begin{new-lemma}
\label{lemma-two-degeneracy-loci-on-Delta-t-are-identical}
\begin{enumerate}
\item
\label{lemma-item-generic-rank-of-alpha-t}
The pullback of $\alpha_t$ to $\Delta_t^{[1]}$ has generic rank 
$k$ along $\Delta_t^{[1]}\cap\Delta_k^{[1]}$, for $0\leq k <t$. 
\item
\label{lemma-item-pullback-of-incidence-divisor-has-multiplicity-t-i-along}
Under the isomorphism $B^{[1]}\X^t\cong \Delta_t^{[1]}$, 
the pullback of $\StructureSheaf{\Delta_t^{[t+1]}}(-I_t)$ 
to $\Delta_t^{[1]}$ is isomorphic to
\[
\StructureSheaf{B^{[1]}\X^t}
\left(-\sum_{i=0}^{t-1}(t-i)B^{[1]}\X^{t\cap i}
\right).
\]
\item
\label
{lemma-item-pushforward-of-ideal-sheaf-is-ideal-sheaf-of-incidence-div+2}
The sheaf theoretic pushforward of the line-bundle
${\displaystyle 
\StructureSheaf{\Delta^{[1]}_t}\left(-\sum_{i=0}^{t-1}
\Delta_i^{[1]}\cap\Delta^{[1]}_t\right)
}$, 
from $\Delta^{[1]}_t$ to $\Delta^{[t+1]}_t$, 
is isomorphic to the ideal sheaf $\StructureSheaf{\Delta^{[t+1]}_t}(-I_t)$ 
of the incidence divisor $I_t$ in the
fiber product of the two dual grassmannian bundles. The higher direct images
vanish. 
\end{enumerate}
\end{new-lemma}

\medskip
{\bf Proof:}
\ref{lemma-item-generic-rank-of-alpha-t}) 
Denote by $\beta^*\alpha_t$ and 
$\beta^*j$ the pullbacks to $\Delta_t^{[1]}$ of $\alpha_t$ and $j$, 
given in (\ref{eq-alpha_t}) over $\Delta_t^{[t+1]}$.
Let $\beta^*\varphi$ and $\beta^*\lambda(\eta)$ be the pullbacks 
to $\Delta_t^{[1]}$ of $\varphi$ and $\lambda(\eta)$, 
given in Lemma \ref{lemma-petri-section-over-B-t-X-t} over $B^{[t]}\X^t$.
(We denote by $\beta$ both of the blow-up morphisms
$\Delta_t^{[1]}\rightarrow\Delta_t^{[t+1]}$ 
and $\Delta_t^{[1]}\cong B^{[1]}\X^t\rightarrow B^{[t]}\X^t$).
Let us prove the equality
\begin{equation}
\label{eq-composition-of-alpha-with-varphi-is-j}
\beta^*\alpha_t\circ \beta^*\varphi \ \ = \ \ 
\frac{t+3}{t+2}\beta^*\lambda(\eta)\cdot \beta^*j, 
\end{equation}
of homomorphisms from
$\tilde{\delta}_1^*W_{v+\vec{t}}$ to 
$\beta^*q\otimes\beta^*\StructureSheaf{B^{[t]}\X^t}(B^{[t]}\M(v)^t)$.
Note, that $\Delta_t^{[t+1]}\setminus I_t$ is naturally isomorphic to
the complement 
$B^{[t]}\X^t \setminus B^{[t]}\M(v)^t$. 
Over this locus, the identity clearly follows from  
Lemma \ref{lemma-petri-section-over-B-t-X-t}.

The pullback
$\beta^*\StructureSheaf{B^{[t]}\X^t}(B^{[t]}\M(v)^t)$ 
is the line bundle 
$\StructureSheaf{B^{[1]}\X^t}(\sum_{i=0}^{t-1}B^{[1]}\X^{t\cap i})$
(Lemma \ref{lemma-criterion-for-normal-crossing}).
Its section $\beta^*\lambda(\eta)$ has zero-divisor  
$\sum_{i=0}^{t-1}B^{[1]}\X^{t\cap i}$. 
Equation 
(\ref{eq-composition-of-alpha-with-varphi-is-j}) implies the inequality
$\rank(\beta^*\alpha_t)+\rank(\beta^*\varphi)\leq t$,
at every point in $\sum_{i=0}^{t-1}B^{[1]}\X^{t\cap i}$. 
Moreover, equality holds, at every smooth point in 
the normal crossing divisor $\sum_{i=0}^{t-1}B^{[1]}\X^{t\cap i}$.

Part
\ref{lemma-item-pullback-of-incidence-divisor-has-multiplicity-t-i-along}
clearly follows from part \ref{lemma-item-generic-rank-of-alpha-t}.

Part 
\ref{lemma-item-pushforward-of-ideal-sheaf-is-ideal-sheaf-of-incidence-div+2}
follows from part 
\ref{lemma-item-pullback-of-incidence-divisor-has-multiplicity-t-i-along}
via an argument analogous to that used in the proof of Claim
\ref{claim-push-forward-of-ideal-sheaf-is-ideal-sheaf}.
\EndProof

\subsection{The line bundle $A_t$}
\label{sec-the-line-bundle-A-t}

Our next task is to calculate the line bundle $A_t$ in corollary 
\ref{cor-lazarsfeld-reflection-of-universal-sheaf-over-Delta-t}. 
The calculation is completed in Lemma \ref{lemma-calculation-of-A-t}. 
The birational involution $\iota:\M(v)\rightarrow\M(v)$ 
lifts to a regular involution $g$ of $B^{[1]}\M(v)$
(Proposition 3.26 in \cite{markman-reflections}). 
We have the following commutative diagram (Remark 3.25 in 
\cite{markman-reflections}) 
\begin{equation}
\label{eq-commutative-diagram-of-two-involutions-g-iota}
\begin{array}{ccccccc}
B^{[1]}\M(v) & \supset & B^{[1]}\M(v)^t & \LongRightArrowOf{\beta_1} &
B^{[t+1]}\M(v)^t & \rightarrow & B^{[1]}M(v+\vec{t})
\\
g\uparrow & & \restricted{g}{} \ \uparrow & & \tilde{\iota}\uparrow & &
\iota \uparrow
\\
B^{[1]}\M(v) & \supset & B^{[1]}\M(v)^t & \LongRightArrowOf{\beta_2} &
\iota^*B^{[t+1]}\M(v)^t & \rightarrow & B^{[1]}M(v+\vec{t}),
\end{array}
\end{equation}
where $\iota$ is the extension of (\ref{eq-sigma-involution-on-lower-strata}) 
to a regular involution 
(part \ref{lemma-item-iota-extends-to-a-regular-involution}
of Lemma \ref{lemma-structure-of-correspondence-plus-2-case}). 
Consider Lemma \ref{lemma-lazarsfeld-reflection-of-a-family-of-sheaves}
with the triple $T,W,\F$ being $B^{[1]}\M(v),W_0,\beta^*\E_v$, 
where $\E_v$ is the universal family over $\M(v)$, normalized 
to satisfy $\det(\E_v)=\pi_S^*\LB$. 
We get the following Lemma as a special case.

\begin{new-lemma}
\label{lemma-lazarsfeld-reflection-of-universal-sheaf}
Choose the universal family $\E_v$ over $\M(v)\times S$ to be the twist of
the ideal sheaf of the universal subscheme of $S^{[g-1]}$ by $\LB$. 
The pullback $g^*(\beta^*\E_v)$ is related to $\E_v$ 
by the following exact sequence over $B^{[1]}\M(v)\times S$
\begin{equation}
\label{eq-relation-between-E-v-and-its-pull-back-of-by-f}
0\rightarrow 
\SheafHom(\beta^*\E_v,\StructureSheaf{}) \LongRightArrowOf{ev^*} 
\pi_1^*W_0^* \rightarrow 
g^*(\beta^*\E_v)\otimes \det(\pi_1^*W_0)^* \rightarrow 
\SheafExt^1_{B^{[1]}\M(v)\times S}(\beta^*\E_v,\StructureSheaf{}) 
\rightarrow 0,
\end{equation}
where $\pi_1$ is the projection from $B^{[1]}\M(v)\times S$ onto the first
factor.
Moreover, the extension sheaves
$\SheafExt^i_{B^{[1]}\M(v)\times S}(\beta^*\E_v,\StructureSheaf{})$
vanish, for $i\geq 2$. Consequently, we have the following equivalence in
the K-group of $B^{[1]}\M(v)\times S$
\begin{equation}
\label{eq-identification-of-pull-back-of-E-v-by-f}
g^*(\beta^*\E_v)\otimes \det(\pi_1^*W_0)^* - 
\pi_1^!\sum_{t=1}^\mu\left[e_{t,!}\tilde{\phi}_t^*U_t\right]^\vee
\ \ \equiv \ \ 
\beta^![-\tau(\E_v)]^\vee, 
\end{equation}
Above, $\tau(\E_v)$ is the K-theoretic representative, given in
(\ref{eq-tau-of-universal-sheaf}),
of the reflection of $\E_v$ with respect to $\StructureSheaf{S}$.
The morphism $\tilde{\phi}_t:B^{[1]}\M(v)^t\rightarrow B^{[t+1]}\M(v)^t$ 
is the natural one (see 
(\ref{eq-decomposition-of-K-pushforward-of-universal-sheaf-over-B-1-M-v})). 
\end{new-lemma}

{\em Note: 
the line bundle $\det(\pi_1^*W_0)$ is {\em not} a pullback 
of a line bundle on $\M(v)$. 
}

\noindent
{\bf Proof}
(of Lemma \ref{lemma-lazarsfeld-reflection-of-universal-sheaf}).
The lemma is the special case $t=0$ in
corollary
\ref{cor-lazarsfeld-reflection-of-universal-sheaf-over-Delta-t},
where $W_{v+\vec{0}}=W_0(-\sum_{i=1}^\mu\Delta_0^{[1]}\cap \Delta_i^{[1]})$
and $\F_0=\E_v(-\sum_{i=1}^\mu\Delta_0^{[1]}\cap \Delta_i^{[1]})$. 

Following is a second proof. 
The construction of the isomorphism $g$ in Theorem 3.21 of
\cite{markman-reflections} implies, that 
the pull back $g^*(\beta^*\E_v)$ is isomorphic to  
$\sigma(\beta^*\E_v,W_0)\otimes A'$, 
for some line bundle $A'$ over $B^{[1]}\M(v)$. 
Note the equality $(g\times id_S)^*\beta^*\det(\E_v)=\beta^*\det(\E_v)$,
which follows from our choice of the universal family.
The determination $A'=\pi_1^*\det(W_0)$ follows from the fact,
that the determinant line bundle of 
$\SheafExt^1_{B^{[1]}\M(v)\times S}(\beta^*\E_v,\StructureSheaf{})$
is trivial, as the sheaf has support of codimension $\geq 2$. 
Equality (\ref{eq-identification-of-pull-back-of-E-v-by-f}) 
follows from 
(\ref{eq-relation-between-E-v-and-its-pull-back-of-by-f}) and equality 
(\ref{eq-decomposition-of-K-pushforward-of-universal-sheaf-over-B-1-M-v}). 
\EndProof

\begin{new-lemma}
\label{lemma-the-image-of-eta-v-under-Z}
Choose the universal family as in Lemma
\ref{lemma-lazarsfeld-reflection-of-universal-sheaf}.
Then the following equality holds:
\begin{equation}
\label{eq-Z-sigma-takes-eta-v-to}
\Z_*\circ D_\M(c_1(\eta_v)) \ = \ c_1(\eta_v)+c_1(p_*\E_v),
\end{equation}
where $\eta_v$ is the class in 
(\ref{eq-invariant-normalized-class-of-chern-character-of-universal-sheaf}). 
Moreover, $c_1(p_*\E_v)$ is $\Z_*$ invariant. 
\end{new-lemma}

\noindent
{\bf Proof:}
The stratum $\M(v)^1$ has codimension $3$. 
Hence, the action of $\Z_*$ on the second cohomology of $\M(v)$ coincides with
the action of the graph of the birational involution. 
Let $E$ be a sheaf on $S$ with Mukai vector $v$. We get
\[
(v,v)\cdot \Z_*(c_1(\eta_v)) \ = \
\Z_*(\theta_v(v)) \ = \ 
c_1p_!\beta_!\left[(\pi_S^*E^\vee)\otimes g^*(\beta^*\E_v)\right].
\]
Lemma \ref{lemma-lazarsfeld-reflection-of-universal-sheaf} 
(ignoring terms supported in codimension $\geq 2$) yields 
\[
c_1p_!\beta_!\left[(\pi_S^*E^\vee)\otimes \beta^!(-\tau(\E_v))^\vee
\otimes \det(\pi_1^*W_0)\right].
\]
The equality $\pi_1=p\circ \beta$ and the projection formula yield 
\[
c_1p_!\left[(\pi_S^*E^\vee)\otimes (-\tau(\E_v))^\vee
\right] - (v,-\tau\circ D(v))c_1(\beta_! W_0).
\]
Note that $-\tau\circ D(v)=\sigma(v)=v$. 
Equation 
(\ref{eq-decomposition-of-K-pushforward-of-universal-sheaf-over-B-1-M-v}) 
implies the equality $c_1(p_*\E_v)=c_1(\beta_! W_0)$.
Equation (\ref{eq-chern-character-of-tau-E-v}), 
the fact that $\tau$ is an isometry of the Mukai lattice,  
and Grothendieck-Riemann-Roch, yield
\[
c_1p_!\left\{\left[(-\pi_S^*\tau(E^\vee))\otimes \E_v
\right]^\vee\right\}   - (v,v)c_1(p_*\E_v).
\]
Finally, the above is equal to
\[
- \theta_v(\sigma(v))  - (v,v)c_1(p_*\E_v),
\]
by Grothendieck-Serre Duality. 
Equality (\ref{eq-Z-sigma-takes-eta-v-to}) follows. 

The invariance of $c_1(p_*\E_v)$ follows from equation
(\ref{eq-Z-sigma-takes-eta-v-to}) and the fact that 
$\Z_*$ and $D_\M$ are two commuting involutions. 
\EndProof

\medskip
Recall, that $\Delta_0^{[1]}$ is the graph in $B^{[1]}\M(v)\times B^{[1]}\M(v)$
of the involution $g$. Identify $\Delta_0^{[1]}$ with $B^{[1]}\M(v)$
via the first projection. Under this identification, we have
$\Delta_0^{[1]}\cap \Delta_t^{[1]}=B^{[1]}\M(v)^t$ (Lemma
\ref{lemma-structure-of-correspondence-plus-2-case} part 
\ref{lemma-item-pull-back-of-BN-divisors-to-correspondence-in-plus-2-case}).
Let $\widetilde{W}_0$ be the subbundle of 
$\tilde{\phi}_t^*W_{v+\vec{t}}$ over $B^{[1]}\M(v)^t$, which projects via
(\ref{eq-W-v-plus-t-as-an-extension}) onto 
the subbundle $W_0(-\sum_{i=t+1}^\mu B^{[1]}\M(v)^{t\cap i})$ of
$(\tilde{\phi}_t^*W_t)(-\sum_{i=t+1}^\mu B^{[1]}\M(v)^{t\cap i})$. 
Consider the following flag of vector bundles over $B^{[1]}\M(v)^t$
\begin{equation}
\label{eq-flag-U-t-subset-W-0-subset-W-v-plus-t}
\tilde{\phi}_t^*U_t \ \subset \ \widetilde{W}_0 \ \subset \ 
\tilde{\phi}_t^*W_{v+\vec{t}}.
\end{equation}
Dualizing, we get the corresponding flag (of the annihilating subbundles)
\begin{equation}
\label{eq-dual-flag-of-annihilators}
\left[(\tilde{\phi}_t^*W_t)/W_0\right]^*
\left(\sum_{i=t+1}^\mu B^{[1]}\M(v)^{t\cap i}\right)  
\ \subset \ 
\tilde{\phi}_t^*W_t^*\left(\sum_{i=t+1}^\mu B^{[1]}\M(v)^{t\cap i}
\right)  
\ \subset \ 
\tilde{\phi}_t^*W_{v+\vec{t}}^*.
\end{equation}

\begin{new-lemma}
\label{lemma-pull-back-of-filtration-is-dual-filtration-twisted-by-A-t}
The pullback, of the flag (\ref{eq-flag-U-t-subset-W-0-subset-W-v-plus-t})
by the involution $g$, is isomorphic to the tensor product of the flag
(\ref{eq-dual-flag-of-annihilators}) 
with the line bundle $A_t^*$ in Corollary
\ref{cor-lazarsfeld-reflection-of-universal-sheaf-over-Delta-t}.
\end{new-lemma}

\noindent
{\bf Proof:}
We have a commutative diagram
\begin{equation}
\label{diagram-relating-the-two-projections-delta-i-by-the-involution-g}
\begin{array}{ccccccccc}
\Delta_0^{[1]} & \supset &
\Delta_0^{[1]}\cap\Delta_t^{[1]} & \subset & \Delta_t^{[1]} & 
\LongRightArrowOf{\beta} & \Delta_t^{[t+1]} & \LongRightArrowOf{\delta_1} & 
B^{[t+1]}\M(v)^t
\\
g \downarrow \hspace{1ex} & & g_\mid \ \downarrow \ \hspace{2ex} & & & & & &
\hspace{1ex} \downarrow =
\\
\Delta_0^{[1]} & \supset &
\Delta_0^{[1]}\cap\Delta_t^{[1]} & \subset & \Delta_t^{[1]} & 
\LongRightArrowOf{\beta} & \Delta_t^{[t+1]} & \LongRightArrowOf{\delta_2} & 
B^{[t+1]}\M(v)^t
\end{array}
\end{equation}
The isomorphism 
$g^*\tilde{\phi}_t^*W_{v+\vec{t}}\cong 
\tilde{\phi}_t^*W_{v+\vec{t}}^*\otimes A_t^*$
follows from the equality 
$\delta_2\circ\beta\circ g = \delta_1\circ \beta$ on
$\Delta_0^{[1]}\cap\Delta_t^{[1]}$ 
and Corollary
\ref{cor-lazarsfeld-reflection-of-universal-sheaf-over-Delta-t}. 
The isomorphism between 
the two flags of subbundles follows, via Lemma
\ref{lemma-structure-of-correspondence-plus-2-case} part 
\ref{lemma-item-components-of-correspondence-intersect-as-flag-variety},
from the fact, that the involution $g$ of $B^{[1]}\M(v)^t$ 
interchanges the two rulings of $\Delta_0^{[1]}\cap\Delta_t^{[1]}$, 
when the latter is identified with $B^{[1]}\M(v)^t$ via Lemma 
\ref{lemma-structure-of-correspondence-plus-2-case} part 
\ref{lemma-item-pull-back-of-BN-divisors-to-correspondence-in-plus-2-case}. 
\EndProof

We are now able to calculate the line bundle $A_t$
introduced in Corollary 
\ref{cor-lazarsfeld-reflection-of-universal-sheaf-over-Delta-t}.
Let $\tilde{\delta}_i:\Delta_t^{[1]}\rightarrow B^{[t+1]}\M(v)^t$ be the 
composition of 
$\delta_i:\Delta_t^{[t+1]}\rightarrow B^{[t+1]}\M(v)^t$
with the blow-down morphism. When $t=0$, 
both $\delta_1$ and $\delta_2$ are isomorphisms and 
$\delta_2=\delta_1\circ g$. 

\begin{new-lemma}
\label{lemma-calculation-of-A-t}
Assume that $\E_v$ is normalized as in 
Lemma \ref{lemma-lazarsfeld-reflection-of-universal-sheaf}.
Over $\Delta_t^{[1]}$ we have the isomorphism
\begin{eqnarray}
\label{eq-A-t-in-terms-of-E-v}
\beta^*A_t &  \cong &
(\beta^*\tilde{\delta}_1^*\det(p_!\E_v))^*\left(
\sum_{i=1}^t -i[\Delta_t^{[1]} \cap \Delta_i^{[1]}]
+
\sum_{i=t+1}^\mu (2-i)[\Delta_t^{[1]} \cap \Delta_i^{[1]}]
\right)
\\
\nonumber
&  \cong &
(\beta^*\tilde{\delta}_1^*\det(p_!\E_v))^*\left(
\sum_{i=0}^{t-1} (t-i)[\Delta_t^{[1]} \cap \Delta_i^{[1]}]
+
\sum_{i=t+1}^\mu (t+2-i)[\Delta_t^{[1]} \cap \Delta_i^{[1]}]
\right).
\end{eqnarray}
Over $\Delta_t^{[t+1]}$, we have
\[
A_t  \ \ \cong \ \ 
(\delta_1^*\det(p_!\E_v))^*(I_t)\left(
\sum_{i=t+1}^\mu (t+2-i)\Delta_{t,i}^{[t+1]}
\right),
\]
where $\Delta_{t,i}^{[t+1]}$ is defined in Corollary 
\ref{cor-lazarsfeld-reflection-of-universal-sheaf-over-Delta-t}
and $I_t$ is the incidence divisor on $\Delta_t^{[t+1]}$,
introduced in lemma \ref{lemma-two-degeneracy-loci-on-Delta-t-are-identical} 
part 
\ref{lemma-item-pushforward-of-ideal-sheaf-is-ideal-sheaf-of-incidence-div+2}. 
The isomorphisms hold for the pullback by $\delta_2$ as well. 
\end{new-lemma}

\noindent
{\bf Proof:}
Observe, that 
$(\delta_1^*W_0\restricted{)}{\Delta_0^{[1]}\cap\Delta_t^{[1]}}
\left(
-\sum_{i=t+1}^\mu \Delta_0^{[1]}\cap\Delta_i^{[1]}\cap\Delta_t^{[1]}
\right)$
is the middle graded summand obtained from
the filtration (\ref{eq-flag-U-t-subset-W-0-subset-W-v-plus-t}).
Similarly, the middle graded summand of the dual filtration
(\ref{eq-dual-flag-of-annihilators}) is
$(\delta_1^*W_0^*\restricted{)}{\Delta_0^{[1]}\cap\Delta_t^{[1]}}
\left(
\sum_{i=t+1}^\mu \Delta_0^{[1]}\cap\Delta_i^{[1]}\cap\Delta_t^{[1]}
\right).$
Lemma \ref{lemma-pull-back-of-filtration-is-dual-filtration-twisted-by-A-t}
implies the isomorphism
\[
(g^*\delta_1^*W_0\restricted{)}{\Delta_0^{[1]}\cap\Delta_t^{[1]}}
\left(
-\sum_{i=t+1}^\mu \Delta_0^{[1]}\cap\Delta_i^{[1]}\cap\Delta_t^{[1]}
\right)
\ \ \cong \ \ 
(\delta_1^*W_0^*\restricted{)}{\Delta_0^{[1]}\cap\Delta_t^{[1]}}
\left(
\sum_{i=t+1}^\mu \Delta_0^{[1]}\cap\Delta_i^{[1]}\cap\Delta_t^{[1]}
\right)
\otimes A_t^{-1}.
\]
Taking determinants, we get
\[
(A_t^{2}\restricted{)}{\Delta_0^{[1]}\cap\Delta_t^{[1]}} \ \ \cong \ \  
(\delta_1^*\det(W_0)^*\otimes 
\delta_2^*\det(W_0)^*\restricted{)}{\Delta_0^{[1]}\cap\Delta_t^{[1]}}
\left(
4\cdot \sum_{i=t+1}^\mu \Delta_0^{[1]}\cap\Delta_i^{[1]}\cap\Delta_t^{[1]}
\right).
\]
Equation 
(\ref{eq-decomposition-of-K-pushforward-of-universal-sheaf-over-B-1-M-v}) 
implies the following isomorphism of line bundles over $B^{[1]}\M(v)$ 
\begin{equation}
\label
{eq-comparison-between-det-W-0-and-det-of-pushforward-of-universal-sheaf}
\det(W_0)^* \ \ \cong \ \ (\beta^*\det(p_!\E_v))^*\left(
-\sum_{t=1}^\mu t\cdot B^{[1]}\M(v)^t
\right). 
\end{equation}
Lemma \ref{lemma-the-image-of-eta-v-under-Z}  implies, that 
$\beta^*\det(p_!\E_v)$ is $g$-invariant. 
Hence, $\det(W_0)$ is $g$-invariant as well.
We conclude, that 
the following isomorphism holds over $\Delta_0^{[1]}\cap \Delta_t^{[1]}$ 
\begin{equation}
\label{eq-solving-for-A-t}
(A_t\restricted{)}{\Delta_0^{[1]}\cap \Delta_t^{[1]}}
\ \  \cong \ \ 
(\delta_1^*\det(W_0)^*\restricted{)}{\Delta_0^{[1]}\cap\Delta_t^{[1]}}
\left(2\cdot 
\sum_{i=t+1}^\mu \Delta_0^{[1]}\cap\Delta_i^{[1]}\cap\Delta_t^{[1]}
\right).
\end{equation}

The homomorphism
$\Pic(\Delta_t^{[t+1]}) \rightarrow \Pic(\Delta_t^{[1]})
\rightarrow \Pic(\Delta_0^{[1]}\cap \Delta_t^{[1]})$
is injective. This is a relative analogue of the statement, that 
the restriction homomorphism 
\[
\Pic(G(t,H)\times G(t,H^*)) \ \ \longrightarrow \ \ \Pic(Flag(t,h-t,H))
\]
is injective, for a vector space $H$ of dimension $h>2t$
(see Lemma \ref{lemma-structure-of-correspondence-plus-2-case} part 
\ref{lemma-item-components-of-correspondence-intersect-as-flag-variety}).
The first equation in (\ref{eq-A-t-in-terms-of-E-v}) follows
from equations 
(\ref
{eq-comparison-between-det-W-0-and-det-of-pushforward-of-universal-sheaf})
and (\ref{eq-solving-for-A-t}).
The second equation follows from 
Lemma \ref{lemma-structure-of-correspondence-plus-2-case} part 
\ref{lemma-item-normal-bundle-of-Delta-t}. 

The equation, in terms of the incidence divisor, follows from 
lemma \ref{lemma-two-degeneracy-loci-on-Delta-t-are-identical} 
part \ref{lemma-item-pullback-of-incidence-divisor-has-multiplicity-t-i-along}.
\EndProof

\subsection{Proof of Theorem \ref{thm-reflection-sigma-satisfies-main-conj}}
\label{sec-proof-of-thm-reflection-sigma}

Lemma \ref{lemma-recovering-f} and equations 
(\ref{eq-relation-between-sigma-and-tau}),
(\ref{eq-chern-character-of-tau-E-v}),
and (\ref{eq-Chow-theoretic-formula-for-gamma-sigma}), 
reduce the proof of Theorem \ref{thm-reflection-sigma-satisfies-main-conj}
to the verification of the relation
\begin{equation}
\label{eq-push-forward-of-normalized-univ-sheaf-by-sigma}
([\Z_*\circ D_\M]\otimes \sigma)\left(
ch(p^*\eta_v)\cdot ch(\E_v)\cdot\pi_S^*\sqrt{td_{S}}
\right) \ = \
ch(p^*\eta_v)\cdot ch(\E_v)\cdot\pi_S^*\sqrt{td_{S}}
\end{equation}
over $\M(v)\times S$, where $\eta_v$ is the class in 
(\ref{eq-invariant-normalized-class-of-chern-character-of-universal-sheaf}).

{\bf Step I:}
Choose the universal family over $\M(v)\times S$ to be the twist by $\LB$ of
the ideal sheaf of the universal subscheme of $S^{[g-1]}$. 
We already know, that $\Z_*$ is a ring isomorphism
(\cite{markman-reflections} Theorem 1.1). 
Using lemma \ref{lemma-the-image-of-eta-v-under-Z}, we see that 
(\ref{eq-push-forward-of-normalized-univ-sheaf-by-sigma})
is equivalent to 
\begin{equation}
\label{eq-push-forward-of-univ-sheaf-by-Z-tensor-sigma}
([\Z_*\circ D_\M]\otimes\sigma)(ch(\E_v)\cdot\pi_S^*\sqrt{td_S}) \ = \
ch[\E_v\otimes p^*\det(p_*\E_v)^\vee]\pi_S^*\sqrt{td_S}.
\end{equation}

Using the relation $\sigma\circ D=-\tau$, 
given in (\ref{eq-relation-between-sigma-and-tau}), 
we see that
equality (\ref{eq-push-forward-of-univ-sheaf-by-Z-tensor-sigma})
is equivalent to 
\[
(\Z_*\otimes -\tau)\left(
ch(\E_v^\vee)\cdot\pi_S^*\sqrt{td_{S}}
\right) \ \ =  \ \ 
ch[\E_v\otimes p^*\det(p_*\E_v)^\vee]\cdot\pi_S^*\sqrt{td_{S}}.
\]
Combined with equality (\ref{eq-chern-character-of-tau-E-v}), 
the last equality 
translates to
\begin{equation}
\label{eq-push-forward-of-univ-sheaf-by-Z}
(\Z_*\otimes 1_S)\left[
ch(-\tau(\E_v^\vee))\pi_S^*\sqrt{td_S}
\right]
\ \  = \ \ 
ch[\E_v\otimes p^*\det(p_*\E_v)^\vee]\cdot\pi_S^*\sqrt{td_S}.
\end{equation}

{\bf Steps II,III:}
Next we translate (\ref{eq-push-forward-of-univ-sheaf-by-Z})
to a K-theoretic equivalence. 
Let $\Z^{[1]}$ be the iterated blow-up of $\Z$ defined in Section
\ref{sec-the-structure-of-the-correspondence-in-the-plus-2-case}.
Using the same argument, as in the proof of Theorem 
\ref{thm-class-of-correspondence-in-stratified-elementary-trans}, 
we get, that (\ref{eq-push-forward-of-univ-sheaf-by-Z}) follows from
\begin{equation}
\label{eq-K-theoretic-image-of-minus-tau-E-v-via-correspondence}
\tilde{\pi}_{12_!}\left(
\tilde{\pi}_{23}^![p^!p_!\E_v-\E_v]^\vee
\right)
\ \ \equiv \ \ 
\E_v\otimes p^*\det(p_*\E_v)^\vee,
\end{equation}
where $\tilde{\pi}_{ij}$ is the projection from
$\Z^{[1]}\times S$ onto the product of the 
$i$-th and $j$-th factors of $\M(v)\times S\times \M(v)$. 

Let $\tilde{\delta}_{t,12}$ and $\tilde{\delta}_{t,23}$ 
be the two morphisms from $\Delta^{[1]}_t\times S$ to 
$\M(v)\times S$. For any vector bundle $E$ on $\M(v)$, we have 
the equivalence, in the K-group of coherent sheaves on 
$\Z^{[1]}\times S$,
\[
\tilde{\pi}_{23}^*E^\vee \ \  \equiv \ \ 
\sum_{t=0}^\mu e_{t,!}\tilde{\delta}_{t,23}^*(E^\vee)\left(
-\sum_{i=0}^{t-1}\Delta^{[1]}_t\cap \Delta_i^{[1]}
\right). 
\]
The equality $\tilde{\pi}_{12}\circ e_t= \tilde{\delta}_{t,12}$ 
gives rise to the equivalence, in the K-group of $\M(v)\times S$, 
\[
\tilde{\pi}_{12_!}\left(
\tilde{\pi}_{23}^!E^\vee\right) \ \ \equiv \ \ 
\sum_{t=0}^\mu \tilde{\delta}_{t,12_!}
\left\{
\left(\tilde{\delta}_{t,23}^!E\right)^\vee
\left(
-\sum_{i=0}^{t-1}\Delta^{[1]}_t\cap \Delta_i^{[1]}
\right)
\right\}.
\]
We conclude, that the left hand side of 
(\ref{eq-K-theoretic-image-of-minus-tau-E-v-via-correspondence}) 
satisfies the equivalence 
\[
\tilde{\pi}_{12_!}\left(
\tilde{\pi}_{23}^![p^!p_!\E_v-\E_v]^\vee
\right)
\ \ \equiv \ \ 
\sum_{t=0}^\mu \tilde{\delta}_{t,12_!}
\left\{\left(\tilde{\delta}_{t,23}^!\left[p^!p_!\E_v-\E_v\right]^\vee\right)
\left(
-\sum_{i=0}^{t-1}\Delta^{[1]}_t\cap \Delta_i^{[1]}
\right)\right\}.
\]

{\bf Step IV:}
Corollary \ref{cor-lazarsfeld-reflection-of-universal-sheaf-over-Delta-t} 
implies, that the above sum is equivalent to the difference 
\begin{eqnarray}
\label{eq-term-involving-F-t-twisted-by-A-t}
\sum_{t=0}^\mu\tilde{\delta}_{t,12_!}
\left\{
\tilde{\delta}_{t,12}^!\F_t\otimes 
A_t\left(
-\sum_{\stackrel{i=0}{i\neq t}}^\mu\Delta_t^{[1]}\cap\Delta_i^{[1]}
\right)
\right\}
\\
\label{eq-term-involving-double-sum-of-U-i}
-\sum_{t=0}^\mu\tilde{\delta}_{t,12_!}
\left\{
\left(
\sum_{i=t}^\mu\tilde{e}_{t,i_!}
\left[
\tilde{\delta}_{t\cap i,3}^* U_{\restricted{i}{B^{[i]}\M(v)^{t\cap i}}}
\right]
\right)^\vee
\left(
-\sum_{j=0}^{t-1}\Delta_t^{[1]}\cap\Delta_j^{[1]}
\right)
\right\},
\end{eqnarray}
where $\tilde{\delta}_{t\cap i,3}:
[\Delta_t^{[1]}\cap \Delta_i^{[1]}]\times S \rightarrow 
B^{[i+1]}\M(v)^{t\cap i}$
is the composition of the inclusion in $\Delta_i^{[1]}\times S$, the
projection to $\Delta_i^{[i+1]}$, 
followed by the second projection from
$\Delta_i^{[i+1]}$ to $B^{[i+1]}\M(v)^{i}$. 
Corollary \ref{cor-lazarsfeld-reflection-of-universal-sheaf-over-Delta-t} 
is stated over $\Delta_t^{[t+1]}$. 
We used above lemma 
\ref{lemma-K-theoretic-pullback-of-a-sheaf-supported-on-a-divisor}, 
in order to restate it in terms of sheaves on $\Delta_t^{[1]}$. 
The first term (\ref{eq-term-involving-F-t-twisted-by-A-t})
decomposes further, using the short exact sequence 
(\ref{eq-the-tautological-extension-of-families-over-BN-locus}), to
the sum
\begin{eqnarray}
\label{eq-term-involving-E-v-twisted-by-A-t}
\sum_{t=0}^\mu\tilde{\delta}_{t,12_!}
\left\{
\tilde{\delta}_{t,12}^!
\E_{\restricted{v}{\M(v)^t}}
\left(
-2\cdot\sum_{i=t+1}^\mu \Delta_{t,i}^{[t+1]}
\right)
\otimes A_t(-I_t)
\right\}
\\
\label{eq-term-involving-U-t-twisted-by-A-t}
+ \sum_{t=0}^\mu\tilde{\delta}_{t,12_!}
\left\{
\tilde{\delta}_{t,1}^! U_t\otimes A_t
\left(
-\sum_{\stackrel{i=0}{i\neq t}}^\mu\Delta_t^{[1]}\cap\Delta_i^{[1]}
\right)
\right\}
\end{eqnarray}
We used lemma \ref{lemma-two-degeneracy-loci-on-Delta-t-are-identical} part
\ref{lemma-item-pushforward-of-ideal-sheaf-is-ideal-sheaf-of-incidence-div+2}
to write the $t$-th term in (\ref{eq-term-involving-E-v-twisted-by-A-t})
involving sheaves over $\Delta_t^{[t+1]}$.
Using lemma \ref{lemma-calculation-of-A-t}, the term
(\ref{eq-term-involving-E-v-twisted-by-A-t})
becomes
\[
\sum_{t=0}^\mu\tilde{\delta}_{t,12_!}
\left\{
\tilde{\delta}_{t,12}^!\E_{\restricted{v}{\M(v)^t}}\otimes 
\tilde{\delta}_1^*\det(p_*\E_v)^*
\left(
\sum_{i=t+1}^\mu (t-i)\Delta_{t,i}^{[t+1]}
\right)
\right\}.
\]

The morphism $\tilde{\delta}_{t,12}:\Delta_t^{[t+1]}\rightarrow \M(v)$
factors through $B^{[t+1]}\M(v)^t$ and the blow-down morphism
$\beta_t:B^{[t+1]}\M(v)^t\rightarrow \M(v)^t$.
Using the projection formula, once for each of the two morphisms, 
the term
(\ref{eq-term-involving-E-v-twisted-by-A-t})
further simplifies to
\[
\left(\E_v\otimes p^*\det(p_*\E_v)^*\right)\otimes
\left[
\sum_{t=0}^\mu 
\beta_{t_!}\StructureSheaf{B^{[t+1]}\M(v)^t}
\left(
-\sum_{i=t+1}^\mu (i-t) B^{[t+1]}\M(v)^{t\cap i}
\right)
\right].
\]
The above expression is equal to
\[
\left(\E_v\otimes p^*\det(p_*\E_v)^*\right)\otimes
\left[
\sum_{t=0}^\mu 
\beta_{t_!}\StructureSheaf{B^{[1]}\M(v)^t}
\left(
-\sum_{i=t+1}^\mu (i-t) B^{[1]}\M(v)^{t\cap i}
\right)
\right],
\]
using the projection formula and the fact that the morphism 
$B^{[1]}\M(v)^t\rightarrow B^{[t+1]}\M(v)^t$
is the composition of blow-ups with smooth centers. 
The latter expression is $\E_v\otimes p^*\det(p_*\E_v)^*$, as seen by
lemma \ref{lemma-pushforward-of-pull-back-of-ideal-sheaf-is-same-as} 
and the decomposition
\[
\StructureSheaf{B^{[1]}\M(v)} \ \ \equiv \ \ 
\sum_{t=0}^\mu \StructureSheaf{B^{[1]}\M(v)^t}
\left(
-\sum_{i=t+1}^\mu B^{[1]}\M(v)^{t\cap i}
\right).
\]

{\bf Step V:}
It remains to prove, that the sum of
(\ref{eq-term-involving-double-sum-of-U-i})
and
(\ref{eq-term-involving-U-t-twisted-by-A-t})
vanishes. 
We would like to reverse the order of summation in 
(\ref{eq-term-involving-double-sum-of-U-i}) using the identity 
$\tilde{\delta}_{i,12}\circ \tilde{e}_{i,t}=
\tilde{\delta}_{t,12}\circ \tilde{e}_{t,i}$.
We will further need the identity
\[
\tilde{\delta}_{t,12_!}
\left\{
\left(
\tilde{e}_{t,i_!}
\left[
\tilde{\delta}_{t\cap i,3}^* U_{\restricted{i}{B^{[i]}\M(v)^{t\cap i}}}
\right]
\right)^\vee
\right\}
\ \ \equiv \ \ 
\tilde{\delta}_{i,12_!}
\left\{
\left(
\tilde{e}_{i,t_!}
\left[
\tilde{\delta}_{t\cap i,3}^* U_{\restricted{i}{B^{[i]}\M(v)^{t\cap i}}}
\right]
\right)^\vee
(\Delta_i^{[1]}-\Delta_t^{[1]})
\right\},
\]
where the dualization on the left is carried over $\Delta_t^{[1]}$
and the one on the right is carried over $\Delta_i^{[1]}$. 
The above identity follows from lemma
\ref{lemma-dualizing-the-pushforward-of-a-vector-bundle-from-a-subvariety},
and the identities 
\begin{eqnarray*}
N_{\Delta_t^{[1]}\cap \Delta_i^{[1]}/\Delta_t^{[1]}} & \cong &
\StructureSheaf{\Delta_t^{[1]}\cap \Delta_i^{[1]}}(\Delta_i^{[1]}),
\\
N_{\Delta_t^{[1]}\cap \Delta_i^{[1]}/\Delta_i^{[1]}} & \cong &
\StructureSheaf{\Delta_t^{[1]}\cap \Delta_i^{[1]}}(\Delta_t^{[1]}).
\end{eqnarray*}

Reversing the order of summation in 
(\ref{eq-term-involving-double-sum-of-U-i}), 
the term (\ref{eq-term-involving-double-sum-of-U-i}) becomes 
\[
-\sum_{i=1}^\mu\tilde{\delta}_{i,12_!}
\left\{
\sum_{t=0}^i
\left(
\tilde{e}_{i,t_!}
\left[
\tilde{\delta}_{t\cap i,3}^* U_{\restricted{i}{B^{[i]}\M(v)^{t\cap i}}}
\right]
\right)^\vee
\left(
\Delta_i^{[1]}-\Delta_t^{[1]}
-\sum_{j=0}^{t-1}\Delta_j^{[1]}
\right)
\right\}
\]
(because $U_0=0$, by part
\ref{property-constant-rank}
of Proposition \ref{prop-decomposition-of-complex-over-interated-blow-up}).
Using lemma 
\ref{lemma-dualizing-the-pushforward-of-a-vector-bundle-from-a-subvariety}
once more, for the terms where $i>t$, 
the above expression can be written in the form 
\[
-\sum_{i=1}^\mu\tilde{\delta}_{i,12_!}
\tilde{\delta}_{3}^* U^*_i 
\left(
-\sum_{j=0}^{i-1}\Delta_j^{[1]}
\right)
\ \ + \ \ 
\sum_{i=1}^\mu\tilde{\delta}_{i,12_!}
\left\{
\sum_{t=0}^{i-1}
\left(
\tilde{e}_{i,t_!}
\left[
\tilde{\delta}_{t\cap i,3}^* U^*_{\restricted{i}{B^{[i]}\M(v)^{t\cap i}}}
\right]
\right)
\left(
\Delta_i^{[1]}
-\sum_{j=0}^{t-1}\Delta_j^{[1]}
\right)
\right\}.
\]
Add to the right term above and subtract from the left term the sum
\[
\sum_{i=1}^\mu\tilde{\delta}_{i,12_!}
\tilde{\delta}_{3}^* U^*_i 
\left(
\Delta_i^{[1]}
-\sum_{j=0}^{i-1}\Delta_j^{[1]}
\right)
\] 
to get, that (\ref{eq-term-involving-double-sum-of-U-i}) simplifies to
\begin{equation}
\label{eq-result-of-reverse-of-order-of-summation}
\sum_{i=1}^\mu\tilde{\delta}_{i,12_!}
\left\{
-\tilde{\delta}_{3}^* U^*_i 
\left(
-\sum_{j=0}^{i-1}\Delta_j^{[1]}
\right)
- \tilde{\delta}_{3}^* U^*_i 
\left(
\Delta_i^{[1]}
-\sum_{j=0}^{i-1}\Delta_j^{[1]}
\right)
+ \left(
\tilde{\delta}_{3}^* U^*_i
\right)(\Delta_i^{[1]})
\right\},
\end{equation}
where $\tilde{\delta}_{3}:\Delta_i^{[1]}\times S\rightarrow B^{[i+1]}\M(v)^i$
is the composition of the blow-down morphism,
with the second projection from $\Delta_i^{[i+1]}$ to $B^{[i+1]}\M(v)^i$.

Next, we simplify the expression (\ref{eq-term-involving-U-t-twisted-by-A-t}).
Lemma \ref{lemma-relationship-between-two-pullcacks-of-U_t-to-Delta-t}
and lemma
\ref{lemma-structure-of-correspondence-plus-2-case}
part \ref{lemma-item-normal-bundle-of-Delta-t} imply, that 
(\ref{eq-term-involving-U-t-twisted-by-A-t})
is equivalent to
\[
\sum_{t=0}^\mu\tilde{\delta}_{t,12_!}
\left\{
\tilde{\delta}_3^*U_t^* 
\left(\Delta_t^{[1]}-\sum_{i=0}^{t-1}\Delta_i^{[1]}\right)
-
\StructureSheaf{\Delta_t^{[1]}}
\left(
\sum_{j=t}^\mu\Delta_j^{[1]}
\right)\otimes
\sum_{k=1}^{t-1}
\left(
e_{t,k_!}\tilde{\delta}_{t\cap k,3}^*U_k
\right)^\vee
\right\}.
\]
Lemma 
\ref{lemma-dualizing-the-pushforward-of-a-vector-bundle-from-a-subvariety}
translates it to
\[
\sum_{t=0}^\mu\tilde{\delta}_{t,12_!}
\left\{
\tilde{\delta}_3^*U_t^* 
\left(
\Delta_t^{[1]}-\sum_{i=0}^{t-1}\Delta_i^{[1]}
\right)
+
\StructureSheaf{\Delta_t^{[1]}}
\left(
\sum_{j=t}^\mu\Delta_j^{[1]}
\right)\otimes
\sum_{k=1}^{t-1}
\left(
e_{t,k_!}\tilde{\delta}_{t\cap k,3}^*U^*_k(\Delta_k^{[1]})
\right)
\right\}.
\]
Reversing the order of summation, we get
\begin{equation}
\label{eq-result-of-reversing-summation-in-term-involving-U-t-twisted-by-A-t}
\sum_{k=1}^\mu
\tilde{\delta}_{k,12_!}
\left\{
\tilde{\delta}_3^*U_k^*\left(
\Delta_k^{[1]}-\sum_{i=0}^{k-1}\Delta_i^{[1]}
\right)
+\sum_{t=k+1}^\mu
\left(
e_{k,t_!}\tilde{\delta}_{t\cap k,3}^*U_k^*
\left(
\Delta_k^{[1]}+\sum_{j=t}^\mu\Delta_j^{[1]}
\right)
\right)
\right\}.
\end{equation}
The $k$-th inner sum, in the expression above, is the restriction
of $U_k^*\left(\sum_{i=k}^\mu\Delta_i^{[1]}\right)$ to the divisor 
$\cup_{t=k+1}^\mu[\Delta_k^{[1]}\cap \Delta_t^{[1]}].$
This restriction is the difference
\[
U_k^*\left(
\sum_{i=k}^\mu\Delta_i^{[1]}
\right)
-
U_k^*(\Delta_k^{[1]}).
\]
Using the rational equivalence 
in lemma \ref{lemma-structure-of-correspondence-plus-2-case} part 
\ref{lemma-item-normal-bundle-of-Delta-t}, 
the expression 
(\ref{eq-result-of-reversing-summation-in-term-involving-U-t-twisted-by-A-t})
becomes
\[
\sum_{k=1}^\mu
\tilde{\delta}_{k,12_!}
\left\{
\tilde{\delta}_3^*U_k^*\left(
\Delta_k^{[1]}-\sum_{i=0}^{k-1}\Delta_i^{[1]}
\right)
+
\tilde{\delta}_3^*U_k^*\left(
-\sum_{i=0}^{k-1}\Delta_k^{[1]}
\right)
-
\tilde{\delta}_3^*U_k^*(\Delta_k^{[1]})
\right\}.
\]
This is precisely the negative of 
(\ref{eq-result-of-reverse-of-order-of-summation}). 
This completes the proof of Theorem 
\ref{thm-reflection-sigma-satisfies-main-conj} 
\EndProof

The following corollary is an easy consequence of Theorem 
\ref{thm-reflection-sigma-satisfies-main-conj}. 

\begin{cor}
\label{cor-action-of-Z-on-individual-strata-plus-2-case}
The correspondence $\Z$, given in
(\ref{eq-correspondence-inducing-involution-of-hilbert-scheme}), 
acts via multiplication by $(-1)^t$ on the class of $\M(v)^t$ 
in the Chow group of $\M(v)$. 
\end{cor}

\noindent
{\bf Proof:} We only sketch the proof, as it is similar to that of Corollary 
\ref{cor-action-of-Z-on-individual-strata}. 
Equality (\ref{eq-K-theoretic-image-of-minus-tau-E-v-via-correspondence}) 
and Proposition \ref{prop-rational-singularities-of-correspondence}
imply the K-theoretic equivalence
\[
\Z_!([p_!\E_v]^\vee) \ \ \ \equiv \ \ \ p_!\E'_v,
\]
where $\E'_v$ is another universal sheaf (compare with equation 
(\ref{eq-correspondence-multiplies-K-theoretic-image-of-E-v-by-minus-1})). 
Equivalence (b) below follows, where
$\Delta_t^{(s)}$ is the determinant, introduced in
the proof of Corollary \ref{cor-action-of-Z-on-individual-strata}. 

\noindent
$
\Z_*[\M(v)^t] 
\ \stackrel{(a)}{\equiv} \
\Z_*\Delta^{(t+\chi(v))}_t\left(c[-p_!\E_v\right]) 
\ \stackrel{(b)}{\equiv} \ 
\Delta^{(t+\chi(v))}_{t}\left(c[-p_!\E'_v\right]^\vee) 
\ \stackrel{(c)}{\equiv} \ 
\Delta_{t+\chi(v)}^{(t)}\left(c[p_!\E'_v\right]) 
\ \stackrel{(d)}{\equiv} \ 
(-1)^{t(t+\chi(v))}\Delta^{(t+\chi(v))}_{t}
\left(c[-p_!\E'_v\right]) 
\ \stackrel{(e)}{\equiv} \ 
(-1)^{t(t+\chi(v))}\M(v)^t.
$

Equivalences (a) and (e) follow from the Porteous formula and proposition
\ref{prop-decomposition-of-complex-over-interated-blow-up}
part \ref{property-equivalent-to-push-forward-of-universal-sheaf}.
Equivalence (c) follows from consequence (2) of lemma 14.5.1
in \cite{fulton}. Equivalence (d) follows from example
14.4.9 in \cite{fulton}. 
In our case, $\chi(v)=2$ and $(-1)^{t(t+\chi(v))}=(-1)^t$. 
\EndProof

\begin{new-lemma}
\label{lemma-dualizing-the-pushforward-of-a-vector-bundle-from-a-subvariety}
Let $X$ be a local complete intersection subvariety of 
codimension $n$ in a smooth variety $Y$, 
$e:X\hookrightarrow Y$ the closed imersion, and $U$ 
a vector bundle over $X$.
Then the following equivalence holds in the K-group of $Y$
\[
[e_!(U)]^\vee \ \ \equiv \ \ 
-(1)^n\cdot e_!\left[U^*\otimes\det\left(
N_{X/Y}
\right)\right].
\]
\end{new-lemma}

\noindent
{\bf Proof:}
We have the equivalence
\[
[e_!(U)]^\vee \ \ \equiv \ \ 
\sum_{i=0}^n (-1)^i
\SheafExt^i_Y(e_*U,\StructureSheaf{Y}).
\]
All extension sheaves vanish, for $0 \leq i\leq n-1$, and
$\SheafExt^n_Y(e_*U,\StructureSheaf{Y})\cong 
e_*\left[U^*\otimes\det\left(
N_{X/Y}
\right)\right]$,
by the Local Duality Theorem.
\EndProof


\begin{new-lemma}
\label{lemma-pushforward-of-pull-back-of-ideal-sheaf-is-same-as}
The following equivalence holds in the K-group of $\M(v)$
\[
e_!\beta_{t_!}\StructureSheaf{B^{[t+1]}\M(v)^t}
\left(
-\sum_{i=t+1}^\mu (i-t) B^{[t+1]}\M(v)^{t\cap i}
\right)\ \ \equiv \ \ 
e_!\beta_{t_!}\StructureSheaf{B^{[t+1]}\M(v)^t}
\left(
-\sum_{i=t+1}^\mu B^{[t+1]}\M(v)^{t\cap i}
\right),
\]
where $\beta_t: B^{[t+1]}\M(v)^t \rightarrow \M(v)^t$ is the natural morphism
and $e:\M(v)^t\hookrightarrow \M(v)$ is the closed imersion.
\end{new-lemma}

\noindent
{\bf Proof:}
The proof consists of two steps. In step A, we reduce the proof to the case
$t=0$. In step B, we reduce the proof to  the case covered by
claim \ref{claim-push-forward-of-ideal-sheaf-is-ideal-sheaf}, 
in which the determinantal variety is a divisor. 

{\bf Step A:} 
Consider the following diagram
\[
\begin{array}{ccccccc}
& & B^{[t+1]}\M(v)^t & \RightArrowOf{\beta_t} & \M(v)^t 
\\
& f \swarrow \hspace{1ex} & & \hspace{2ex} \searrow \alpha_t & & \nwarrow 
\\
B^{[1]}\M(v+\vec{t}) & & & & Flag(t,\chi,v+\vec{t}) & \rightarrow &
G(t,v+\vec{t})
\\
 & \psi \searrow \hspace{1ex} & &  \hspace{1ex} \swarrow \phi & 
\\
& & G(\chi,v+\vec{t})
\\
& & \hspace{1ex} \downarrow g
\\
& & \M(v+\vec{t})
\end{array}
\]
Above, $\chi:=\chi(v+\vec{t}) = \chi(v)+2t$ is the Euler characteristic 
of sheaves in the moduli space $\M(v+\vec{t})$. 
We denote by $G(k,v+\vec{t})$, for $k=t,\chi$, the moduli spaces of 
$k$-dimensional coherent systems. These moduli spaces parametrize pairs
$(F,W)$ of a stable sheaf $F$ on the K3 surface $S$, with Mukai vector 
$v+\vec{t}$, and a $k$-dimensional subspace $W$ of $H^0(S,F)$.
They were constructed by Le Potier \cite{le-potier} (with a more general 
definition of stability, which we do not need). 
We denote by $\phi:Flag(t,\chi,v+\vec{t})\rightarrow G(\chi,v+\vec{t})$ the 
Grassmannian bundle parametrizing 
triples $(F,W_t\subset W_\chi)$, where $(F,W_\chi)$ is a coherent system 
in $G(\chi,v+\vec{t})$, and $W_t$ is a $t$-dimensional subspace of $W_\chi$. 
The existence of the morphisms $f$, $\psi$, and $\alpha_t$, 
was proven in \cite{markman-reflections} Proposition 3.18
(see also parts \ref{prop-item-U-and-W-are-direct-images} and 
\ref{prop-item-morphism-is-a-classifying-morphism-of-F-t}
of Proposition 
\ref{prop-decomposition-of-complex-over-interated-blow-up} above).  
We let $g$ be the forgetful morphism.

Set 
\begin{eqnarray*}
F_1 & := & \StructureSheaf{B^{[1]}\M(v+\vec{t})}
\left(
-\sum_{i=1}^{\mu-t} B^{[1]}\M(v+\vec{t})^i
\right),
\\
F_2  & := & \StructureSheaf{B^{[1]}\M(v+\vec{t})}
\left(
-\sum_{i=1}^{\mu-t} i\cdot B^{[1]}\M(v+\vec{t})^i
\right). 
\end{eqnarray*}
The statement of the claim translates to the equality
\[
e_!\beta_{t_!}f^*F_1 \ \ \equiv \ \ e_!\beta_{t_!}f^*F_2, 
\]
because the exceptional divisors in $B^{[t+1]}\M(v)^t$ are the pullback of 
those in $B^{[1]}\M(v+\vec{t})$ (see \cite{markman-reflections} 
proposition 3.18 or equation 
(\ref{eq-conditioned-equality-of-catier-divisor}) above).
Since $\beta_t$ factors through $\alpha_t$, it suffices to prove 
the equality
\[
\alpha_{t_!}f^*F_1 \ \ \equiv \ \ \alpha_{t_!}f^*F_2. 
\]

It was proven in 
 \cite{markman-reflections} Proposition 3.18, that 
$f:B^{[t+1]}\M(v)^t\rightarrow B^{[1]}\M(v+\vec{t})$ is
a Grassmannian bundle, which is the pullback of 
$\phi : Flag(t,\chi,v+\vec{t})\rightarrow G(\chi,v+\vec{t})$
(see also Proposition
\ref{prop-decomposition-of-complex-over-interated-blow-up} part 
\ref{prop-item-morphism-is-a-classifying-morphism-of-F-t}
and Lemma \ref{lemma-W-v-plus-t-is-a-direct-image} above). 
Since $\phi$ is flat, and cohomology commutes with flat base change, 
then $\alpha_{t_!}f^!(F)$ is equivalent to 
$\phi^!\psi_!(F)$, for any sheaf $F$ on $B^{[1]}\M(v+\vec{t})$. 
The claim reduces to the equivalence
\begin{equation}
\label{eq-psi-pushes-forward-F-1-and-F-2-to-the-same-class}
\psi_!(F_1) \ \ \equiv \ \ \psi_!(F_2). 
\end{equation}

{\bf Step B:} 
We will show next, that $B^{[1]}\M(v+\vec{t})$ is the iterated blow-up of 
both $\M(v+\vec{t})$ and $G(\chi,v+\vec{t})$, 
with respect to closely related determinantal stratifications. 
Both staratifications are ``transversal''. 
Moreover, the determinantal variety $G(\chi,v+\vec{t})^1$ is a divisor
in $G(\chi,v+\vec{t})$. Hence, we would obtain the desired
reduction of the lemma to the case covered in
claim \ref{claim-push-forward-of-ideal-sheaf-is-ideal-sheaf}.

Recall, that the Brill-Noether stratification of 
$\M(v+\vec{t})$ is the determinantal stratification of 
a homomorphism $\rho : V_0\rightarrow V_1$
between vector bundles, with $\rank(V_0)=\rank(V_1)+\chi$
(see Proposition \ref{prop-decomposition-of-complex-over-interated-blow-up} 
above or section 3.5 in \cite{markman-reflections}). 
The homomorphism $\rho$ exists globally, only 
if there is a universal sheaf over $\M(v+\vec{t})$, 
but it always exists locally. 
Over $G(\chi,v+\vec{t})$, there is a tautological 
rank $\chi$ subbundle $W_\chi$ of $g^*V_0$, which is contained in the kernel
of $g^*\rho$. We get the induced homomorphism
\[
\bar{\rho} \ : \ \bar{V}_0 \ \longrightarrow \ g^*V_1,
\]
where $\bar{V}_0$ is the quotient bundle $g^*V_0/W_\chi$. 

The vector bundles $\bar{V}_0$ and $g^*V_1$ have the same rank. 
Hence, the determinantal variety of $\bar{\rho}$ is a divisor. 
The section $\rho$ is transversal to the 
determinantal variety in the total space of 
$\Hom(V_0,V_1)$ (see \cite{markman-reflections} Proposition 3.18). Hence, 
$G(\chi,v+\vec{t})$ is smooth 
(see \cite{a-c-g-h} Chapter IV section 4). 
We claim, that the image of $G(\chi,v+\vec{t})$, via $\bar{\rho}$, 
is transversal to the determinantal variety in the total space 
of $\Hom(\bar{V}_0,g^*V_1)$. 
Furthermore, the morphism $\psi$
is the iterated blow-up, of the determinantal stratification 
induced by $\bar{\rho}$ on $G(\chi,v+\vec{t})$. 
The equality 
(\ref{eq-psi-pushes-forward-F-1-and-F-2-to-the-same-class}) 
would then follow, by the argument used in the proof of
Claim \ref{claim-push-forward-of-ideal-sheaf-is-ideal-sheaf}. 

It remains to prove the transversality of $\bar{\rho}$. 
Let $\pi: G(\chi,V_0) \rightarrow \M(v+\vec{t})$ be the Grassmannian bundle
of $\chi$-dimensional subspaces in the fibers of $V_0$, 
and let $q$ be the tautological quotient bundle of $\pi^*V_0$ over 
$G(\chi,V_0)$. 
Then $G(\chi,v+\vec{t})$ is a subvariety of $G(\chi,V_0)$.
The bundle $\Hom(q,\pi^*V_1)$ over $G(\chi,V_0)$ 
has a section $s$, induced by the pullback $\pi^*\rho$ of $\rho$.
The section $\bar{\rho}$ is the restriction of $s$ to $G(\chi,v+\vec{t})$.
The construction fits in the following {\em cartesian} diagram.
\[
\begin{array}{ccc}
G(\chi,v+\vec{t}) & \LongRightArrowOf{\restricted{s}{}} & \Hom(q,\pi^*V_1)
\\
g \downarrow \hspace{1ex} & & \downarrow
\\
\M(v+\vec{t}) & \LongRightArrowOf{\rho} & \Hom(V_0,V_1). 
\end{array}
\]
The cartesian nature of the diagram implies, that $\rho$ is
transversal, if and only if $\bar{\rho}=\restricted{s}{}$ is
transversal. 
This is seen as follows. Choose a trivialization of $V_0$ and $V_1$
in a neighborhood $U$ of a point $x$ in $\M(v+\vec{t})$. 
The section $\rho$ and the trivialization induce a morphism 
$\rho_x$ from $U$ to the fiber $\Hom(V_0,V_1)_x$. 
Similarly, the section $s$ induces a morphism $s_x$,
from the open subset $g^{-1}(U)\subset G(\chi,v+\vec{t})$ to 
$\Hom(q_x,V_{1,x})$, where $q_x$ is the tautological quotient bundle, over
the grassmanian $G(\chi,V_{0,x})$. 
The cartesian nature of the diagram implies, that $\rho_x$ is submersive,
if and only if $s_x$ is submersive.
\EndProof

\section{Monodromy of self-dual symplectic singularities}
\label{sec-monodromy-of-A1-singularities}

We exhibit the involutions $\gamma_{\tau_{v_0}}$ in Theorem 
\ref{thm-class-of-correspondence-in-stratified-elementary-trans} 
and $D_{\M(v)}\circ \gamma_{\sigma_{u_0}}$ in Theorem 
\ref{thm-reflection-sigma-satisfies-main-conj} 
as monodromy operators.
Let $\rho$ be either the reflection $\tau_{v_0}$ in Theorem 
\ref{thm-class-of-correspondence-in-stratified-elementary-trans} 
or the involution $-\sigma_{u_0}$, where $\sigma_{u_0}$ is the
reflection in Theorem 
\ref{thm-reflection-sigma-satisfies-main-conj}. 
Let $Q\subset \PP{}H^2(\M(v),\ComplexNumbers)$ be the quadric
cut out by the Beauville pairing and $\Omega\subset Q$ the open analytic 
period domain 
\begin{equation}
\label{eq-period-domain-Omega}
\{\omega \ | \ (\omega,\bar{\omega}) > 0\}/\ComplexNumbers^\times.
\end{equation} 
Regard $H^{2,0}(\M(v))$ as a point $\ell$ in $\Omega$.
The isometry $\rho$ acts on $v^\perp$ and,
via the isomorphism 
$v^\perp\otimes_\Integers\ComplexNumbers\cong H^{2}(\M(v),\ComplexNumbers)$, 
also on $Q$ (Theorem \ref{thm-irreducibility}). 
The point $\ell$ is $\rho$-invariant because 
the latter is a Hodge isometry. 
Let $W\subset H^1(T\M(v))$ be a simply connected open neighborhood of $0$, 
over which we have a universal (Kuranishi) family. 
Let $U\subset Q$ be the image of $W$ via the period map and 
\begin{equation}
\label{eq-universal-family-of-hyperkahler-varieties-over-moduli}
\X  \ \rightarrow \ U
\end{equation}
the (push forward of the) universal  family 
of hyperk\"{a}hler deformations of $\M(v)$. We assume, that $U$ is
symmetric with respect to $\rho$ and denote the fixed locus 
by $U^{\rho}$. Its tangent space is the $(+1)$-eigenspace 
$(T_\ell Q)^+$. 
When $\rho=\tau_{v_0}$, the fixed locus $U^{\rho}$ has codimension $1$, 
while  
the $(-1)$-eigenspace $(T_\ell Q)^-$ is 
the line $\Hom(H^{2,0}(\M(v)),\mbox{span}\{v_0\})$.
When $\rho=-\sigma_{u_0}$, the fixed locus $U^{\rho}$ is one-dimensional.

\begin{new-lemma}
\label{lemma-correspondence-is-a-monodromy-operator}
\begin{enumerate}
\item
\label{lemma-item-the-involution-lifts-to-the-kuranishi-family}
Each of the correspondences 
$\Z \subset \M(v)\times \M(v)$, of Theorems
\ref{thm-class-of-correspondence-in-stratified-elementary-trans} 
and 
\ref{thm-reflection-sigma-satisfies-main-conj}, 
deforms to an isomorphism of the two families over 
$U\setminus U^{\rho}$
\begin{equation}
\label{eq-family-is-tau-invariant}
\restricted{\X}{[U\setminus U^{\rho}]} \ \IsomRightArrow \ 
\rho^*\restricted{\X}{[U\setminus U^{\rho}]}.
\end{equation}
Consequently, we get a universal family 
\begin{equation}
\label{eq-universal-quotient-family}
\X/\rho \ \longrightarrow \ [U\setminus U^{\rho}]/\rho. 
\end{equation}
\item
\label{lemma-item-gamma-tau-is-the-monodromy}
Moreover, the monodromy group, of the family over
$[U\setminus U^{\rho}]/\rho$, 
is generated by the involution $\Z_*$ 
of $H^*(\M(v),\Integers)$ induced by $\Z$.
\end{enumerate}
\end{new-lemma}

\noindent
{\bf Proof:}
\ref{lemma-item-the-involution-lifts-to-the-kuranishi-family})
Let $C$ be a smooth analytic curve in $U$ through $\ell$, 
which is transversal to the tangent hyperplane
$A:=\Hom[H^{2,0}(\M(v)),\{v_0,v\}^\perp\cap H^{1,1}(\M(v))]$. 
(When $\rho=\tau_{v_0}$, then $A$ is the tangent space
of the fixed locus $U^{\tau_{v_0}}$ at $\ell$).
We denote the restriction of the universal family to $C$ by 
\begin{equation}
\label{eq-restricted-Kuranishi-family-X}
\pi' \ : \ \X' \ \rightarrow \ C.
\end{equation}
We show next, that this family satisfies Condition
\ref{cond-extension-class-is-non-trivial}.
Let $e:\PP^n\hookrightarrow \M(v)$ be a fiber of
$\M(v)^1\setminus \M(v)^{2}\rightarrow \M(v+\vec{1})$.
Consider the composition
\[
H^1(\M(v),T\M(v)) \ \LongIsomRightArrowOf{dp} \ 
\Hom(H^{2,0}(\M(v)),H^{1,1}(\M(v)) \ \LongRightArrowOf{e^*} \ 
\Hom(H^{2,0}(\M(v)),H^{1,1}(\PP^n)),
\]
where $dp$ is the differential of the period map, as identified by Griffiths, 
and $e^*$ is the pullback. The middle vector space is identified with 
$T_\ell Q$. 
Lemma \ref{lemma-degree-of-determinantal-line-bundles-on-fibers-of-M-1} 
identifies the kernel of $e^*$ with the hyperplane $A$.
Thus, Condition \ref{cond-extension-class-is-non-trivial}, for 
the family (\ref{eq-restricted-Kuranishi-family-X}), follows from
the transversality of $C$ and $A$. 
Let 
\[
\pi'' \ : \ \X'' \ \rightarrow \ C
\]
be the birational stratified transform of $\X'$ 
and let
\begin{equation}
\label{eq-g-is-the-identity-away-from-special-fiber}
g \ : \ \restricted{\X'}{C\setminus\{\ell\}} \ \LongIsomRightArrow \ 
\restricted{\X''}{C\setminus\{\ell\}}
\end{equation}
be the corresponding isomorphism of families
(see \cite{markman-reflections} Theorem 2.4 or equation (\ref{eq-Y}) above). 
Denote the special fibers by $X'_\ell$ and $X''_\ell$. They are
both isomorphic to $\M(v)$. We get the two classifying 
maps $\kappa'$ and $\kappa'': C\rightarrow W$ of $\X'\rightarrow C$ 
and $\X''\rightarrow C$. 
We know already, that the correspondence 
$\Z \subset X'_\ell\times X''_\ell$ 
deforms as the graph of the isomorphism $g_t$ 
of the the fibers $X'_t$ and  $X''_t$, $t\in C\setminus\{\ell\}$ 
(see \cite{markman-reflections} Theorem 2.4 or Proposition
\ref{prop-rational-singularities-of-correspondence} 
part \ref{prop-item-Z-is-a-fiber} above). 
On the other hand, the correspondence $\Z$ induces the involution $\rho$
on $H^2(\M(v),\ComplexNumbers)$. Thus, we get the following relation between
the compositions, of each of the classifying maps $\kappa'$ and 
$\kappa''$, with the period map $p: W\rightarrow U$ of the universal family: 
\[
p\circ \kappa' \ = \ \rho\circ p\circ \kappa''.
\]
It follows, that $\X''$ is isomorphic to the restriction of the 
universal family
(\ref{eq-universal-family-of-hyperkahler-varieties-over-moduli})
to the image $\rho(C)$ of $C$ under the reflection $\rho$. 

It remains to prove, that the construction of the isomorphism 
$X_t\cong X_{\rho(t)}$, $t\in C$, is independent of the choice of the 
curve $C$, and varies analytically with $t$. We perform 
the above construction over $U\setminus U^{\rho}$ as follows. 
Let $\beta:\widetilde{U}\rightarrow U$ be the blow-up of $U$ at $\ell$ and
$E$ the exceptional divisor. 
Then $\PP{A}$ is a hyperplane of $E$. 
(When $\rho=\tau_{v_0}$, then $\PP{A}$ is the intersection of $E$ with 
the proper transform of $U^{\tau_{v_0}}$). 
The differential $d\beta$ embeds the normal bundle
$N_{E/\widetilde{U}}$ in $\beta^*TU$ and $\PP{A}$ is the locus in $E$,
where $d\beta(N_{E/\widetilde{U}})$ is contained in $A$. 
The stratification of $X_\ell$, by Grassmannian bundles, 
pulls back to a stratification of $X_\ell\times [E\setminus \PP{A}]$ 
by Grassmannian bundles. The family 
$\restricted{\X}{\widetilde{U}\setminus \PP{A}}$ has the following property.
The normal bundle, to each stratum in $X_\ell\times [E\setminus \PP{A}]$, 
restricts to each Grassmannian fiber 
as a non-trivial extension
\[
0\rightarrow T^*_{G(t,k)}\rightarrow N \rightarrow 
\StructureSheaf{G(t,k)}
\rightarrow 0.
\] 
This is precisely the condition needed, in order to carry out the stratified
elementary transformation of the family 
$\restricted{\X}{\widetilde{U}\setminus \PP{A}}$, using a 
straightforward generalization of Theorem 2.4 in \cite{markman-reflections}
(compare also with Condition \ref{cond-extension-class-is-non-trivial} 
above). 

\smallskip
\ref{lemma-item-gamma-tau-is-the-monodromy})
Choose a smooth analytic curve $C$ in $U$, which is 
tangent to $(T_\ell Q)^-$ at $\ell$ 
and is invariant with respect to $\rho$. 
Assume, that $C$ is simply connected, possibly after replacing it by an open 
subset containing $\ell$. 
Denote by $\iota:C\rightarrow C$ the involution induced by $\rho$.
Let $g : \restricted{\X'}{C\setminus\{\ell\}} \rightarrow 
(\iota^*\X'\restricted{)}{C\setminus\{\ell\}}$
be the restriction of the isomorphism 
(\ref{eq-family-is-tau-invariant}). 
We get an induced isomorphism 
$g_*:(R^i_{\pi'_*}\Integers\restricted{)}{C\setminus\{\ell\}}\rightarrow
\iota^*(R^i_{\pi''_*}\Integers\restricted{)}{C\setminus\{\ell\}}$
of the local systems of relative integral cohomologies.
Both local systems are trivial, since the family $\X'$ is defined on the 
whole of $C$. Thus, $g_*$ extends to an isomorphism 
\begin{equation}
\label{eq-extended-g-star}
R^i_{\pi'_*}\Integers \ \ \LongIsomRightArrow \ \ 
\iota^*(R^i_{\pi''_*}\Integers)
\end{equation}
of the local systems over the whole of $C$. 
By the construction of the correspondence $\Z$, the homomorphism 
$\Z_*:H^*(\M(v),\Integers) \rightarrow H^*(\M(v),\Integers)$ 
is equal to the limit ${\displaystyle \lim_{t\rightarrow \ell}g_{t_*}}$.
Thus, upon trivialization of both local systems in 
(\ref{eq-extended-g-star}), modeled after the fiber over $\ell$,
the isomorphism $g_*$ in (\ref{eq-extended-g-star}) becomes
\[
\iota \times \Z_* \ : \ C\times H^*(\M(v),\Integers) 
\ \ \longrightarrow \ \ C\times H^*(\M(v),\Integers).
\]

Let $\bar{C} := C/\iota$ be the quotient and $\bar{\ell}\in \bar{C}$ the 
branch point corresponding to $\ell$. 
We conclude, that the restriction to $\bar{C} \setminus \bar{\ell}$, 
of the local system of relative cohomologies of 
(\ref{eq-universal-quotient-family}), 
is isomorphic to the quotient, of the trivial local system 
$(C\setminus \{\ell\})\times H^*(\M(v),\Integers)$, 
by the action of $\iota \times \Z_*$. Consequently, 
$\Z_*$ is the monodromy operator of the former local system. 
Part \ref{lemma-item-gamma-tau-is-the-monodromy} follows from
the surjectivity of $\pi_1(\bar{C} \setminus \bar{\ell})\rightarrow
\pi_1([U\setminus U^{\rho}]/\rho)$.
\EndProof


\end{document}